\documentclass[a4paper,times,3p]{elsarticle}
\usepackage[utf8]{inputenc}
\usepackage[T1]{fontenc}

\usepackage{microtype}

\usepackage{booktabs}

\usepackage{amsmath,amsfonts,amsthm,enumerate,graphicx,epstopdf,url,etoolbox,algorithm,listings}
\usepackage[noend]{algpseudocode}
\usepackage{subcaption}
\usepackage{caption}
\usepackage{diagbox}
\usepackage{graphbox}
\usepackage{enumitem}
\usepackage{dirtytalk}
\usepackage{tikz}
\usepackage[export]{adjustbox}
\usetikzlibrary{shapes, arrows}
\usepackage{tcolorbox}
\usepackage{mathtools}
\usepackage{wrapfig}
\usepackage{picins}

\definecolor{cmykcyan}{cmyk}{1,0,0,0}
\definecolor{cmykred}{cmyk}{0,1,1,0}
\definecolor{cmykblack}{cmyk}{0,0,0,1}

\newtheorem*{remark}{Remark}
\newtheorem{theorem}{Theorem}

\newcommand{\br}{{\mathbf{r}}}

\newcommand{\bx}{{\mathbf{x}}}
\newcommand{\bxi}{{\boldsymbol{\xi}}}
\newcommand{\bs}{{\mathbf{s}}}
\newcommand{\bmu}{{\boldsymbol{\mu}}}

\newcommand{\hOm}{{\hat{\Omega}}}
\newcommand{\Om}{{\Omega}}

\usepackage[font=footnotesize,labelfont=bf]{caption}
\usepackage{multirow}

\usepackage{pxfonts}

\usepackage{import}

\graphicspath{./graphics/}

\usepackage[textsize=scriptsize]{todonotes}
\usepackage{setspace}

\makeatletter
\def\ps@pprintTitle{%
 \let\@oddhead\@empty
 \let\@evenhead\@empty
 \def\@oddfoot{}%
 \let\@evenfoot\@oddfoot}
\makeatother

\journal{Computer Methods in Applied Mechanics and Engineering}

\begin{document}

\begin{frontmatter}

\title{On the Spline-Based Parameterisation of Plane Graphs via Harmonic Maps}

\author[add1]{Jochen Hinz\corref{cor1}}
\ead{jochen.hinz@epfl.ch}

\cortext[cor1]{Corresponding author}

\address[add1]{Institute of Mathematics, Ecole Polytechnique F\'ed\'erale de Lausanne, 1015 Lausanne, Switzerland.}

\begin{abstract}
This paper presents a spline-based parameterisation framework for plane graphs. The plane graph is characterised by a collection of curves forming closed loops that fence-off planar faces which have to be parameterised individually. Hereby, we focus on parameterisations that are conforming across the interfaces between the faces. Parameterising each face individually allows for the imposition of locally differing material parameters which has applications in various engineering disciplines, such as elasticity and heat transfer. For the parameterisation of the individual faces, we employ the concept of harmonic maps. The plane graph's spline-based parameterisation is suitable for numerical simulation based on isogeometric analysis or can be utilised to extract arbitrarily dense classical meshes. Application-specific features can be built into the geometry's mathematical description either on the spline level or in the mesh extraction step.
\end{abstract}

\end{frontmatter}

\section{Introduction}
\label{sect:Introduction}

In practical applications, the geometry is typically represented by no more than a CAD-based description of the geometry's boundary contours, which is itself given by a collection of spline curves or dense point clouds. As such, virtually all numerical techniques rely heavily on routines that address the problem of finding a valid mathematical description of the geometry's interior from no more than the provided CAD input. Besides nondegeneracy, a parameterisation pipeline furthermore strives for parametric descriptions with application-dependent numerically favourable properties such as cell size homogeneity or orthogonalised isolines. \\
The parameterisation problem is, despite several decades of research, still considered a major robustness bottleneck in computational simulation workflows. It often contributes substantially to the overall computational costs and is notoriously difficult to automate. This is especially critical in light of the fact that numerical artifacts manifest themselves in the form of cell degeneracy, which rules out the possibility to perform the analysis step completely. \\
The associated difficulties are further exacerbated in settings that require a mathematical description of several geometries, in the planar case referred to as \textit{faces}. Dividing a planar geometry into a set of faces is necessary in settings that impose, for instance, locally differing material parameters. Since the faces may intersect at shared interfaces, most numerical techniques favour a parametric description that preserves interface conformity. That is, wherein bordering faces' interface elements share the same (often piecewise linear) interface curve. However, in geometrically complex settings element conformity may be difficult to accomplish which is why numerical schemes may have to be combined with mortaring techniques \cite{bernardi2005basics}, which are themselves associated with a number of challenges. \\
Since the onset of isogeometric analysis \cite{hughes2005isogeometric}, the spline-based parameterisation problem, including multipatch region parameterisation, has received an increased amount of interest from the scientific community \cite{buchegger2017planar, falini2015planar, gravesen2012planar, xu2013constructing, xu2010optimal, xu2011parameterization}. Algorithms are based on algebraic techniques \cite{farin1999discrete, varady2011transfinite, salvi2014ribbon}, quality cost function optimisation \cite{falini2015planar, gravesen2012planar, buchegger2017planar, xu2010optimal, xu2011parameterization, xu2013constructing, ji2022penalty, ji2021constructing, wang2021smooth}, frame field \cite{hiemstra2020towards, shepherd2023quad, zhang2024multi} and PDE-based approaches \cite{shamanskiy2020isogeometric, hinz2018elliptic, hinz2024use}. \\

\noindent Unfortunately, so far, most of these efforts have been restricted to the parameterisation of a single geometry from one or several patches. \\

\noindent To overcome these limitations, this paper presents a general framework for the spline-based parameterisation of plane graphs, respecting the requirement that each face be parameterised individually, while furthermore retaining conforming interfaces. For the parameterisation of individual faces, we adopt the harmonic map-based techniques from \cite{hinz2024use} and discuss several techniques capable of bringing the plane graph into a form that allows it to be parameterised using quad layouts. \\

\noindent While this paper focuses on techniques that yield spline-based parametric descriptions for use in isogeometric analysis applications, we nevertheless assume the plane graph to come in the form of a collection of dense point sets which are then converted back into a set of spline curves via curve fitting. We have found point-set based inputs to be better suited for this paper's techniques since they are more easily manipulated and can be combined with standard polygon routines. Furthermore, interface conformity is facilitated by the use of dyadically refined knotvectors in the fitting stage, a requirement that may not be respected by a CAD-based spline input. \\
Clearly, for spline-based inputs the presented techniques remain applicable since each spline can be converted into an arbitrarily dense point set via collocation.

\subsection{Notation}
This paper denotes vectors in boldface. The $i-th$ entry of a vector is denoted by $x_i$. Similarly, the $ij$-th entry of a matrix is denoted by $A_{ij}$. Let $\mathbf{y}: \Om \rightarrow \mathbb{R}^m$ and $\bx: \Om \rightarrow \mathbb{R}^n$. We frequently work with vector spaces $\mathcal{V}$. By default, we employ the abuse of notation
\begin{align}
    \mathcal{V}^n = \underbrace{ \mathcal{V} \times \cdots \times \mathcal{V} }_{n \text{ terms}}
\end{align}
and similarly for tensorial spaces, i.e., $\mathcal{V}^{n \times n}$. Sobolev spaces are denoted by $H^s(\Omega)$ and vectorial Sobolev spaces by $H^s(\Omega, \mathbb{R}^n)$. \\
By $\operatorname{Int}\left(D \right)$, we denote the interior of a closed domain $D$, while $\overline{\Om}$ denotes the closure of an open domain $\Om$. The operator $| \, \cdot \, |$ applied to a set returns its number of elements.

\subsection{Problem Statement}
\label{subsect:problem_statement}
\begin{figure}[h!]
\centering
\includegraphics[width=0.8\textwidth]{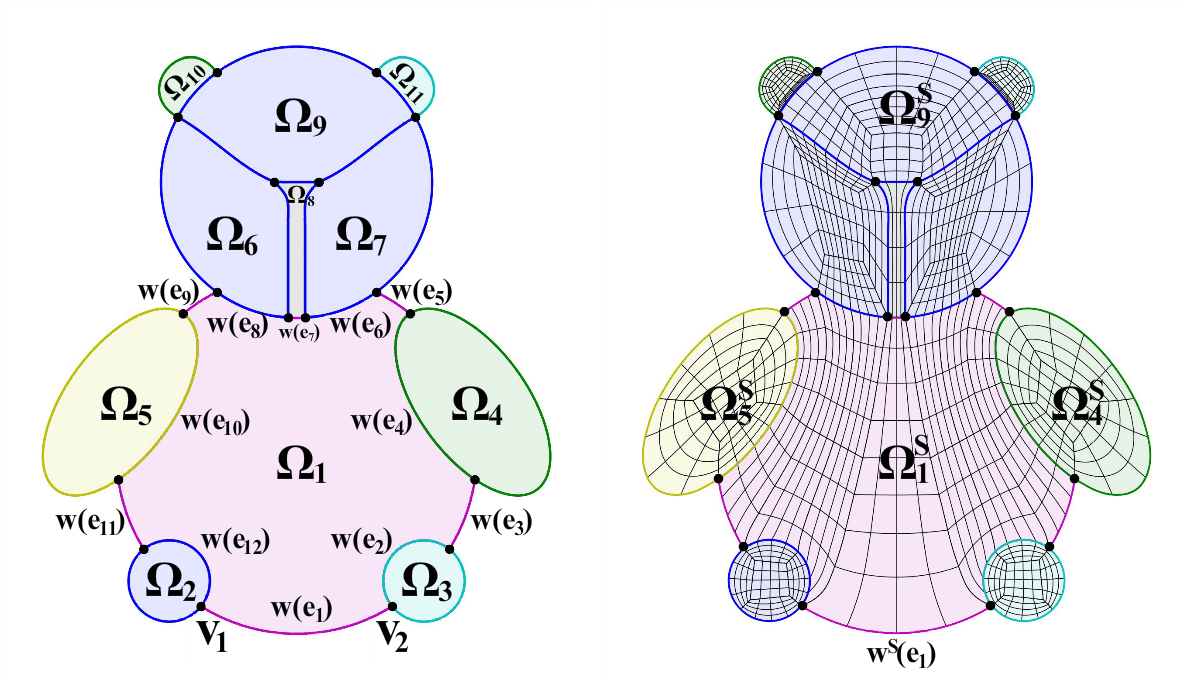}
\caption{An example of a plane graph (left) that is visualised via its weight function $w(\, \cdot \,)$. The various polygons $\Omega_i$ that are bounded by the point sets $w(e), \, e \in F_i$ are highlighted in different colors. The graph contains several examples of two edges connecting the same two vertices, for example $e_{10}$ and $e_{13}$. The figure shows that $F_1 = \left(\pm e_1, \pm e_2, \ldots, \pm e_{12} \right)$. The right figure shows a possible spline-based parameterisation wherein each point set $w(e)$ is replaced by an accurate spline approximation $w^S(e)$ fencing off the spline domains $\Omega_i^S \approx \Omega_i$. Each face is parameterised from one or several patches using a harmonic map, while neighboring faces use the same spline edges $w^S(e)$, thus leading to a conforming interface.}
\label{fig:example_plane_graph}
\end{figure}

We are given a graph G associated with a set $V \subset \mathbb{R}^2$ of vertices and a set $E$ containing directed edges. Letting $\mathcal{I}_v = \{ 1, 2, \ldots, N_v \}$ denote the index-set of vertices, we have $V = \{v_i\}_{i \in \mathcal{I}_v}$, where each $v_i \in V$ is a unique vertex $v_i \in \mathbb{R}^2$. The set $E = \{e_1, \ldots, e_{N_e} \}$ contains edges incident to vertices. The operator $\iota: E \rightarrow \mathbb{N} \times \mathbb{N}$ returns the indices of the vertices incident to an edge $e \in E$. The set of edges connecting two vertices $\{v_i, v_j \}$ is denoted by $E_{i, j}$. Typically, $\vert E_{i, j} \vert \in \{0, 1\}$ but we note that $\vert E_{i, j} \vert > 1$ is possible, see Figure \ref{fig:example_plane_graph}. The same operator may be applied to (ordered) sets of edges element-wise, i.e., $\iota(\{e_1, \ldots, e_N\}) = \{\iota(e_1), \ldots, \iota(e_N) \}$. \\
There is a weight function $w(e), \, e \in E$ that assigns an ordered point set to each edge. If $\iota(e) = (i, j)$, this means that
$$w(e) = \left(p_1, \ldots, p_{N} \right), \quad \text{with} \quad p_1 = v_i \in V \quad \text{and} \quad p_{N} = v_j \in V.$$ By default, we assume that the point sets $w(e)$, $e \in E$ are sufficiently smooth to be well-approximated by a $C^1$ spline curve with a reasonable number of degrees of freedom (DOFs). Given $\iota(e) = (i, j)$, we denote by $-e$ the directed edge with the property that $\iota(-e) = (j, i)$ and use the convention that $w(\, \cdot \,)$ assigns to $-e$ the point set that reverses the order of the $p_i \in w(e)$. Applying the operator $U( \, \cdot \, )$ to an edge $\pm e$, with $e \in E$, always returns the associated edge $e \in E$, i.e., $U(\pm e) = e \in E$. As before, the same operator can be applied to sets / sequences of edges, i.e., $U(\{ \pm e_{i_1}, \ldots, \pm e_{i_N} \} ) = \{e_{i_1}, \ldots, e_{i_N} \} \subset E$. For a graph $G = (V, E)$, the ordered set $\partial E = (\pm e_{i_1}, \ldots, \pm e_{i_{N_\partial}} )$ represents the sequence of edges $e_{i_j} \in E$ that lie on the graph's boundary, ordered in counterclockwise orientation. \\

\noindent The graph $G$ is furthermore associated with a collection of $N_F$ faces $\mathcal{F} = \{F_1, \ldots, F_{N_F}\}$, i.e. $G = (V, E, \mathcal{F})$. The faces are characterised by ordered sets $F_i = \left( e_{i_1}, \ldots, e_{i_{N_i}} \right)$, with the property that if  $e \in F_i$, then $U(e) \in E$. The faces $F_i \in \mathcal{F}$ represent (dense) non self-intersecting  and mutually disjoint polygons $\Omega_i$ whose boundaries result from concatenating the point sets $w(e), \, e \in F_i$ in the order as they appear in $F_i$ and removing duplicate vertices. We assume that each $F_i$ is oriented in counterclockwise orientation and that the sequences form a closed loop, i.e., if $\iota(e_{i_k}) = (p, q)$, then $\iota(e_{i_{k+1}}) = (q, r)$ (where the subscript is taken modulus $\vert F_i \vert$). As such, each $\Omega_i$ is bounded by a piecewise linear Jordan curve $\partial \Omega_i$. \\

\noindent Since the point sets $w(e), \, e \in E$ define the vertices of the non self-intersecting and mutually disjoint $\Omega_i$, we have 
$$w(e_q) \cap w(e_r) \in V \quad \text{or} \quad w(e_q) \cap w(e_r) = \emptyset, \quad \text{ while} \quad U(F_i) \cap U(F_j) \subset E \quad \text{or} \quad U(F_i) \cap U(F_j) = \emptyset.$$
Finally, the operator $\mathbb{V}(\, \cdot \,)$, when applied to an $F_i \in \mathcal{F}$, returns the sequence of vertices $v_j \in V$ in the order they appear in $F_i$, i.e., if $\iota(F_i) = [(p, q), (q, r), \ldots, (\alpha, p)]$, then $\mathbb{V}(F_i) = (v_p, v_q, \ldots, v_{\alpha}) \subset V$. The same operator can also be applied to $\partial E$. \\

\noindent Given a planar graph $G = (V, E, \mathcal{F})$ with the aforementioned properties, the purpose of this paper is providing a framework for the spline-based parameterisation of $G$, respecting the requirement that each face $F_i \in \mathcal{F}$ be parameterised individually, see Figure \ref{fig:example_plane_graph} right. For this purpose, we assign an $N_{F_i} := \vert F_i \vert$-sided polygonal parametric domain $\hOm_i \subset \mathbb{R}^2$ to each face $F_i$. By $\hOm(N)$ we denote the regular $N$-sided polygon of radius one. Unless stated otherwise, we utilise $\hOm_i = \hOm(N)$ with $N = \vert F_i \vert$. \\

\noindent We denote the set of edges that make up $\partial \hOm_i$ by the ordered set $\partial E_i$ (ordered in counterclockwise direction). Assigning the $N_{F_i}$ edges of $\partial \hOm_i$ to the $N_{F_i}$ edges of $F_i$ in ascending order induces a canonical bijective edge correspondence $\phi_i: \partial E_i \rightarrow F_i$. This also induces a correspondence between the vertices $\mathbb{V}(\partial E_i)$ and $\mathbb{V}(F_i)$ and we utilise the same symbol $\phi_i: \mathbb{V}(\partial E_i) \rightarrow \mathbb{V}(F_i)$ for convenience. We note that the edge /vertex correspondence is unique up to a cyclic permutation of $\partial E_i$ or $F_i$. \\

\noindent The weight function $w^S(\, \cdot \,)$ assigns to each $e \in E$ an open spline curve $s \in C^1([0, 1], \mathbb{R}^2)$ that connects the vertices $v_i$ and $v_j$, with $\iota(e_q) = (i, j)$. As before, the spline curve $-s$ reverses the orientation of $s$. We denote the domains that are fenced off by the spline curves $s = w^S(e)$, associated with the faces $F_i \in \mathcal{F}$, by $\Omega_i^S$. With that, the purpose of this paper is presenting a framework for parameterising $G$ such that:
\begin{enumerate}[label=\Alph*.]
    \item The spline curve $s = w^S(e)$ approximates the point set $p = w(e)$ to a user-defined accuracy;
    \item Each $s = w^S(e)$ is associated with a knotvector that is a (possibly repeated) dyadic refinement of a user-defined base knotvector $\Xi^0$;
    \item The map $\mathbf{x}_i: \hOm_i \rightarrow \Om_i^S$ is bijective and maps the side of $\partial \hOm_i$ represented by $\hat{e} \in \partial E_i$ onto the spline curve $s = w^S(e)$, where $e = \phi_{i}(\hat{e})$;
    \item The parametric properties of $\mathbf{x}_i: \hOm_i \rightarrow \Om_i^S$ can be tuned.
\end{enumerate}

\noindent The workflow presented in this paper is comprised of the following stages:
\begin{enumerate}
    \item Graph preparation and templatisation;
    \item Curve fitting and face parameterisation;
    \item (Optional) Post processing.
\end{enumerate}

\noindent Stage 1. performs operations on the graph that allow each face $F_i \in \mathcal{F}$ to be parameterised from a parametric domain $\hOm_i$ with a multipatch layout. For instance, each face requires an even number of edges which may require an edge $e \in F_i$ to be split in two. The same stage assigns a graph $T_i = (V_i, E_i, \mathcal{Q}_i)$, a so-called \textit{template}, with four-sided faces $q \in \mathcal{Q}_i$ to each $F_i \in \mathcal{F}$. The graph $T_i$ represents a quadrangulation of a regular $N_{F_i} := \vert F_i \vert$-sided polygon $\hOm_i = \hOm(N_{F_i})$. Upon finalisation of Stage 1., the updated plane graph $G = (V, E, \mathcal{F})$ typically differs from the initial graph. In order to be compatible with this paper's parameterisation strategy, this stage also removes concave corners between adjacent point sets $(w(e_1), w(e_2))$ with $\{e_1, e_2\} \subset F_i$ by constructing a splitting curve originating from the concave vertex that splits the face in two. Here the notion of concave corner is based on a suitable discrete notion of the jump in the tangent vector on both sides of the vertex incident to both $e_1$ and $e_2$ (see Section \ref{sect:Discrete_Tangent}). A suitable template $T_i = (V_i, E_i, \mathcal{Q}_i)$ is chosen from a large pre-computed catalogue $\mathbb{T}$ (see Section \ref{sect:patch_adjacency_graphs}). \\

\noindent Stage 2. generates the weight function $w^S(\, \cdot \,)$ which assigns a spline curve $s \in C^1([0, 1], \mathbb{R}^2)$ to all $e \in E$. Each $s = w^S(e)$ results from a (possibly repeated) dyadic refinement of a (typically uniform) open base knotvector $\Xi^0$ over the unit interval. Refinements are performed to ensure that $s = w^S(e)$ approximates $p = w(e)$ sufficiently accurately, while the refinement's dyadic nature reduces the interior knot density of the knotvectors we assign to the $\hat{e} \in E_i$ of $T_i$ to capture the induced boundary correspondence $\mathbf{f}_i: \partial \hOm_i \rightarrow \partial \Om_i^S$. \\
\noindent Upon finalisation of $w^S(\, \cdot \,)$, Stage 2. employs a suitable numerical algorithm to compute the bijective maps $\mathbf{x}_i: \hOm_i \rightarrow \Om_i^S, \, \, i \in \{1, \ldots, N_{F}\}$ each satisfying $\mathbf{x}_i\vert_{\partial \hOm_i} = \mathbf{f}_i$. In this paper, the parameterisation stage is based on the concept of harmonic maps. \\ 

\noindent The (optional) final stage optimises the parameterisations starting from an initial folding-free map for each individual face. Optimisation may target, for instance, cell size homogeneity and is based on variants of plain harmonic maps. The generated surface splines may then be used for analysis using IGA tools or may optionally be converted into a (orders of magnitude) denser classical mesh via collocation.


\subsection{Related Work}
\label{sect:related_work}
The spline-based multipatch region parameterisation problem has garnered renewed attention with the advent of IGA \cite{hughes2005isogeometric}. The prevailing methods are broadly divided into two broad categories: 1) segmentation-based approaches or 2) direct multipatch parameterisation approaches. Type 1) methods partition the domain into disjoint pieces, each diffeomorphic to the unit quadrilateral. The resulting curvilinear quad segments are then parameterised individually using a singlepatch parameterisation algorithm. The segmentation itself may rely on techniques such as, cross / frame fields \cite{viertel2019approach, bommes2009mixed, hiemstra2020towards, shepherd2023quad, zhang2024multi} or triangulation- / quadrangulation-assisted approaches \cite{xiao2018computing, xu2018constructing}. \\
Conversely, type 2) methods choose a multipatch layout that is compatible with the number of boundary curve segments after which an algorithm for the interior parameterisation is employed. In this context, the multipatch interface control points become degrees of freedom in the formulation. The interior parameterisation relies on multipatch generalisations of singlepatch techniques such as optimisation- \cite{buchegger2017planar, bastl2021planar} or PDE-based approaches \cite{hinz2024use, shamanskiy2020isogeometric} while the multipatch layout is either a manual choice \cite{hinz2024use, shamanskiy2020isogeometric} or based on more deterministic methods such as patch adjacency graphs \cite{buchegger2017planar} or skeleton-based approaches \cite{bastl2021planar}. \\

\noindent These approaches primarily address the parameterisation of single domains, either simply- or multiply-connected, often featuring complex boundary curves. Extending aforementioned techniques to encompass collections of such regions, each with prescribed interface curves, introduces new challenges. To the best of our knowledge, these challenges have not yet been systematically addressed in the literature. \\

\noindent This work builds upon the PDE-based parameterisation techniques introduced in \cite{hinz2024use} and extends them with a mechanism that facilitates a deterministic selection of the multipatch layout. These techniques are then integrated into a framework that enables the paramterisation multiple, mutually disjoint regions while maintaining a conforming interface.

\subsection{Discrete Tangent and Normal Vectors}
\label{sect:Discrete_Tangent}
For introducing a discrete tangent, we treat the point set $p = w(e)$ with $U(e) \in E$ as a piecewise linear spline curve that crosses all $p_i \in p$ parameterised over $[0, 1]$ using the chord-length abscissae associated with $p$. Its length is denoted by $L(p)$. The discrete tangent $t_{i}$ on the linear interval between $\left(p_i, p_{i+1} \right) \subseteq p$ is the vector of length $\| p_{i+1} - p_i \|_2$ pointing in the direction of $p_{i+1} - p_i$. Its normalised counterpart is denoted by $\hat{t}_i$. Given $\mathcal{F} \ni F = (\ldots, e_\alpha, e_\beta, \ldots)$, with $\iota(e_\alpha) = (p, q)$ and $\iota(e_\beta) = (q, r)$, we denote by $\hat{t}_{-}(v_{q})$ and $\hat{t}_{+}(v_{q})$ the discrete unit tangents on either side of the vertex $v_q \in V$ as defined by the adjacent linear edge in $p_\alpha = w(e_\alpha)$ and $p_\beta = w(e_\beta)$, respectively. The discrete unit outward normal on the face $F \in \mathcal{F}$ at vertex $v_q \in V$ follows from the average $\hat{t}(v_q) := (\hat{t}_{-}(v_{q}) + \hat{t}_{+}(v_{q})) / 2$ and is given by $\hat{n}(v_q) = \left( \hat{t}_2(v_q), -\hat{t}_1(v_q) \right)^T$. The discrete interior angle at vertex $v_q \in \mathbb{V}(F)$ in $F \in \mathcal{F}$ is denoted by $\angle(v_q)$. The discrete angle $\angle(\, \cdot \,)$ will be utilised to flag a vertex $v_q \in \mathbb{V}(F)$ as concave while the discrete normal is employed to construct splitting curves for the removal of concave corners in Section \ref{sect:prep_temp}. Moreover, the discrete angle $\angle( \, \cdot \,)$ plays a central role in the selection of a suitable template $T \in \mathbb{T}$ for each $F \in \mathcal{F}$.

\subsection{Harmonic Maps}
\label{sect:harmonic_maps}
Let $F \in \mathcal{F}$ be a face with associated spline domain $\Om^S$. \\
To find a valid parameterisation $\mathbf{x}: \hat{\Omega} \rightarrow \Omega^S$ between an $N$-sided ($N$ even) convex polygonal parametric domain $\hat{\Omega}$ and the $N$-sided target domain $\Omega^S$ given no more than a (piecewise smooth, spline-based) boundary correspondence $\mathbf{f}: \partial \hat{\Omega} \rightarrow \partial \Omega^S$, this paper adopts the concept of harmonic maps. The general idea of harmonic maps is requiring the inverse map $\mathbf{x}^{-1}: \Omega^S \rightarrow \hat{\Omega}$ to be a pair of harmonic functions in $\Omega^S$. This choice is motivated by the following famous result:

\begin{theorem}[Rad\'o-Kneser-Choquet]
\label{thrm:RKC}
    The harmonic extension of a homeomorphism from the boundary of a Jordan domain $\Om^S \subset \mathbb{R}^2$ onto the boundary of a convex domain $\hOm \subset \mathbb{R}^2$ is a diffeomorphism in $\Om$. \\
\end{theorem}
For proofs, we refer to \cite{kneser1926losung, choquet1945type, sauvigny1991embeddedness, castillo1991mathematical}. \\

\noindent Since we have the luxury of choosing $\hat{\Omega}$ convex, Theorem \ref{thrm:RKC} applies to the pair $\left(\hat{\Omega}, \, \Omega^S \right)$ even though $\Omega^S$ is generally nonconvex. However, we mention that despite being diffeomorphic in the interior, a harmonic $\mathbf{x}^{-1}: \Omega^S \rightarrow \hat{\Omega}$ may fail to extend diffeomorphically to the closure $\overline{\Omega^S}$, and similarly for its inverse $\mathbf{x}: \hat{\Omega} \rightarrow \Omega^S$. As such, the map's Jacobian determinant
\begin{align}
    \det J(\mathbf{x}) := \det \frac{\partial \mathbf{x}}{\partial \boldsymbol{\xi}} \quad \text{may fail to satisfy} \quad c < \det J(\mathbf{x}) < C \quad \text{almost everywhere in } \overline{\hat{\Omega}},
\end{align}
for some $0 < c \leq C < \infty$. We regard this as undesirable from a parameterisation quality perspective. Since the vertex corner pattern of $\Omega^S$ is equivalent to that of the vertices $v \in \mathbb{V}(F)$ by construction, the classification of the domain $\Omega^S$ follows from the classification of the polygon $\Omega$ associated with the point sets $w(e), \, e \in F$. \\

\noindent Clearly, a necessary condition for (an inversely harmonic) $\mathbf{x}: \hOm \rightarrow \Om^S$ to extend diffeomorphically to the closure of the domain under the boundary correspondence $\mathbf{f}: \partial \hOm \rightarrow \partial \Omega^S$, is the existence of a diffeomorphism $\mathbf{F}: \Omega_1 \rightarrow \Omega_2$
$$\text{ that maps a neighborhood } \Omega_1 \supset \overline{\hat{\Omega}} \text{ onto some } \Omega_2 \supset \overline{\Omega^S} \text{ while satisfying } \mathbf{F} \vert_{\partial \hOm} = \mathbf{f}.$$
As such, $\mathbf{f}: \partial \hOm \rightarrow \Omega^S$ must map convex corners onto convex corners and similarly for concave corners. However, the existence of concave corners in $\partial \hOm$ violates the requirement that $\hOm$ be convex, hence the requirement that $\partial \Omega^S$ contain only convex corners. \\

\noindent To enable the computation of a uniformly diffeomorphic harmonic map, we introduce a control map $\mathbf{r}: \hOm \rightarrow \hOm^{\mathbf{r}}$ which maps onto a convex domain $\hOm^{\mathbf{r}}$ that is diffeomorphic to $\Omega^S$. Depending on the number of vertices $v \in \mathbb{V}(F)$ that we model as (convex) corners, we utilise different layouts. Given a user-defined threshold value $\mu_\angle > 0$, we model a vertex $v \in \mathbb{V}(F)$ as a corner if $\angle(v) < \pi - \mu_\angle$. Typically, we take $\mu_\angle \ll 1$. Depending on the number 
$$N_{\text{convex}} := \left \vert \{v \in \mathbb{V}(F) \, \big \vert \, \angle(v) < \pi - \mu_{\angle} \} \right \vert,$$
we construct $\hOm^{\mathbf{r}}$ from a total of five different layouts. \\

\begin{figure}[h!]
\centering
    \begin{subfigure}[b]{0.98\textwidth}
        \centering
        \includegraphics[align=c, height=3cm]{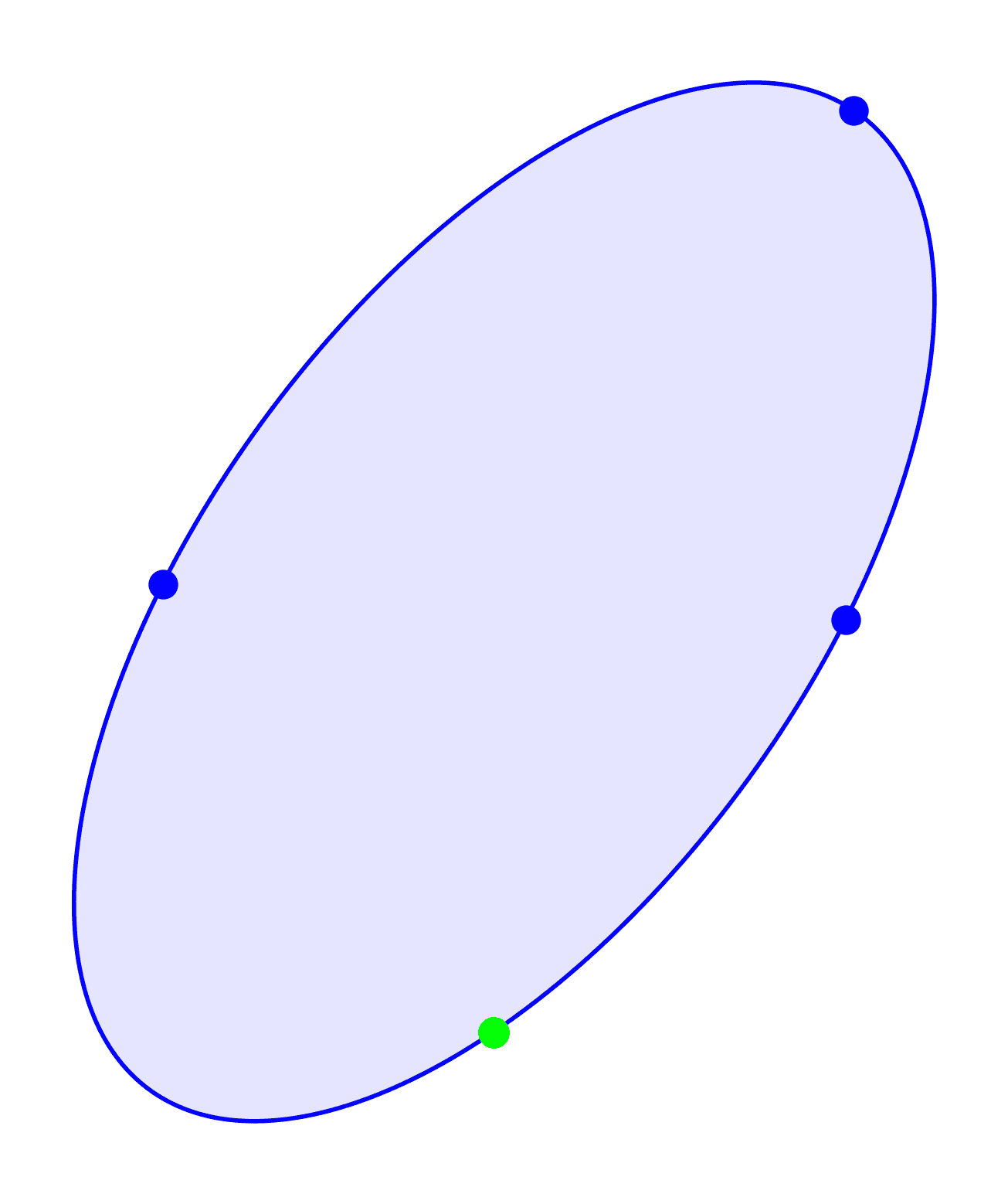}
        \includegraphics[align=c, height=3cm]{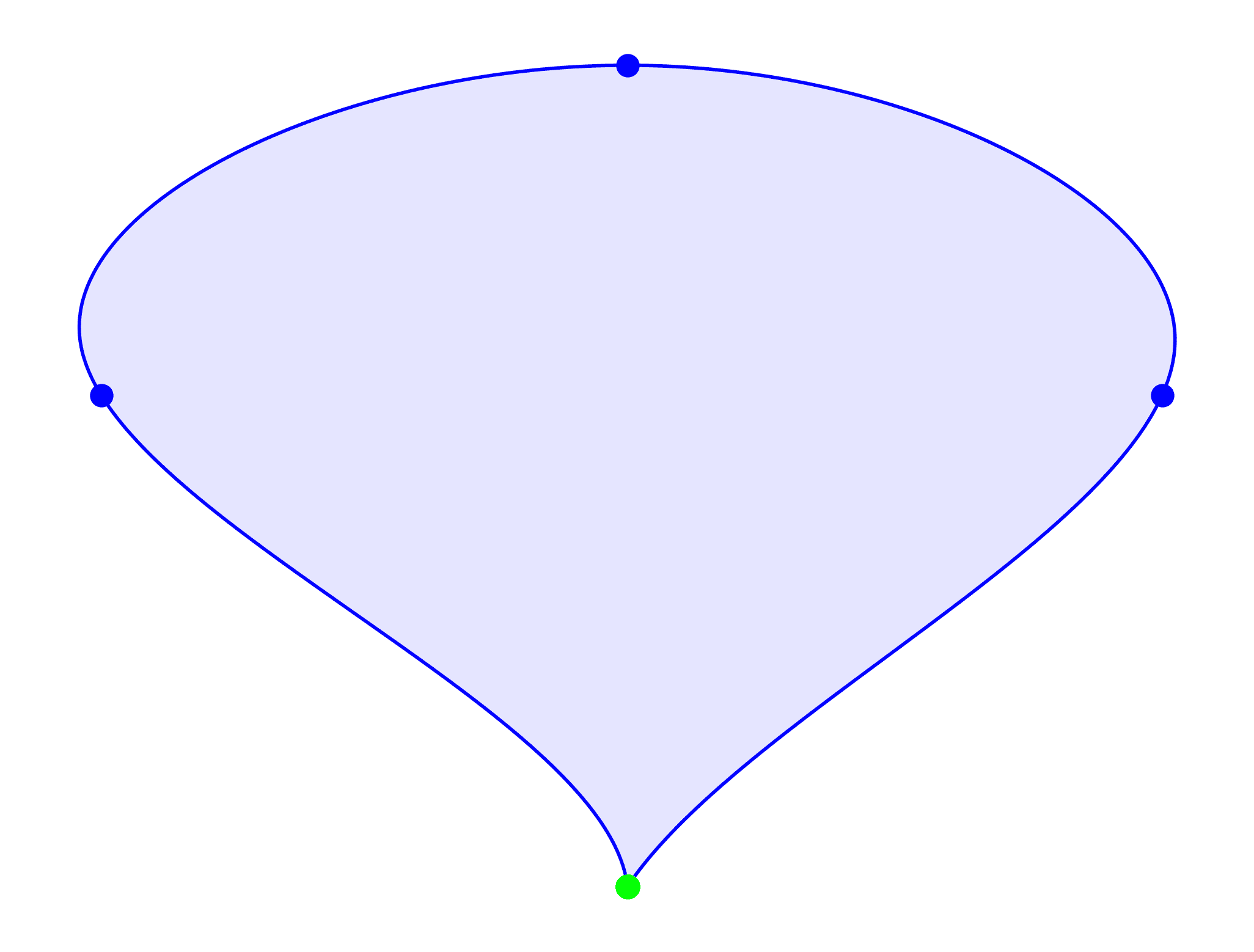}
        \includegraphics[align=c, height=3cm]{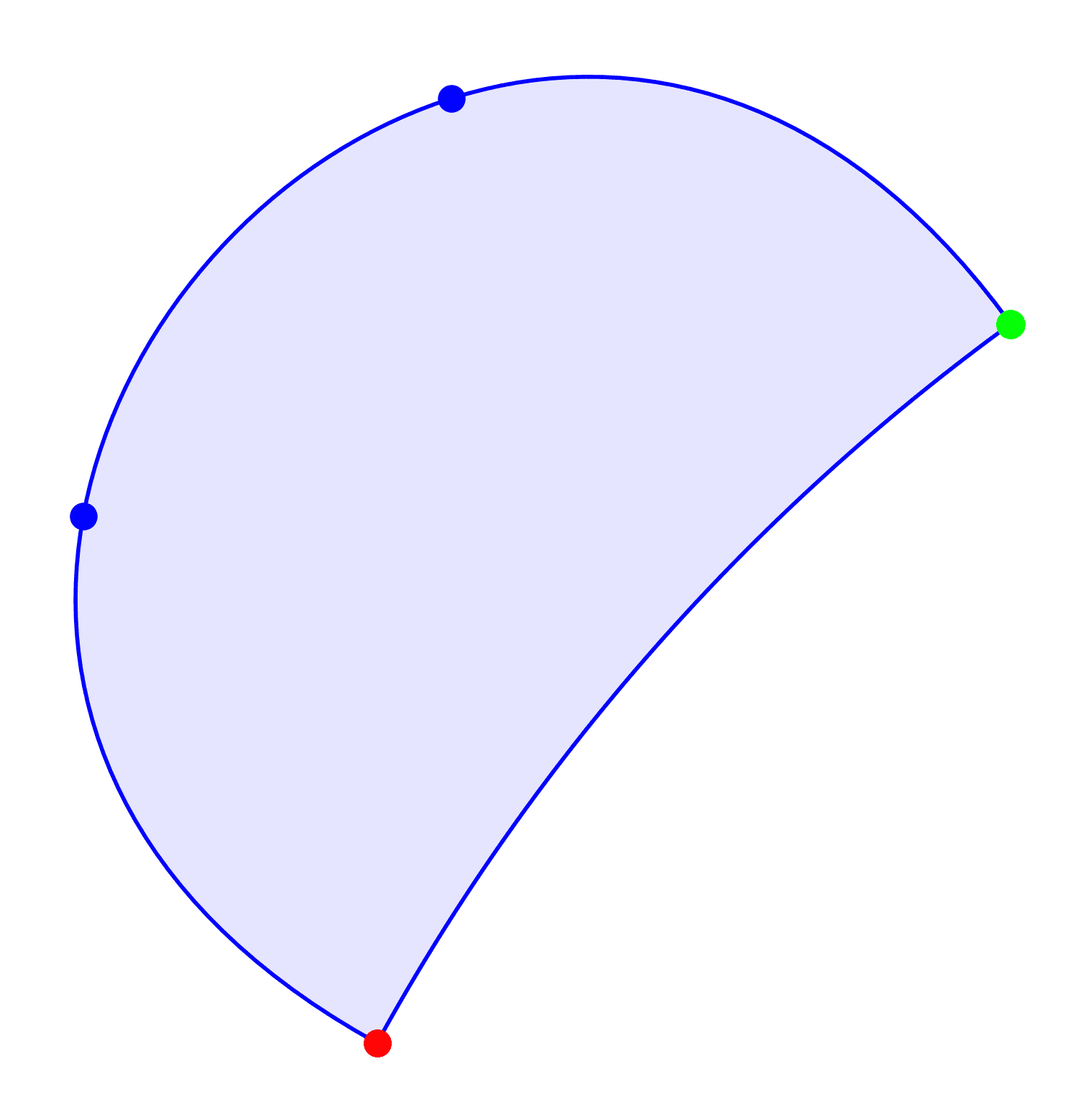}
        \includegraphics[align=c, height=3cm]{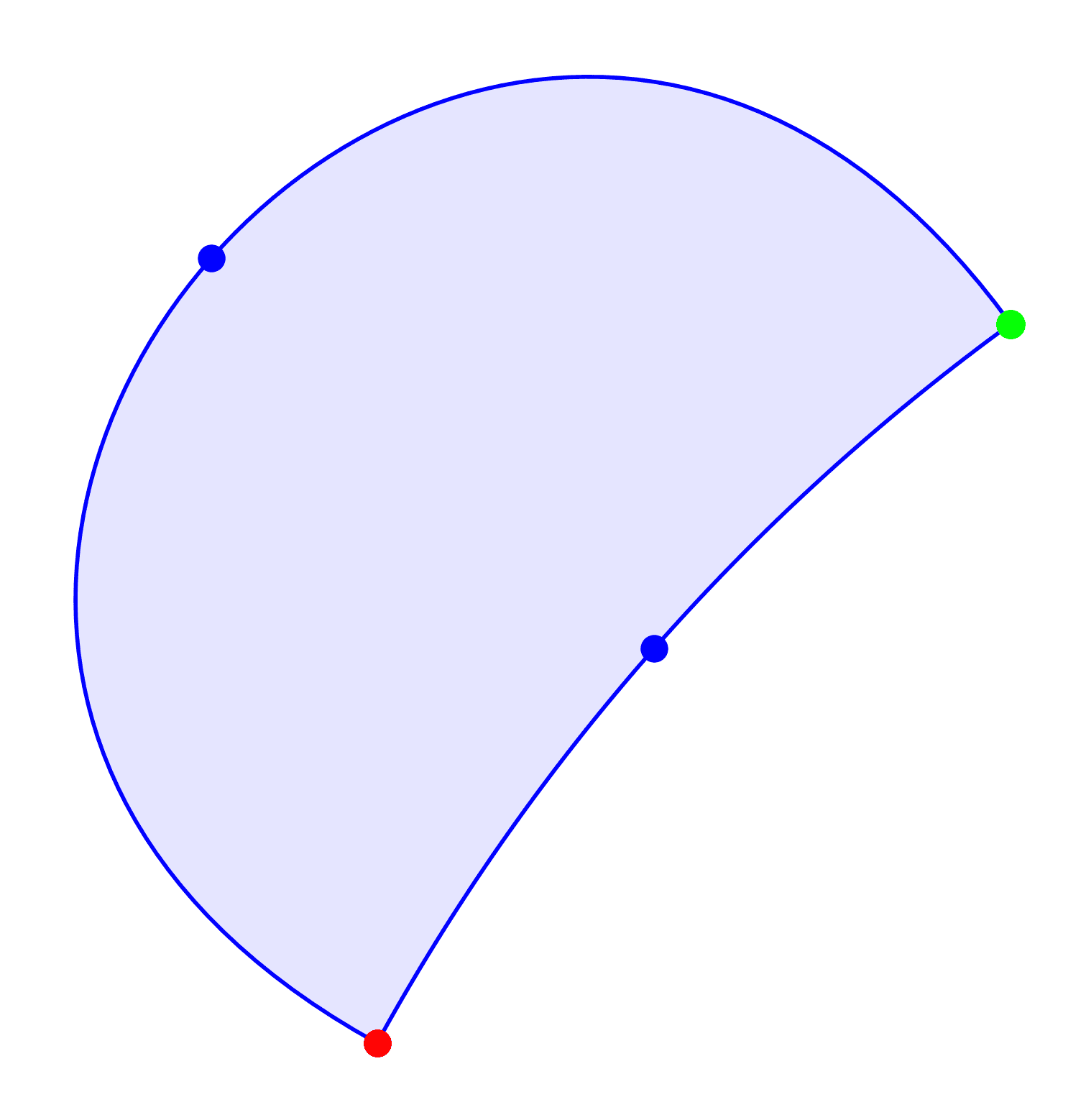}
        \includegraphics[align=c, height=3cm]{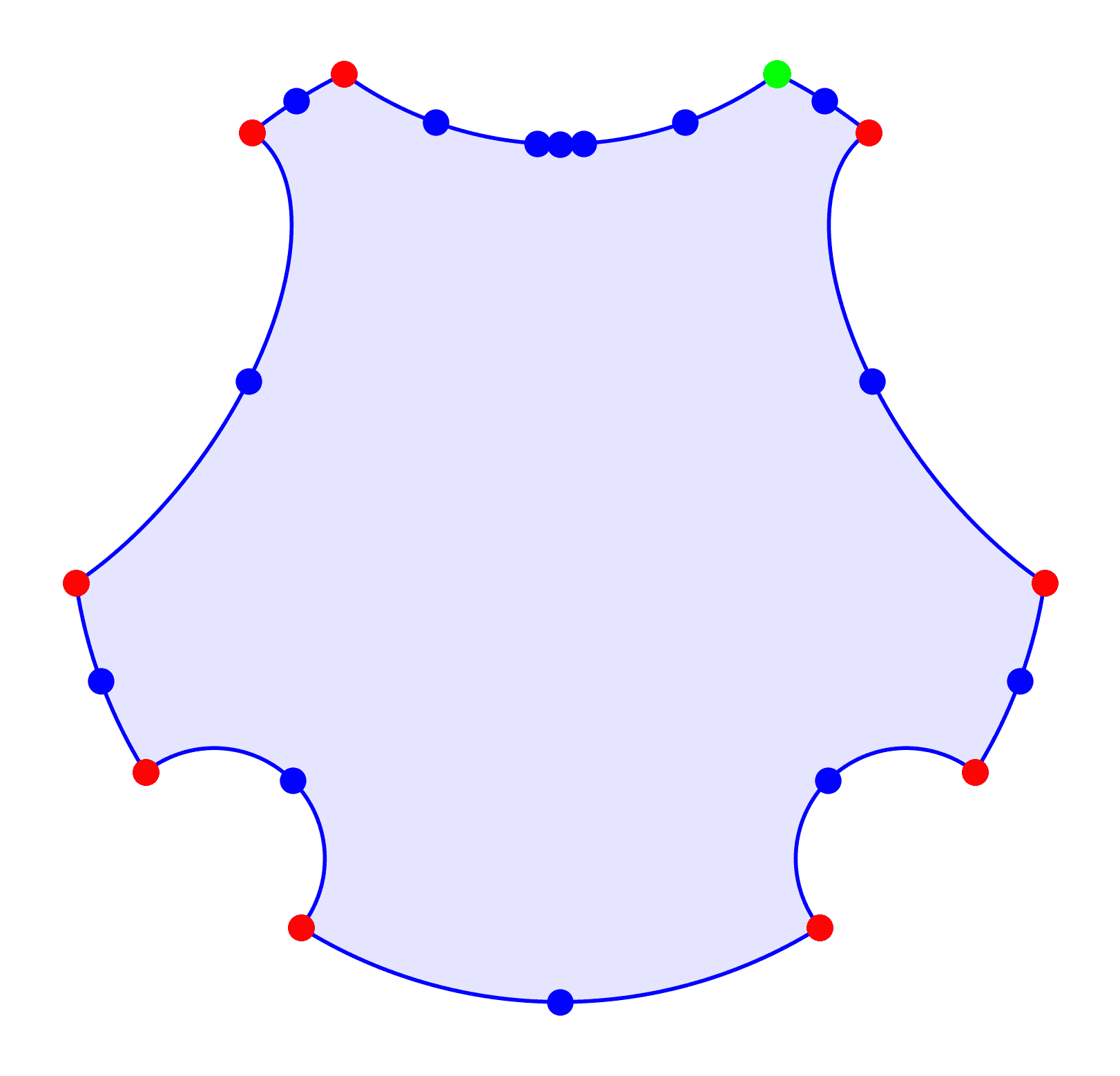}
    \end{subfigure}
    \begin{subfigure}[b]{0.95\textwidth}
        \centering
        \includegraphics[align=c, height=2.8cm]{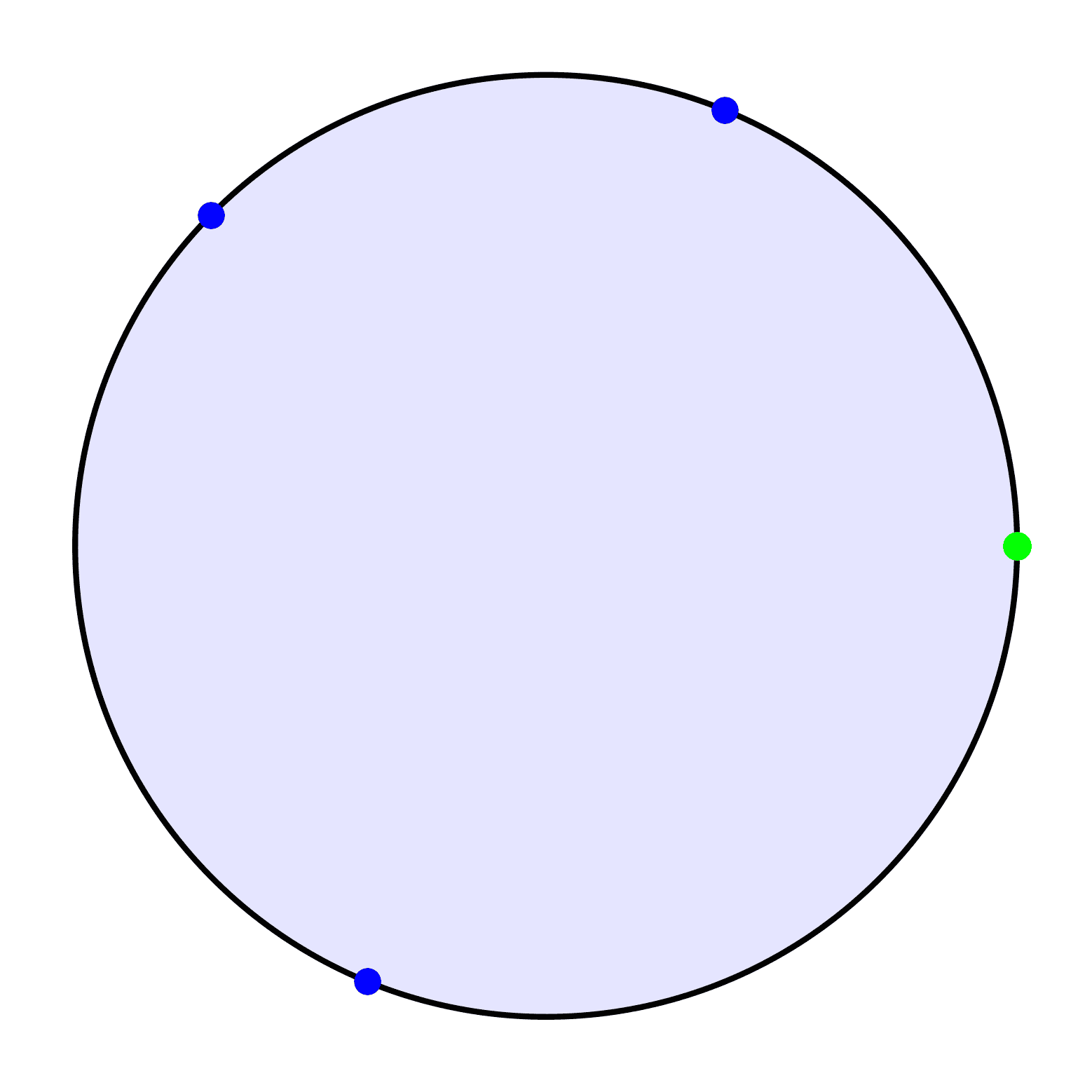}
        \includegraphics[align=c, height=1.7cm]{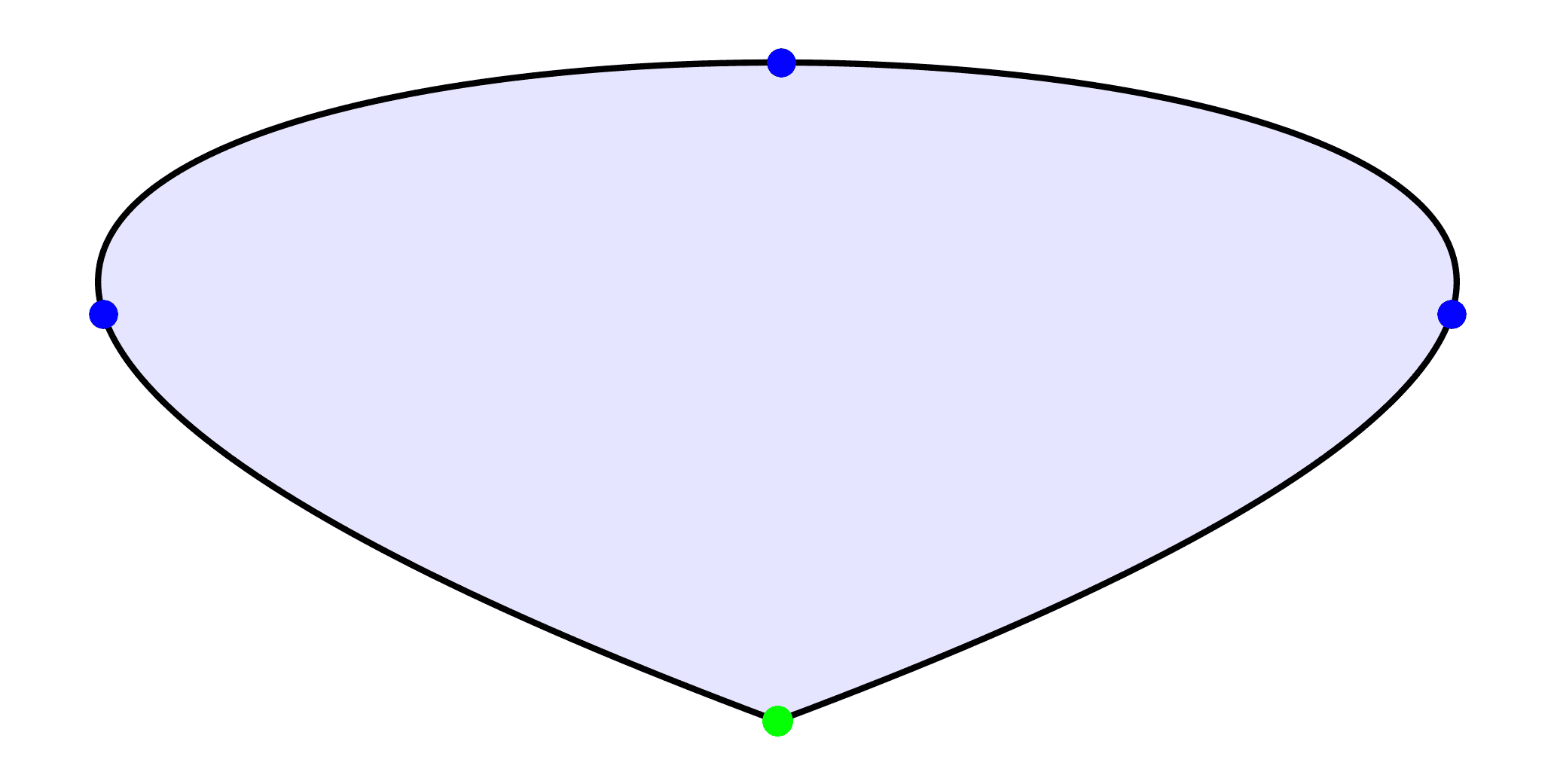}
        \includegraphics[align=c, height=1.5cm]{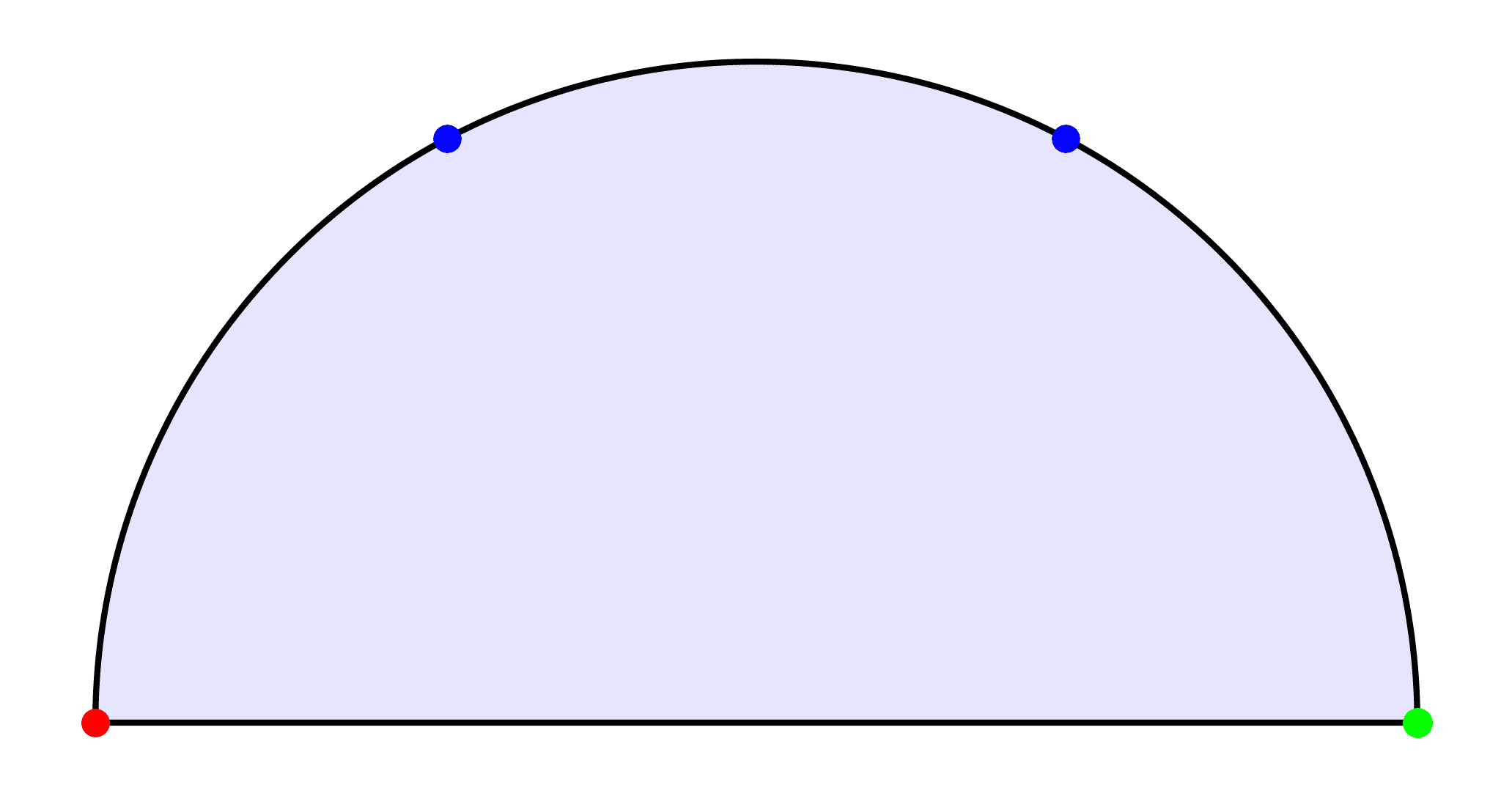}
        \includegraphics[align=c, height=1.2cm]{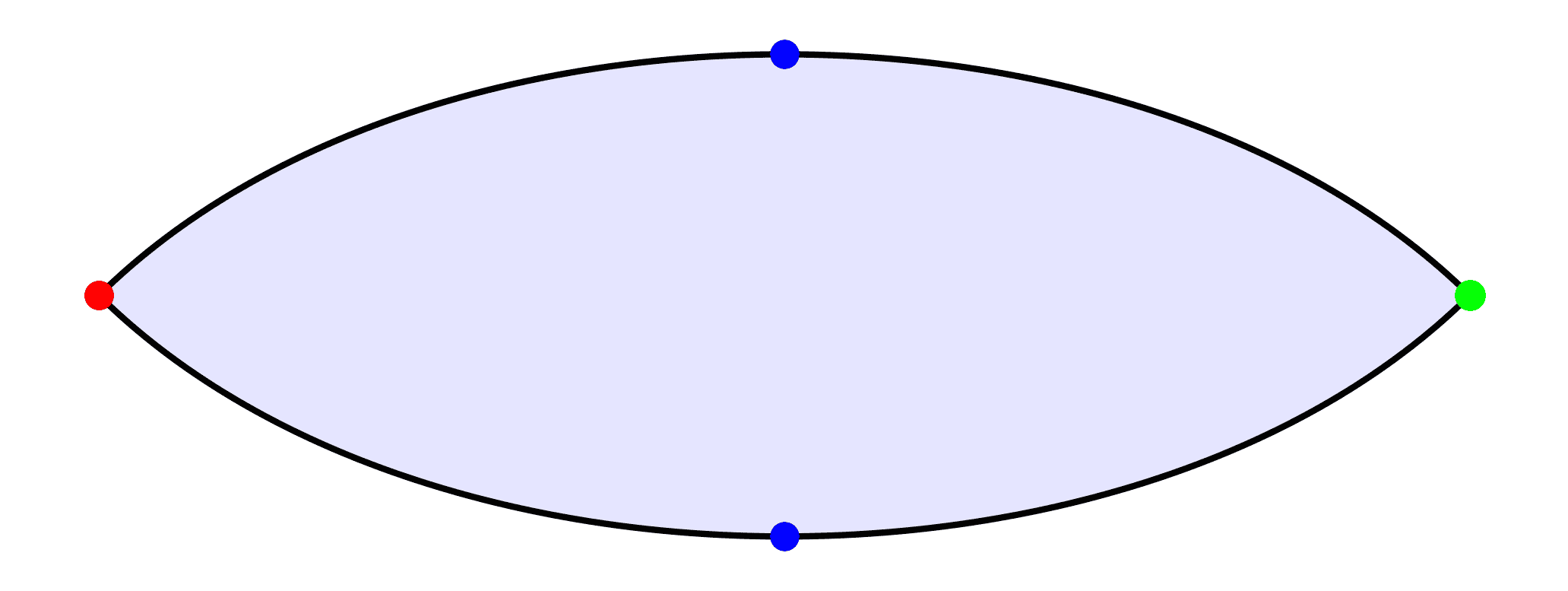}
        \includegraphics[align=c, height=3cm]{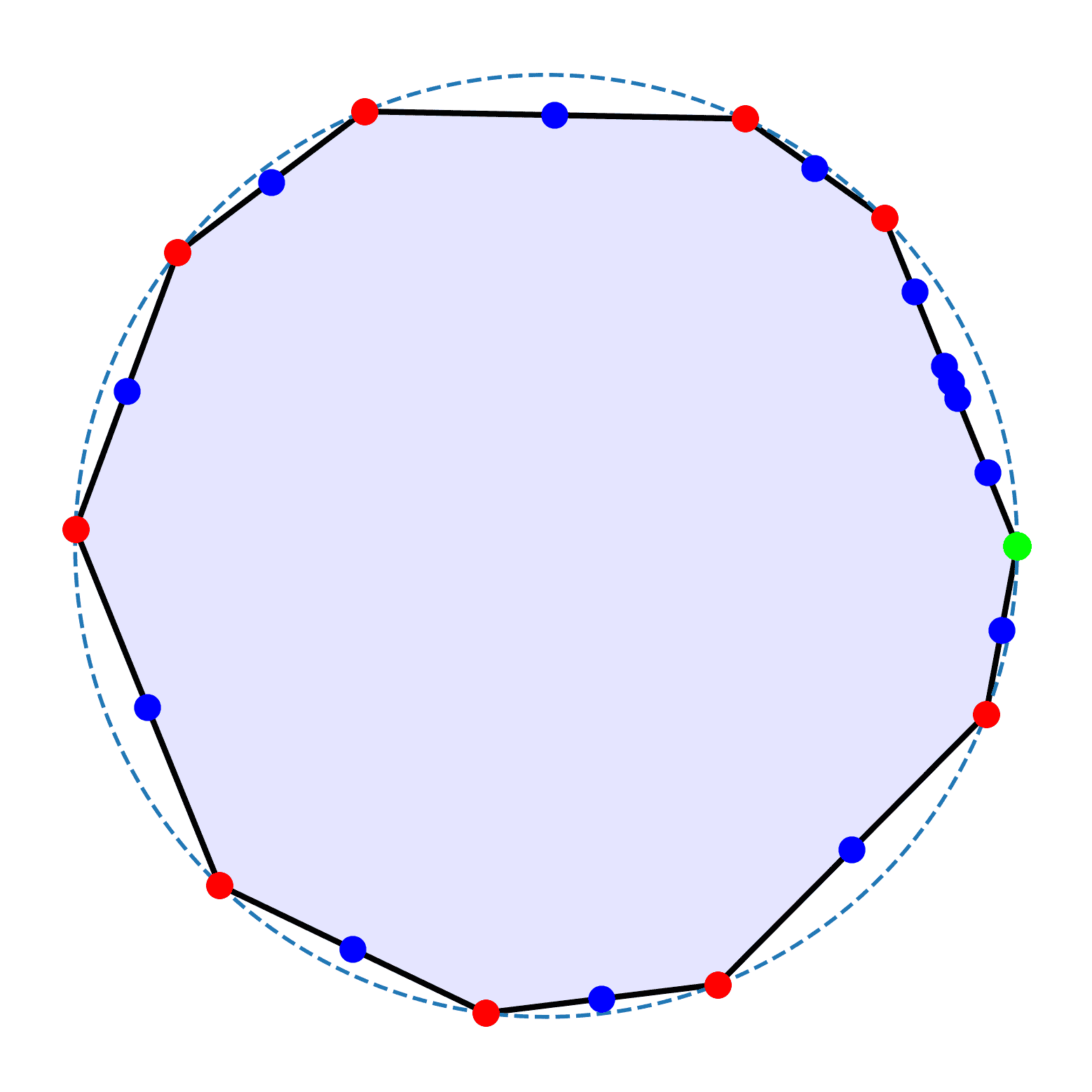}
    \end{subfigure}
\caption{Figure showing various domains $\Omega$ (top) along with the associated parametric domain $\hOm^{\mathbf{r}}$ we assign to each domain type (bottom). Dots indicate vertices that fence off the edges, wherein the root vertex is given in green while vertices that are modelled as a corner vertex are red and non-corner vertices are blue.}
\label{fig:controlmap_domains}
\end{figure}

\noindent For $N_{\text{convex}} = 0$, we take $\hOm^{\mathbf{r}}$ as the unit disc where we divide $\partial \hOm^{\mathbf{r}}$ into $\vert F \vert$ parts whose lengths match the relative lengths of the $p = w(e)$ with $e \in F$ in $\partial \Omega$. Hence, let $L(\partial \Om)$ denote the discrete length of the point set comprised of all vertices of $\partial \Omega$. If $p = w(e), e \in F$ has length $L(e)$, the corresponding segment of $\partial \hOm^{\mathbf{r}}$ has length $2 \pi L(e) / L(\partial \Om)$ (see Figure \ref{fig:controlmap_domains}, left). In this way, we ensure that the induced correspondence $\mathbf{f}^{\mathbf{r} \rightarrow \mathbf{x}}: \partial \hOm^{\mathbf{r}} \rightarrow \partial \Omega^S$ has a (nearly) continuous tangent along $\partial \hOm^{\mathbf{r}}$.
\begin{remark}
    While a harmonic map $\mathbf{x}: \hOm^{\mathbf{r}} \rightarrow \Omega^S$ can not be diffeomorphic under the correspondence $\mathbf{f}^{\mathbf{r} \rightarrow \mathbf{x}}: \partial \hOm^{\mathbf{r}} \rightarrow \partial \Omega^S$, in practice, the minuscule jumps in the tangent can be ignored for discrete approximations.
\end{remark}
For $N_{\text{convex}} = 1$, we employ the so-called teardrop domain (see Figure, second from left) which creates $\partial \hOm^{\mathbf{r}}$ using a cubic Hermite curve with specified tangent(s) while the starting and ending vertices are both placed in the origin. As before, the Hermite Curve is segmented based on the lengths associated with the $e \in F$. By default, the opening angle at the boundary's sole corner is chosen so as to match the (convex) angle at the sole corner of $\partial \Omega^S$. \\

\noindent For $N_{\text{convex}} = 2$, we take $\hOm^{\mathbf{r}}$ as the half disc or as a lens-shaped domain (Figure \ref{fig:controlmap_domains} third and fourth, respectively). The former is chosen when the two corner vertices $(v_\alpha, v_\beta)$ are adjacent in $\mathbb{V}(F)$ while the connecting edge is assigned to the half disc's straight segment. In both cases, the remaining breaks are segmented such that the tangent is approximately continuous. In the lens domain, we choose the two opening angles based on the corresponding interior angles in $F$. \\

\noindent Finally, for $N_{\text{convex}} \geq 3$, let the corner vertices be given by $V_{\text{corner}}(F) = \left(v_{i_1}, v_{i_2}, v_{i_3}, \ldots \right) \subset \mathbb{V}(F)$. The sequence $V_{\text{corner}}(F)$ naturally divides $F$ into $n \geq 3$ batches $b_1, \ldots, b_n$, with $b_j \subset F$ comprised of the edges encountered in the sequence from $v_{i_j} \in V_{\text{corner}}(F)$ to $v_{i_{j+1}} \in V_{\text{corner}}(F)$ (with periodic continuation if applicable). We draw an $n$-sided polygon $\partial \hOm^{\mathbf{r}}$ wherein the vertices are placed on the unit circle such that the lengths of the edges relative to $L(\partial \hOm^{\mathbf{r}})$ equal the lengths of the associated batches $b_j$ relative to $L(\partial \Omega)$. Then, if a batch $b_j$ is comprised of more than one edge, the corresponding side is again segmented based on the relative lengths of the contained edges (see Figure \ref{fig:controlmap_domains}, right). This ensures that the induced correspondence $\mathbf{f}^{\br \rightarrow \bx}: \partial \hOm^\br \rightarrow \partial \Om^S$, with $\partial \Om^S \approx \partial \Om$ the spline approximation, has a nearly continuous tangent on the straight segments of $\partial \Om^\br$. \\

\noindent The multipatch covering of $\hOm^{\mathbf{r}}$, which follows from the covering of $\hOm$ under a suitably constructed $\mathbf{r}: \hOm \rightarrow \hOm^{\mathbf{r}}$, is discussed in Section \ref{sect:patch_adjacency_graphs}.

Having discussed the various geometry types, this paper adopts a total of three different methods for the computation of an inversely harmonic map $\mathbf{x}^{\mathbf{r}}: \hOm^{\mathbf{r}} \rightarrow \Omega^S$, with $\Omega^S \approx \Omega$:
\begin{enumerate}
    \item The $C^0$-DG approach based on a quasi-linear elliptic PDE in nondivergence form \cite{hinz2024use}.
    \item The \textit{regularised weak-form discretisation} \cite{hinz2024use}.
    \item Floater's algorithm \cite{floater2005surface}.
\end{enumerate}

Approaches 1. and 2. create a spline-based parameterisation $\bx: \hOm^\br \rightarrow \Om^S$ by operating on the tuple $(\hOm^\br, \Om^S)$ directly. The discussion of Approach 1. is postponed to Section \ref{sect:parameterisation} while Approach 2. is encountered in Section \ref{sect:post_processing}. \\
Contrary to approaches 1. and 2., approach 3. does not operate on the spline input but on discrete point cloud approximations $(\partial \hOm^\br_h, \partial \Omega_h^S)$ sampled from $(\partial \hOm^\br, \partial \Omega^S)$. The algorithm then utilises a triangulation of $\Omega^S_h$ to provide a piecewise-linear map $\mathbf{x}_h^\br: \hOm^\br_h \rightarrow \Om^S_h$ whose inverse is approximately harmonic. Approach 3. leads to a linear problem that is computationally inexpensive compared to approaches 1. and 2. It finds applications as a surrogate for $\mathbf{x}^\br: \hOm^\br \rightarrow \Omega^S$ in Section \ref{sect:prep_temp} and is discussed in greater detail in \ref{sect:appendix_floater}. \\
While a spline map could be fit to the vertex pairs (potentially with regularisation terms) created by approach 3., we have found approaches 1. and 2. to be a more robust choice in practice because $\left( \mathbf{x}_h^{\mathbf{r}} \right)^{-1}: \Om_h^S \rightarrow \hOm^{\mathbf{r}}$ may map a homogeneous triangulation $\mathcal{T}_h$ to a very heterogeneous triangulation $\hat{\mathcal{T}}_h^r$, which often leads to an underdetermined problem in the absence of regularisation terms. 

\subsection{Creation of a Template Catalogue}
\label{sect:patch_adjacency_graphs}
\begin{figure}[h!]
\centering
\includegraphics[align=c, width=0.9 \textwidth]{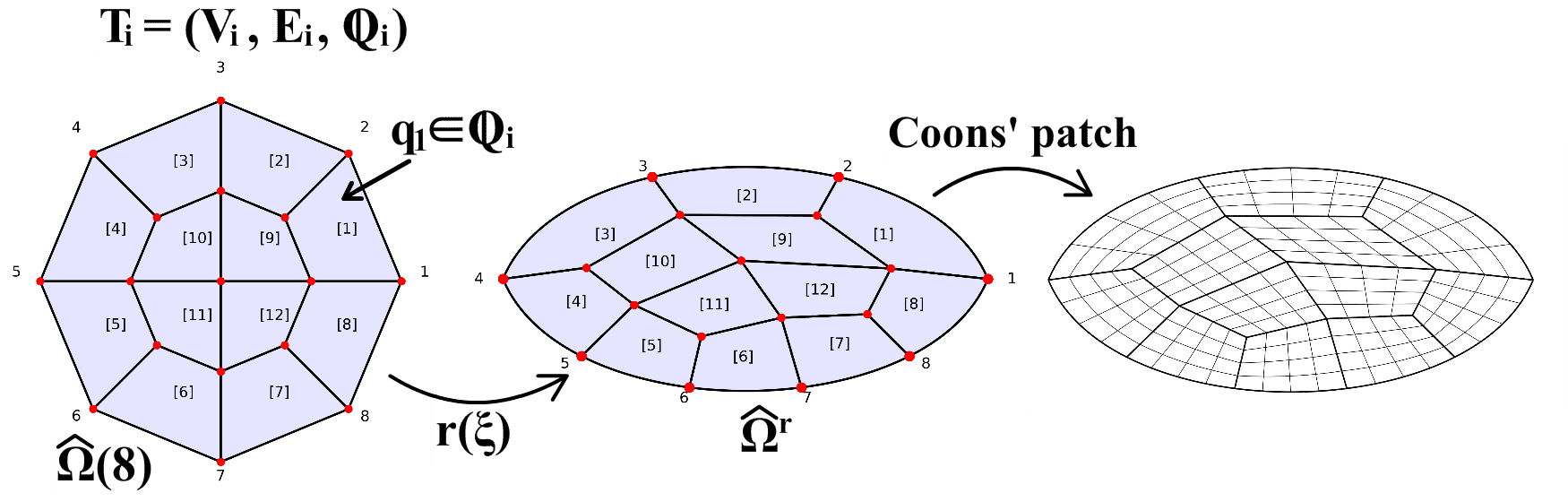}
\caption{Figure showing a template $\mathbb{T} \ni T_i = (V_i, E_i, \mathcal{Q}_i)$ for the $8$-sided regular polygon of radius one. The figure additionally shows the corresponding multipatch covering of a lens-shaped control domain $\hOm^\br$ (see Section \ref{sect:harmonic_maps}). The control domain is parameterised by the controlmap $\mathbf{r}: \hOm(8) \rightarrow \hOm^\br$ that is constructed using optimisation of the inner patch vertices along with the Coons' patch approach. For details, see \ref{sect:appendix_PAG}.}
\label{fig:PAG_to_controlmap}
\end{figure}

Approaches 1. and 2. from Section \ref{sect:harmonic_maps} operate on a multipatch spline space $\mathcal{V}_h \subset H^1(\hOm_i)$ defined over a quadrangulation of $\hOm_i$ which is itself represented by a template $\mathbb{T} \ni T_i = (V_i, E_i, \mathcal{Q}_i)$ with four-sided faces $q \in \mathcal{Q}_i$, see Figure \ref{fig:PAG_to_controlmap}.

\noindent Using the concept of patch adjacency graphs (PAGs), as introduced in \cite{buchegger2017planar}, for this paper we generated a template catalogue $\mathcal{T}$ comprised of $\approx 6 \times 10^4$ templates covering the regular $N$-sided polygon $\hOm(N)$ of radius one for $N \in \{4, 6, \ldots, 16\}$. We note that the computation of all PAGs along with the generation of a concrete template $T_i = (V_i, E_i, \mathcal{Q}_i)$ from the PAG is an \textit{offline} computation which does not contribute to the computational costs of what follows. \\
Since each domain $\Omega$ that follows from some $F \in \mathcal{F}$ is accompanied by a canonical control domain $\hOm^\br$ (see Section \ref{sect:harmonic_maps}), a template $T_i \in \mathbb{T}$ covering $\hOm(N)$ additionally requires the creation of a piecewise smooth (over the quadrangular faces) control map $\mathbf{r}: \hOm(N) \rightarrow \hOm^\br$. Unlike the creation of $\mathbb{T}$ this step constitutes an \textit{online} computation and is tackled via the computationally inexpensive minimisation of a functional. \\
The concept of PAGs, the creation of the catalogue $\mathbb{T}$ and the computation of a valid controlmap $\mathbf{r}: \hOm(N) \rightarrow \Om^\br$ for layouts $T_i \in \mathbb{T}$ is discussed in \ref{sect:appendix_PAG}.

\section{Graph Preparation and Templatisation}
\label{sect:prep_temp}
The goal of this stage is to pre-process $G = (V, E, \mathcal{F})$ so that each $F \in \mathcal{F}$ is transformed into a form attainable through parameterization based on harmonic maps. Since we utilise multipatch quad layouts to cover each $\hOm$ of face $F \in \mathcal{F}$, this means in particular that the graph's faces must be comprised of an even number of edges. Furthermore, the parameterisation strategy in stage 3. requires all of the face's corners to be convex. \\

\noindent In practice, the present stage is the most time consuming and still requires a substantial amount of manual intervention in many cases. However, we will present a number of automated strategies that we have found to be robust on a wide range of input graphs. Improving this stage's degree of automation will the topic of a forthcoming publication. \\
Ensuring that all faces $F \in \mathcal{F}$ are comprised of an even number of edges may require splitting select edges $e \in F$ while avoiding concave corners may necessitate creating a new edge that splits an existing face in two. \\
We note that the main difficulty associated with the splitting of a face's edge is the possibility of it, in turn, introducing an uneven number of edges in a neighbouring face. This includes faces $F \in \mathcal{F}$ to which a suitable template $T \in \mathbb{T}$ had previously been assigned. In order to handle this situation, we introduce the concept of template refinement in Section \ref{sect:templatisation}.

\subsection{Removing concave corners}
\label{sect:remove_concave_corners}
A harmonic map between a convex domain and a domain containing a concave corner will only be homeomorphic (see Section \ref{sect:harmonic_maps}). As such, the map may not be differentiable in a boundary point which we regard as undesirable. While optimisation-based parameterisation methods that can potentially handle concave corners in the multipatch setting have been studied in, for instance, \cite{buchegger2017planar}, they suffer from other drawbacks. More precisely, optimisation-based methods often require nonconvex optimisation that may converge to a local minimum or, more critically, fail to provide a nonsingular map even if the global minimiser over the finite-dimensional spline space is found. This is mainly due to optimisation-based approaches typically not possessing the same nonsingularity property (in the infinite-dimensional setting) of harmonic maps. As such, even the infinite-dimensional global optimiser may be folded. However, we emphasise that the techniques presented in this section are applicable if another parameterisation technique is adopted. If this technique is compatible with concave corners, this section's steps may be skipped. \\

\noindent To avoid concave corners, the first step of the preparation stage is repeatedly generating curves that originate in a concave vertex $v \in \mathbb{V}(F), F \in \mathcal{F}$ while splitting $F$ into two parts $F^+$ and $F^-$ that no longer contain the concave $v \in \mathbb{V}(F)$. In general, we assume that the curve connects two existing vertices $\{v_\alpha, v_\beta\} \subset \mathbb{V}(F_i)$ and is strictly contained in the interior of $\partial \Om = \cup w(e), e \in F$. We denote the newly-created edge between $v_{\alpha}$ and $v_{\beta}$ by $e_{\alpha \beta}$. \\

\begin{figure}[h!]
\centering
\includegraphics[align=c, width=0.6\linewidth]{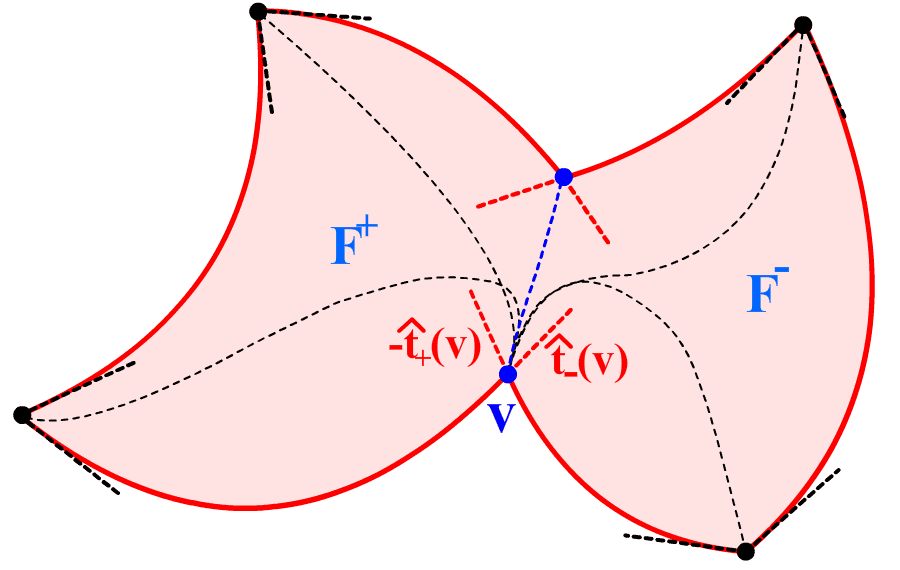}
\caption{An example of removing a concave vertex (blue) via the use of Hermite splitting curves whose starting and end point tangents are contained in the endpoints' convex cones. The figure shows all possible Hermite curves originating in the bottom concave vertex $v \in V$. While all curves are valid, selection based on~\eqref{eq:splitting_curve_quality_mu} would favour the curve highlighted in blue since it is the straightest and removes two concave vertices at once.}
\label{fig:convex_cone_splitting_curve}
\end{figure}

\noindent We call a curve $C: (0, 1) \rightarrow \mathbb{R}$, with $C \in C^1([0, 1], \mathbb{R}^2)$, a valid splitting curve between $v_{\alpha}$ and $v_{\beta}$ if it satisfies:
\begin{align}
\left \{ \begin{array}{l}
    C(0) = v_{\alpha}, \\
    C(1) = v_{\beta}, \\
    v_{\alpha} \text{ concave } \implies C^{\prime}(0) \in CC(\hat{t}_{-}(v_{\alpha}), -\hat{t}_{+}(v_\alpha)),  \\
    v_{\beta} \text{ concave } \implies -C^{\prime}(1) \in  CC(\hat{t}_{-}(v_{\beta}), -\hat{t}_{+}(v_\beta)),  \\
    C \vert_{(0, 1)} \subset \Omega.
    \end{array} \right.
\end{align}
Here, $CC(v, w)$ refers to the convex cone generated by two vectors $\{v, w\} \subset \mathbb{R}^2$, see Figure \ref{fig:convex_cone_splitting_curve}. \\

\noindent In the simplest case, a valid splitting curve is found by simply connecting two vertices linearly. However, we generally employ Hermite interpolation \cite[Chapter~5]{liseikin1999grid} to generate a (piecewise) polynomial curve with specified tangent at the vertices. Given the  vertex and tangent information $C(0) = v_\alpha$, $C^{\prime}(0)$, $C(1) = v_\beta$ and $C^{\prime}(1)$, the curve $C: [0, 1] \rightarrow \mathbb{R}^2$ is given by
\begin{align}
\label{eq:cubic_hermite_curve}
    C(t) = (2 t^3 - 3 t^2 + 1) v_\alpha + (t^3 - 2t^2 + t) C^{\prime}(0) + (-2 t^3 + 3 t^2) v_\beta + (t^3 - t^2) C^{\prime}(1).
\end{align}

\noindent When employing a method that allows for enforcing tangent information, we typically require it to satisfy 
\begin{align}
\label{eq:cubic_hermite_tangents}
    C^{\prime}(0) = -\eta \hat{n}(v_{\alpha}) \quad \text{and} \quad C^{\prime}(1) = \gamma \hat{n}(v_{\beta}),
\end{align}
with $\eta, \gamma > 0$ and $\hat{n}(\, \cdot \,)$ as defined in Section \ref{sect:Discrete_Tangent}. Heuristically, a good choice of scaling parameters is given by $\eta = \gamma = \|v_{\beta} - v_{\alpha} \|$. We flag a vertex $v_\alpha \in \mathbb{V}(F)$ as concave whenever the discrete interior angle $\angle(v_\alpha)$ satisfies $\angle(v_{\alpha}) \geq \pi + \varepsilon_{\angle}$, where $\varepsilon_{\angle} \ll 1$ is a user-specified parameter. \\
In practice, a new point set $p_{\alpha \beta} := w(e_{\alpha \beta})$ is formed by evaluating $C: [0, 1] \rightarrow \overline{\Omega}$ in a large number of discrete values $\{t_j\}_j$ uniformly-spaced over the unit interval. The validity of the discrete splitting curve follows from the requirement that all interior points of $p_{\alpha \beta}$ be strictly contained in $\Omega$, which can be verified using standard polygon routines \cite{shimrat1962algorithm}. \\

\noindent The new edge $e_{\alpha \beta}$ and the associated $p_{\alpha \beta}$ are added to the domain and range, respectively, of $w(\, \cdot \,)$. Then, we update the following components of $G$:
\begin{align}
    E & \leftarrow E \cup \{e_{\alpha \beta} \} \nonumber \\
    \mathcal{F} & \leftarrow \left(\mathcal{F} \setminus \{F\} \right) \cup \{F^+, F^- \}.
\end{align}
Letting $F = \left(e_{i_1}, \ldots, e_{\alpha}^{-}, e_{\alpha}^{+}, \ldots, e_{\beta}^{-}, e_\beta^{+}, \ldots, e_{i_{N_i}} \right)$, where $\{e_{\alpha}^{-}, e_{\alpha}^{+} \}$ and $\{e_{\beta}^{-}, e_{\beta}^{+} \}$ are the two positively-oriented edges incident to $v_\alpha$ and $v_\beta$, respectively, we have
\begin{align}
    F^+ := \left( e_{i_1}, \ldots, e_{\alpha}^{-}, e_{\alpha \beta}, e_{\beta}^{+}, \ldots, e_{i_{N_i}} \right) \quad \text{and} \quad F^- := \left( e_\alpha^+, \ldots, e_{\beta}^{-}, -e_{\beta \alpha} \right).
\end{align}
In what follows, we present a strategy aimed at the automated removal of concave vertices. In rare cases however, the methodology may fail in which case a splitting curve is generated by hand, thus requiring a modest amount of manual intervention in practice. \\

\noindent We introduce the operator $\mathbb{V}^{\varepsilon_{\angle}}_{\text{conc}}(\, \cdot \,)$ which, when applied to an $F \in \mathcal{F}$, returns all $v_\alpha \in \mathbb{V}(F)$ that have been flagged concave (according to the user-specified parameter $\varepsilon_{\angle} \geq 0$). We introduce the set $\mathcal{F}_{\text{conc}}^{\varepsilon_{\angle}} := \{F \in \mathcal{F} \quad \vert \quad \vert \mathbb{V}^{\varepsilon_{\angle}}_{\text{conc}}(F) \vert > 0 \}$ as the set of faces containing at least one concave corner. The default strategy is based on a greedy algorithm that iterates over the $F \in \mathcal{F}_{\text{conc}}^{\varepsilon_{\angle}}$ in an outer loop and selects pairs of vertices $\{v_\alpha, v_{\beta}\} \subset \mathbb{V}(F)$, with $\vert \{v_\alpha, v_\beta \} \cap \mathbb{V}^{\varepsilon_{\angle}}_{\text{conc}}(F) \vert \in \{1, 2\}$, as candidates for splitting. \\
Given some $F \in \mathcal{F}^{\varepsilon_\angle}_{\text{conc}}$, let $E_{\text{valid}}$ be the set of edges $e_{\alpha \beta}$ between $\{v_\alpha, v_\beta \} \subset \mathbb{V}(F)$, with $|\{v_\alpha, v_\beta\} \cap \mathbb{V}^{\varepsilon_\angle}_{\text{conc}}(F)| \in \{1, 2\}$ that result in a valid splitting curve $C_{\alpha \beta}$ (sampled from~\eqref{eq:cubic_hermite_curve}). We measure the quality of $C_{\alpha \beta}$ using the cost function
\begin{align}
\label{eq:splitting_curve_quality}
    \mathcal{Q}(C_{\alpha \beta})^2 = \frac{1}{\vert C_{\alpha \beta} \vert^2} \int \limits_{[0, 1]} \left \| \frac{\partial^2 C_{\alpha \beta}(t)}{\partial t^2} \right \|^2 \mathrm{d}t,
\end{align}
where $\vert C_{\alpha \beta} \vert$ denotes the curve's total length and smaller values imply a higher quality. We favour target vertices $v_{\beta}$ that are themselves concave since this removes two concave vertices instead of one. As such, we introduce a parameter $\mu \geq 1$ and define the adjusted quality cost function
\begin{align}
\label{eq:splitting_curve_quality_mu}
    \mathcal{Q}^{\mu}(C_{\alpha \beta}) := \left \{ \begin{array}{ll} \mathcal{Q}(C_{\alpha \beta}) & v_{\beta} \notin V^{\varepsilon_{\angle}}_{\text{conc}}(F) \\
    \frac{1}{\mu} \mathcal{Q}(C_{\alpha \beta}) & \text{else} \end{array} \right..
\end{align}
Heuristically, a good choice is $\mu = 2$. The greedy strategy now splits $F \in \mathcal{F}^{\varepsilon_\angle}_{\text{conc}}$ using the edge $e_{\alpha \beta} \in E_{\text{valid}}$ that minimises~\eqref{eq:splitting_curve_quality_mu}. The graph and weight function $w(\, \cdot \,)$ are updated to reflect the changes and the greedy strategy is called on the updated graph until no more concave vertices that are eligible for splitting are found. Should algorithm fail to remove a concave corner, a valid splitting curve is created manually. \\
The greedy strategy is summarised in Algorithm \ref{algo:concave_corners}.

\begin{algorithm}
	\caption{Remove concave corners and return the set of faces that require manual intervention}
	\label{algo:concave_corners}
	\small{
	\begin{algorithmic}[1]
		
		\Procedure{Remove concave corners}{Graph $G = (V, E, \mathcal{F})$, discrete concavity threshold $\varepsilon_{\angle} \geq 0$, favour concave targets parameter $\mu \geq 1$}
		\Statex
            \State $\mathcal{F}_{\text{manual}} \leftarrow \emptyset$
		\While{$\mathcal{F}_{\text{conc}}^{\varepsilon_{\angle}}(\mathcal{F}) \setminus \mathcal{F}_{\text{manual}} \neq \emptyset$}
		      \State $F \leftarrow$ first element of $\mathcal{F}_{\text{conc}}^{\varepsilon_{\angle}}(\mathcal{F}) \setminus \mathcal{F}_{\text{manual}}$
                \State $E_{\text{valid}} \leftarrow$ the set of $e_{\alpha \beta}$ ending in at least one $v \in \mathcal{V}_{\text{conc}}^{\varepsilon_{\angle}}(F)$ while creating a valid splitting curve $C_{\alpha \beta}$
                \If{$E_{\text{valid}} = \emptyset$}
                    \State Update $\mathcal{F}_{\text{manual}} \leftarrow \mathcal{F}_{\text{manual}} \cup \{F \}$
                \Else
                    \State $e \leftarrow e_{\alpha \beta} \in E_{\text{valid}}$ such that $C_{\alpha \beta}$ minimises~\eqref{eq:splitting_curve_quality_mu}
                    \State Update $\mathcal{F} \leftarrow \mathcal{F} \setminus F \cap \{F^+, F^-\}$, $E \leftarrow E \cup e_{\alpha \beta}$ and $w(\, \cdot \,)$
                \EndIf
		\EndWhile
            \Return{$\mathcal{F}_{\text{manual}}$}
		\EndProcedure{}
		
	\end{algorithmic}
	}
\end{algorithm}

\begin{figure}[h!]
\centering
\captionsetup[subfigure]{labelformat=empty}
    \begin{subfigure}[t]{0.31\textwidth}
        \centering
        \includegraphics[align=c, width=0.95\linewidth]{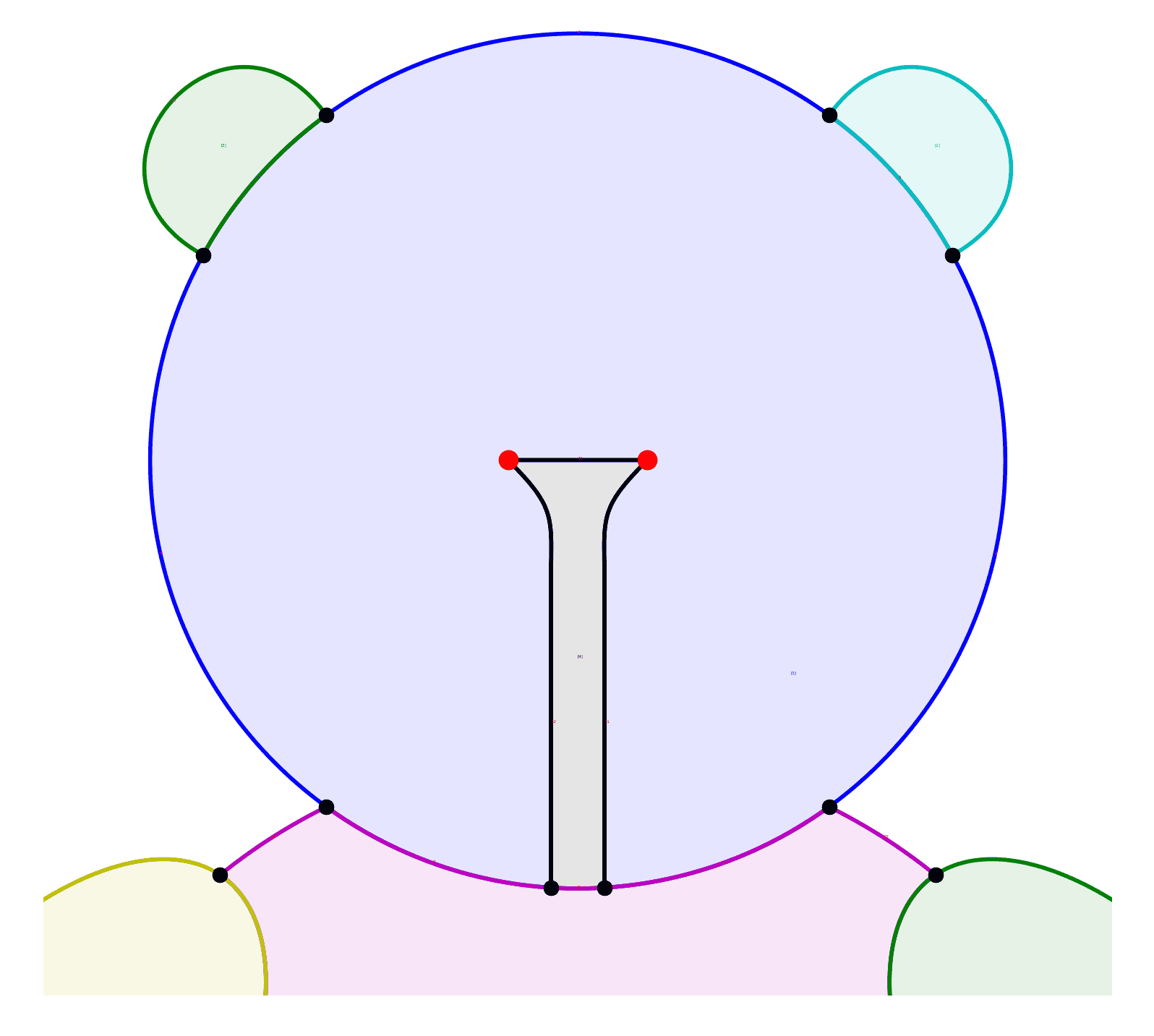}
    \end{subfigure} $\,$
    \begin{subfigure}[t]{0.31\textwidth}
        \centering
        \includegraphics[align=c, width=0.95\linewidth]{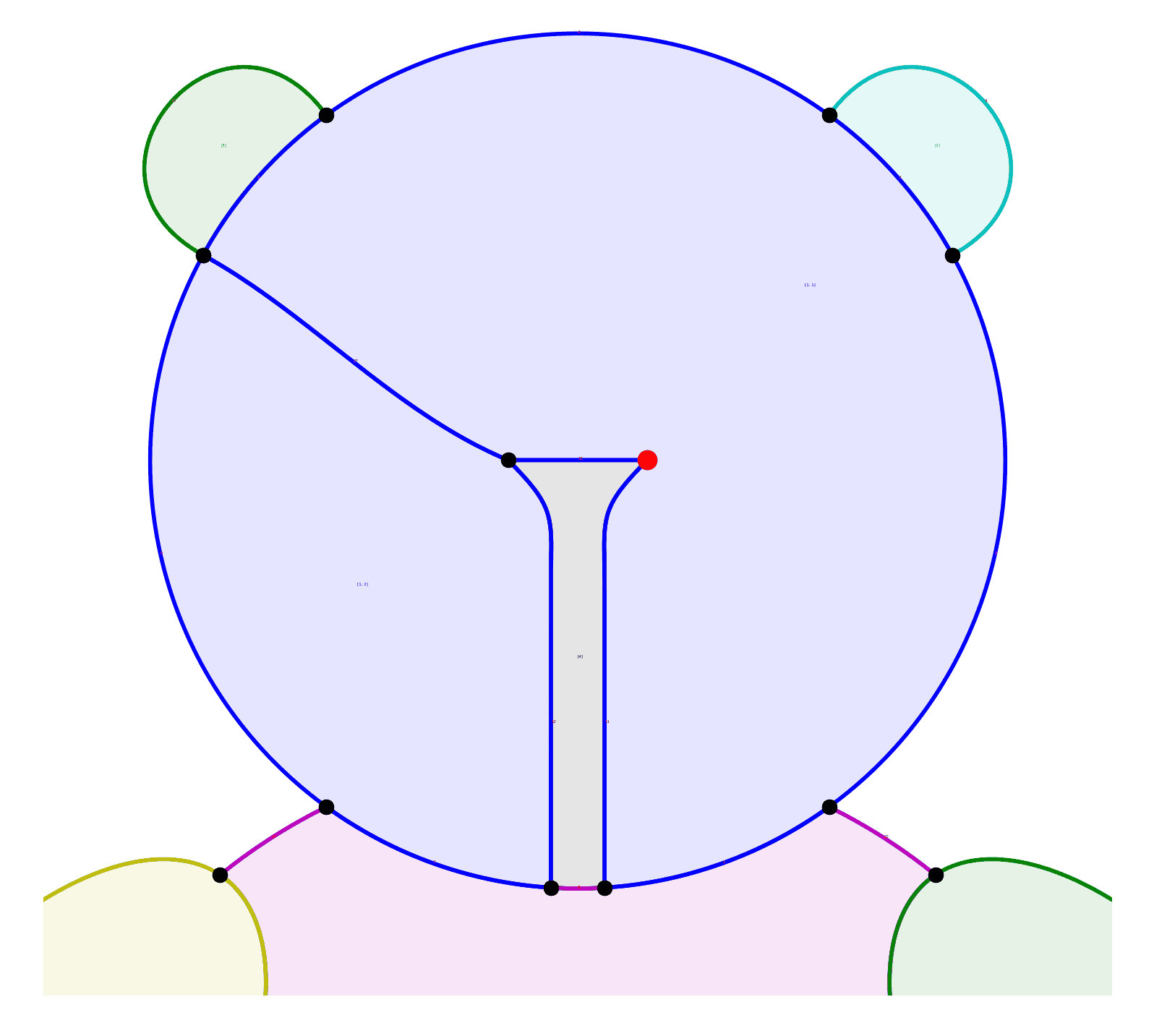}
    \end{subfigure}
    \begin{subfigure}[t]{0.31\textwidth}
        \centering
        \includegraphics[align=c, width=0.95\linewidth]{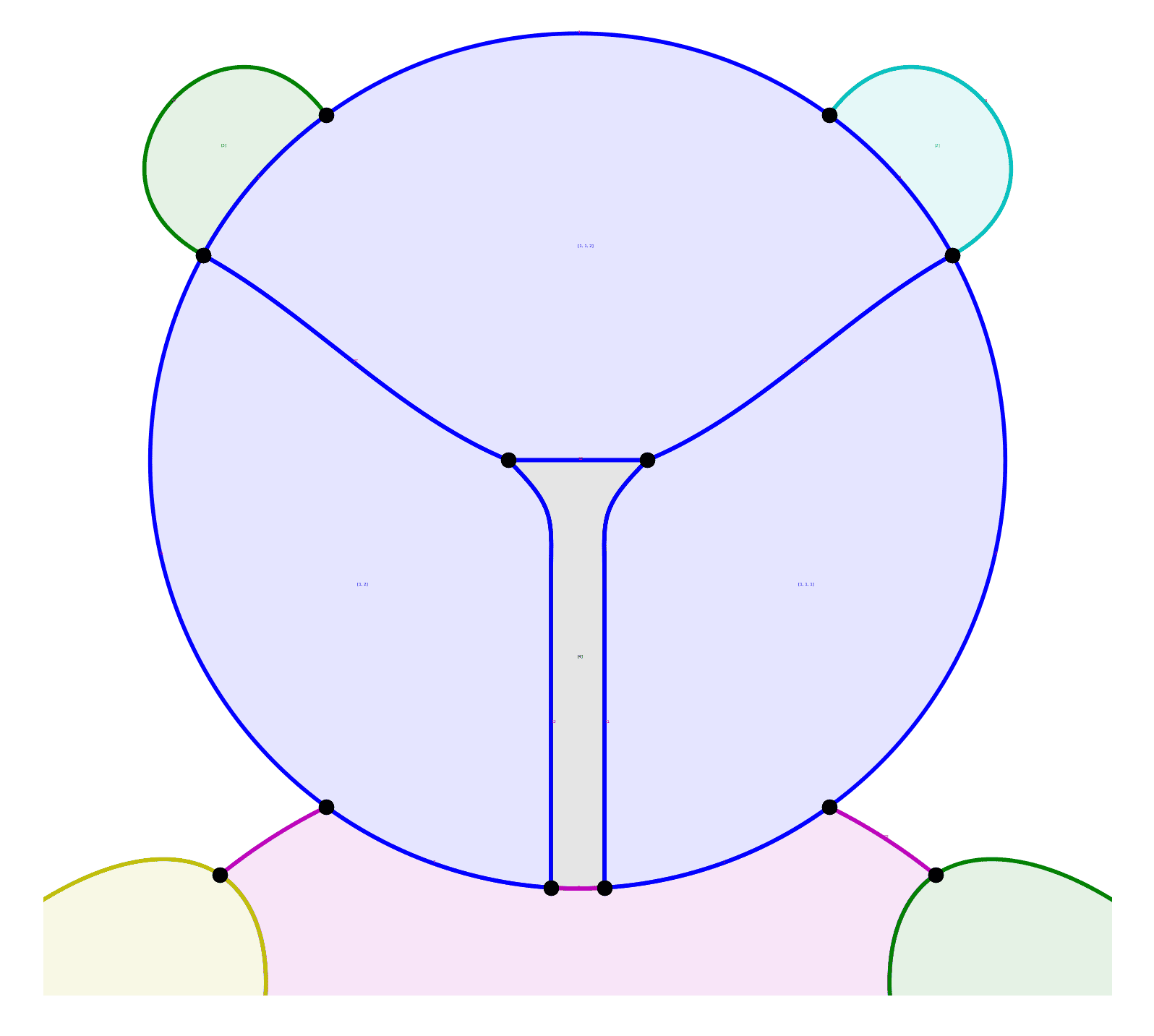}
    \end{subfigure}
\caption{An example showing the splitting of a face with concave corners (blue) by drawing a Hermite curve from a concave corner (red) to another vertex (red or black). Algorithm \ref{algo:concave_corners} selects vertex pairs $\{ v_{\alpha}, v_{\beta} \} \subset \mathbb{V}(F)$, where $v_\alpha$ is a concave vertex in a greedy fashion and connects them using a Hermite curve, thus splitting the face in two. The same routine is then applied to the updated plane graph until no more eligible vertex pairs are found. In this example, the routine succeeds in fully automatically removing all concave vertices.}
\label{fig:concave_corner_removal}
\end{figure}
Figure \ref{fig:concave_corner_removal} shows an example of Algorithm \ref{algo:concave_corners} applied to a geometry that contains a face with two concave vertices. The algorithm succeeds in removing both concave corners without manual intervention.

\subsection{Templatisation}
\label{sect:templatisation}
Upon removal of all concave corners, we are in the position to assign a polygonal parametric domain to each face $F \in \mathcal{F}$. Since we are opting for a quad layout parameterisation of the plane graph, each face $F \in \mathcal{F}$ has to be comprised of an even number of edges. While creating a convex parametric domain that is covered by a curvilinear quad layout is possible for $\vert F \vert = 2$ (see \cite{buchegger2017planar}), such a layout is not homeomorphic to a convex polygon with the same number of sides. As such, for templating faces $F \in \mathcal{F}$, we require $\vert F \vert \geq 4$. \\

The template that is assigned to a face $F_i \in \mathcal{F}$ with $\vert F_i \vert \geq 4$ and $\vert F_i \vert$ even is denoted by $T_i = (V_i, E_i, \mathcal{Q}_i)$. The $T_i \in \mathbb{T}$ represent quadrangulations of the regular polygon $\hOm(N)$ with $N = \vert \partial E_i \vert$. The pairing between $T_i$ and $F_i$ is accompanied by the bijective edge correspondence $\phi_i: \partial E_i \rightarrow F_i$, where we utilise the same operator $\phi_i$ to represent the correspondence between the $\mathbb{V}(\partial E_i) \ni \hat{v} \rightarrow v \in \mathbb{V}(F_i)$. 
For convenience, we require that $\phi_i: \partial E_i \rightarrow F_i$ pair edges in a cyclic fashion starting from the root edge in both $\partial E_i$ and $F_i$. To accommodate choosing other starting edge pairs, we may simply renumber the template's boundary vertices. \\
In what follows, we denote the set of templates $T_i$ that have already been assigned to some $F_i \in \mathcal{F}$ by $\mathcal{T}$. As such, we have $G = (V, E, \mathcal{F}, \mathcal{T})$ with $\vert \mathcal{T} \vert \leq \vert \mathcal{F} \vert$. \\

\noindent An automated templatisation strategy requires three main ingredients:
\begin{enumerate}[label=\textbf{\Alph*)}]
    \item A strategy to decide the order in which the faces $F \in \mathcal{F}$ are treated.
    \item A strategy for selecting an appropriate edge $e \in F$ for splitting in case $\vert F \vert$ is odd.
    \item An automated strategy for selecting an appropriate template $T \in \mathbb{T}$ for an even-sided face $F \in \mathcal{F}$.
\end{enumerate}

For \textbf{B)}, we note that splitting a face's edge $e \in F_i$ for the purpose of creating an even-sided face may, in turn, introduce an odd number of edges on a neighbouring face $F_j \in \mathcal{F}$ that shares the edge $e \in F_i$. More critically, if a template had already been assigned to another face $F_j \in \mathcal{F}$, it will no longer be compatible with the updated graph. To avoid breaking the template structure, this section introduces the concept of \textit{template refinement}. \\

\noindent For \textbf{C)}, this section introduces several strategies for selecting a suitable template from $\mathbb{T}$. By default, if $T \in \mathbb{T}$, the collection $\mathbb{T}$ also contains all graphs that result from a cyclic renumbering of the boundary vertices of $T$. Clearly, the optimal choice is application-dependent. \\

\begin{figure}[h!]
\centering
\includegraphics[align=c, width=0.95\linewidth]{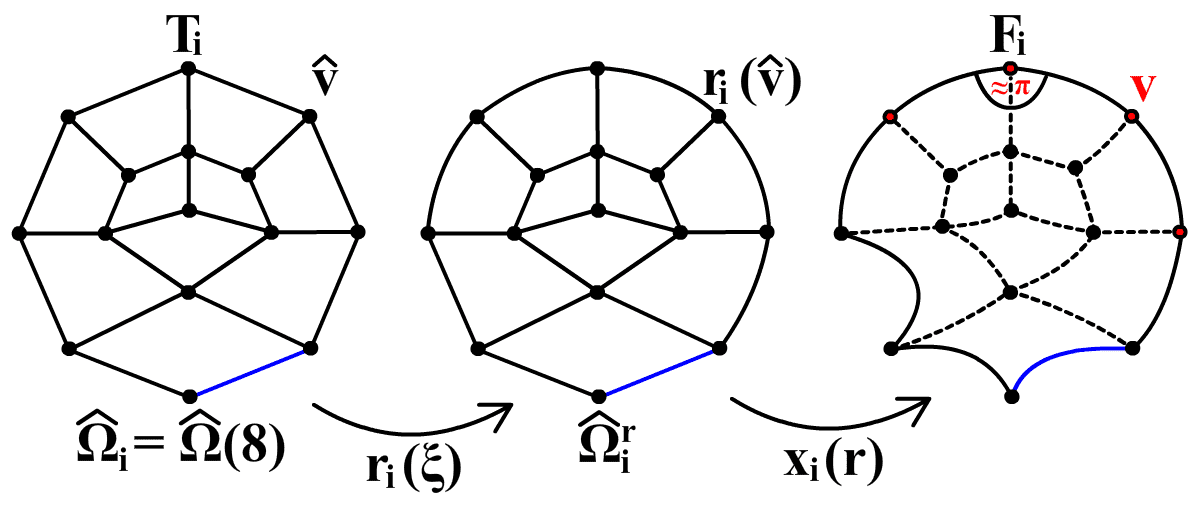}
\caption{Figure showing a geometry (right) with several vertices (highlighted in red) that create an interior angle close to $\theta = \pi$. We retain only templates (left) that assign to vertices $v \in \mathbb{V}(F_i)$ with $\angle(v) \approx \pi$ some $\hat{v} \in V_i$ with $\operatorname{val}(\hat{v}) \geq 3$. The figure additionally shows the convex control domain $\widehat{\Omega}_i^{\mathbf{r}}$ that is compatible with the corner pattern of $\Omega_i$ (see Section \ref{sect:harmonic_maps}).}
\label{fig:valence}
\end{figure}

Apart from the obvious filtering step that keeps only the templates $T_i \in \mathbb{T}$ with $\vert \partial E_i \vert = \vert F_i \vert$, a second filter retains only templates that are compatible with the corner pattern of $F_i \in \mathcal{F}$. More precisely, if $\vert \angle(v) - \pi \vert $, with $v \in \mathbb{V}(F_i)$, is small in $F_i$, we retain only templates $T_i \in \mathbb{T}$ that result in a correspondence wherein the vertex $\phi^{-1}(v) = \hat{v} \in \mathbb{V}(\partial E_i)$ has valence $\operatorname{val}(\hat{v}) \geq 3$ to avoid (nearly) singular corners, see Figure \ref{fig:valence}. \\

\begin{figure}[h!]
\centering
    \includegraphics[width=0.2\linewidth, valign=c]{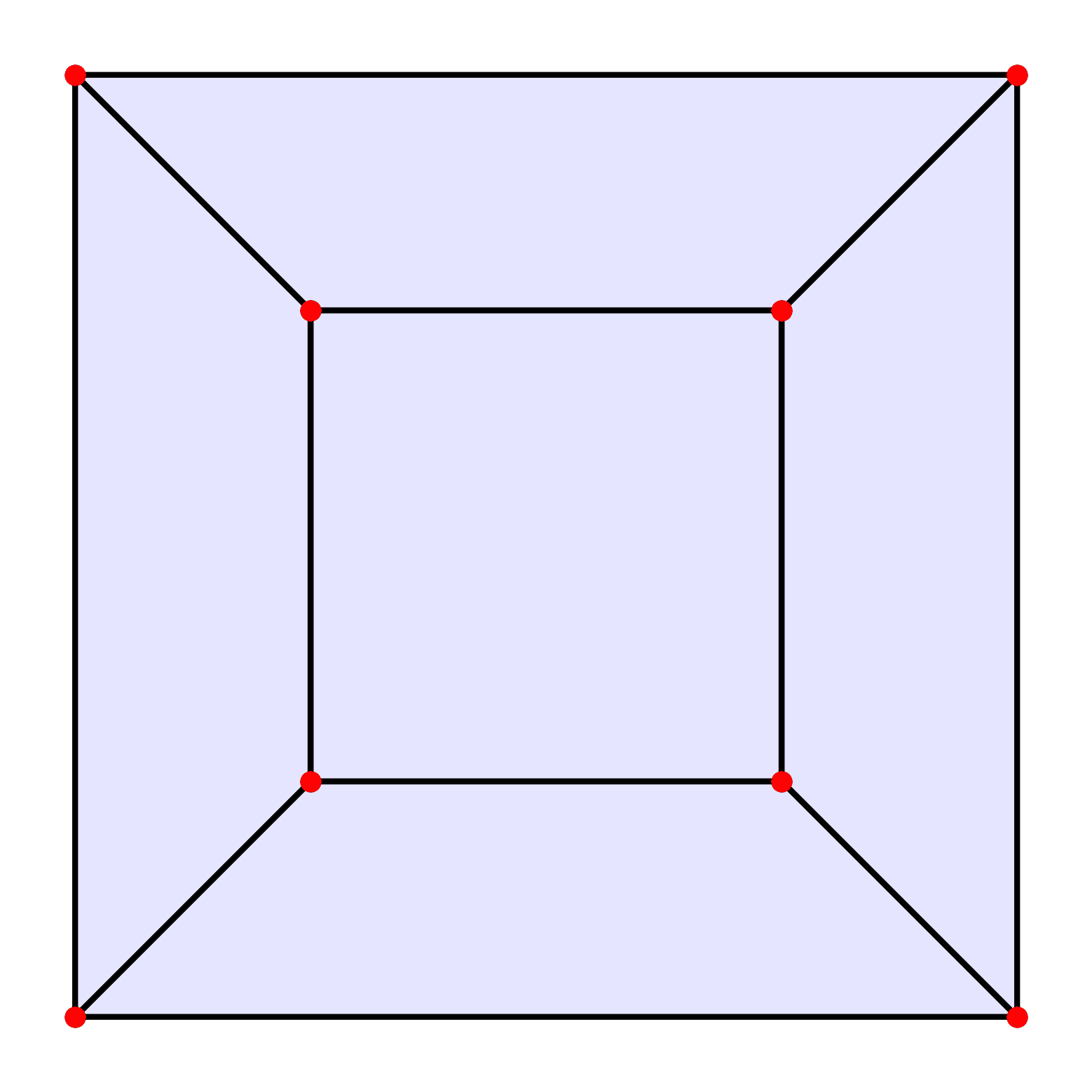}
    \includegraphics[width=0.2\linewidth, valign=c]{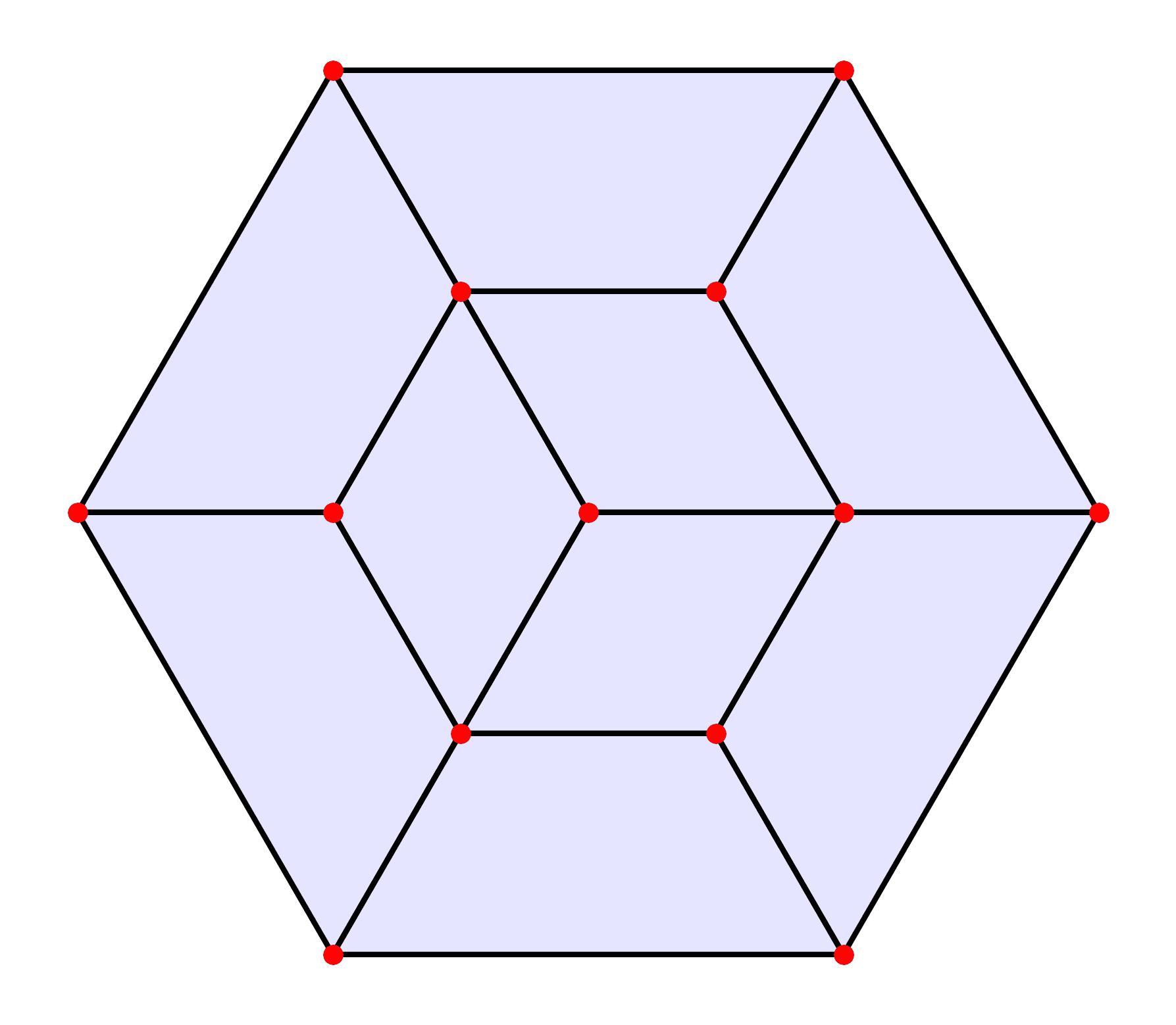}
    \includegraphics[width=0.2\linewidth, valign=c]{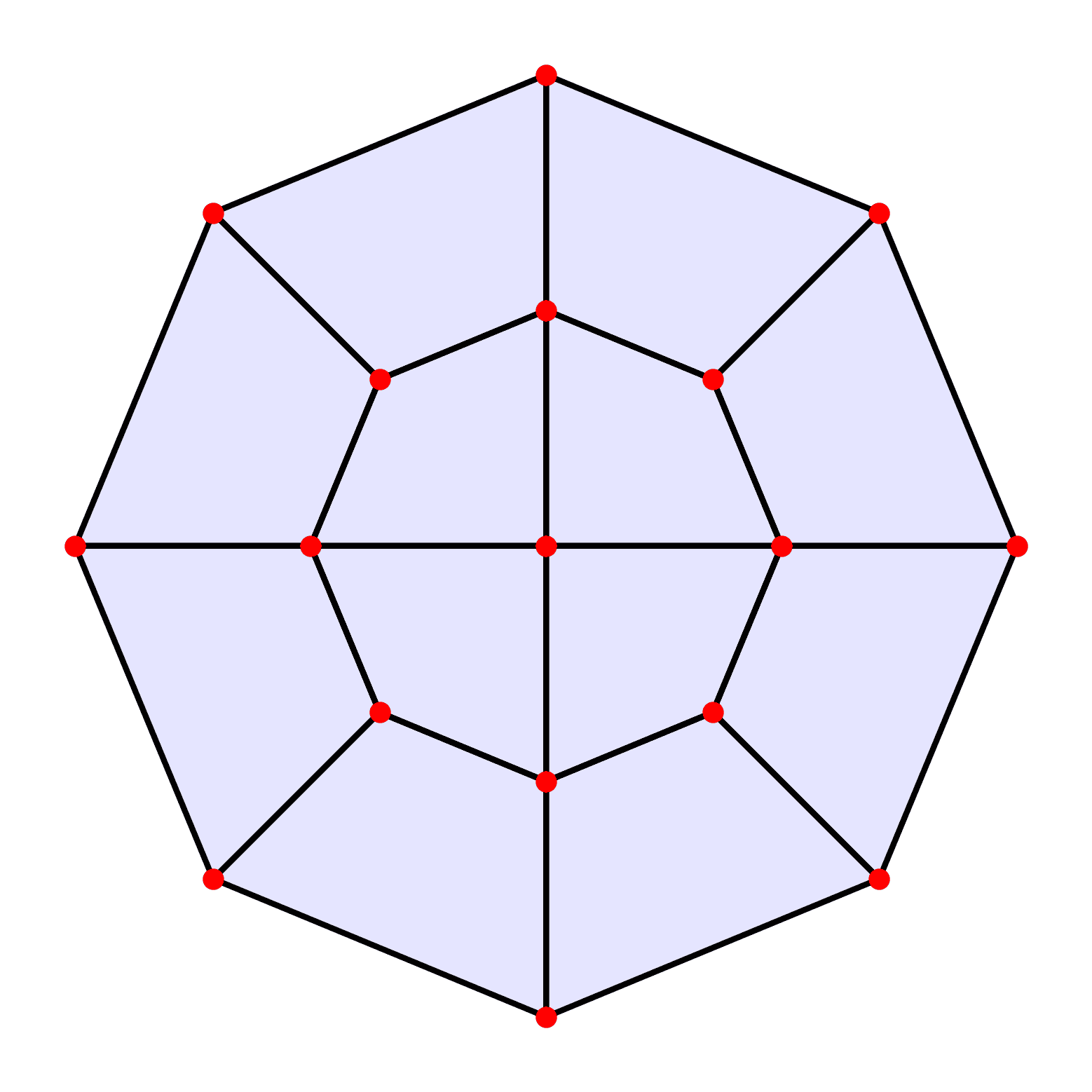}
    \includegraphics[width=0.2\linewidth, valign=c]{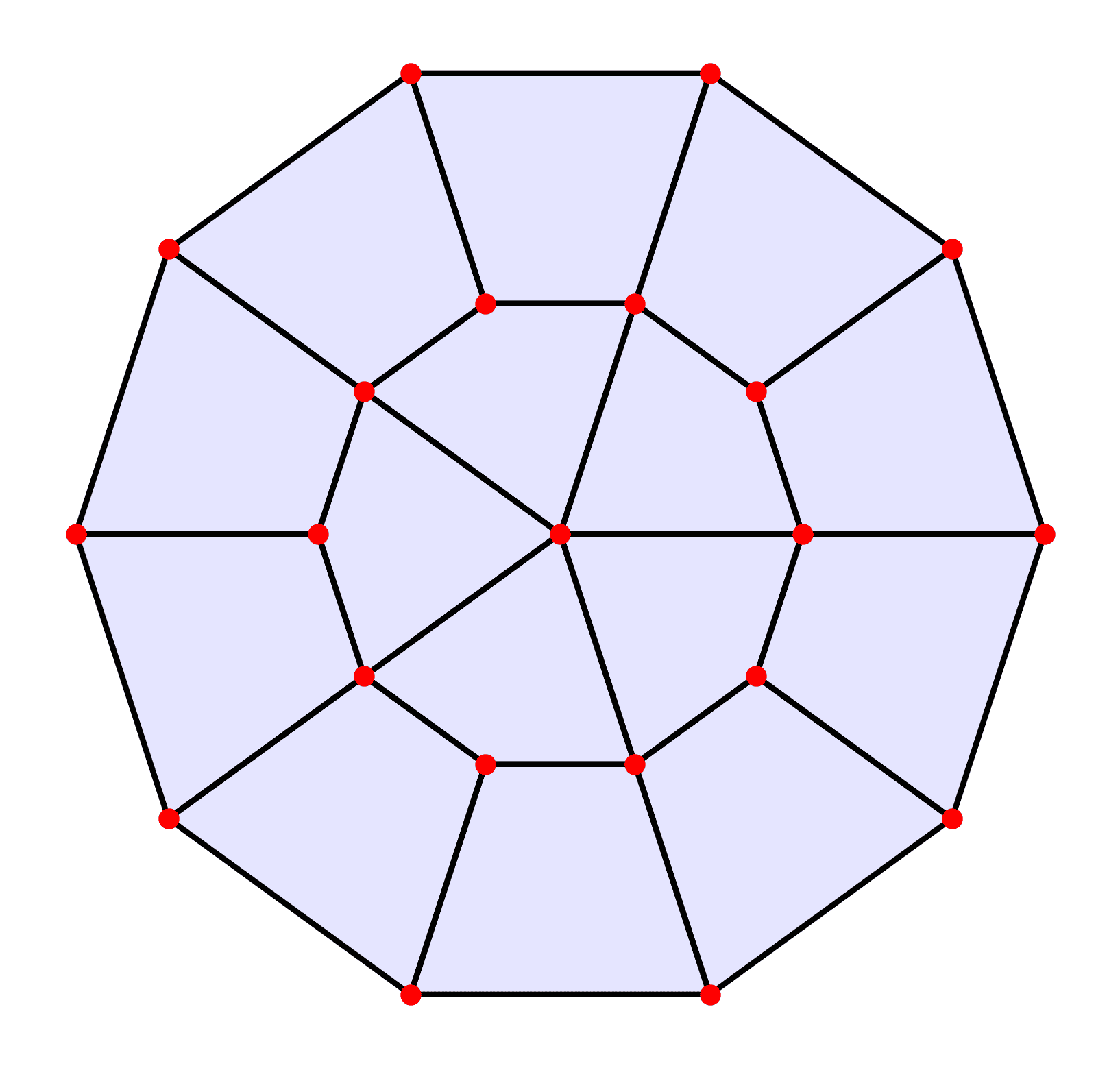}
\caption{Figure showing the radial-centric $N$-leaf template for $N \in \{4, 6, 8, 10\}$.}
\label{fig:trampoline}
\end{figure}

\noindent In rare cases (and depending on the template selection strategy), in step $\textbf{C)}$, the ordered set of eligible templates may be empty. This is in particular true for $N$-sided faces with $N > 16$. In this case, by default, we fall back on the so-called \textit{radial-centric (RC) $N$-leaf} template. Figure \ref{fig:trampoline} shows the RC $N$-leaf template for $N \in \{4, 6, 8, 10 \}$. For $N \geq 12$ the RC template continues the pattern observed for $N \in \{6, 8, 10\}$.  
We note that The RC $N$-leaf template is compatible with all vertex angle patters since each boundary vertex $\hat{v} \in \mathbb{V}(\partial E_i)$ satisfies $\operatorname{val}(\hat{v}) = 3$. \\



\begin{figure}[h!]
\centering
\captionsetup[subfigure]{labelformat=empty}
    \begin{subfigure}[t]{0.3\textwidth}
        \centering
        \includegraphics[align=c, width=0.95\linewidth]{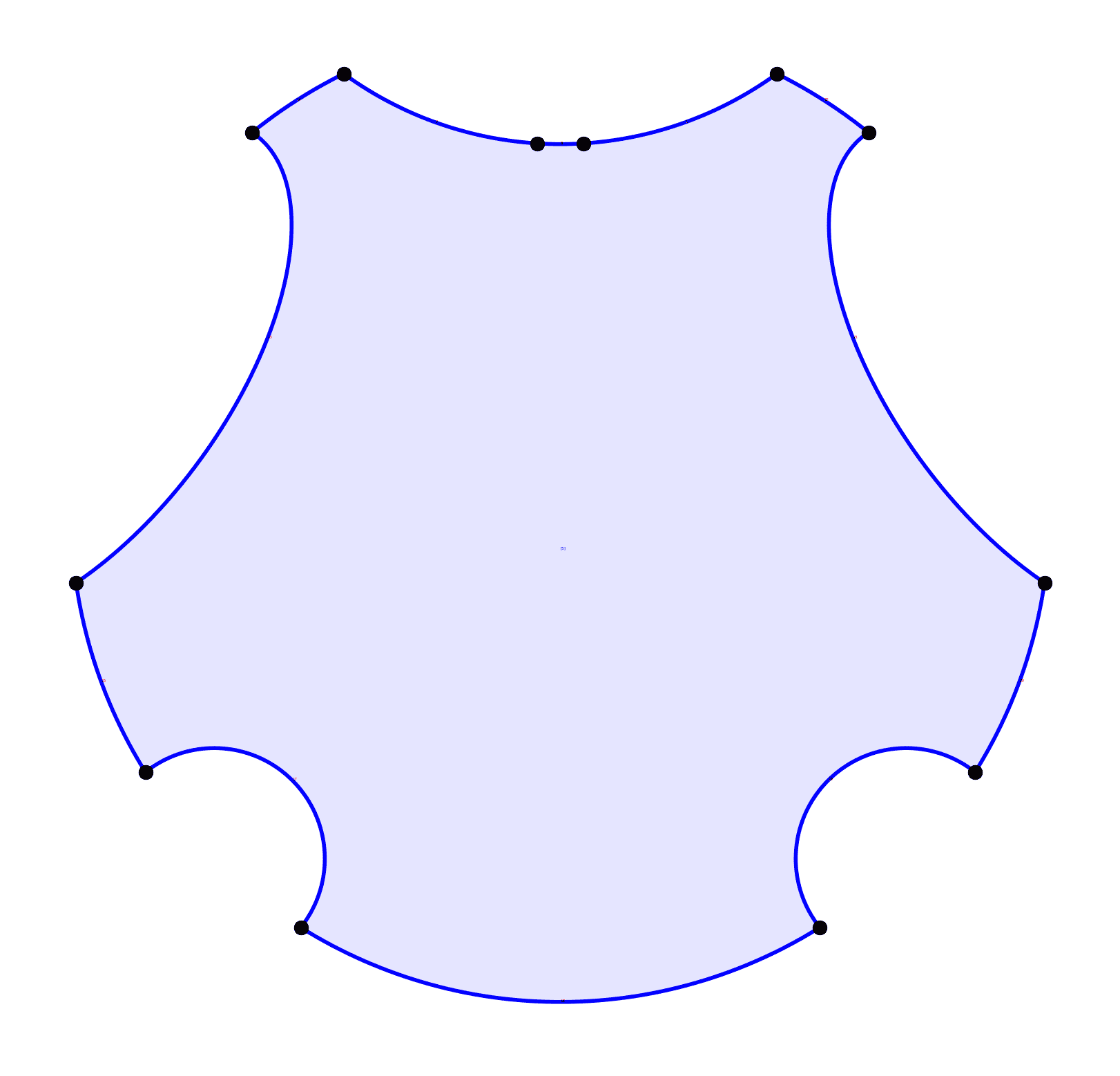}
    \end{subfigure} $\,$
    \begin{subfigure}[t]{0.3\textwidth}
        \centering
        \includegraphics[align=c, width=0.95\linewidth]{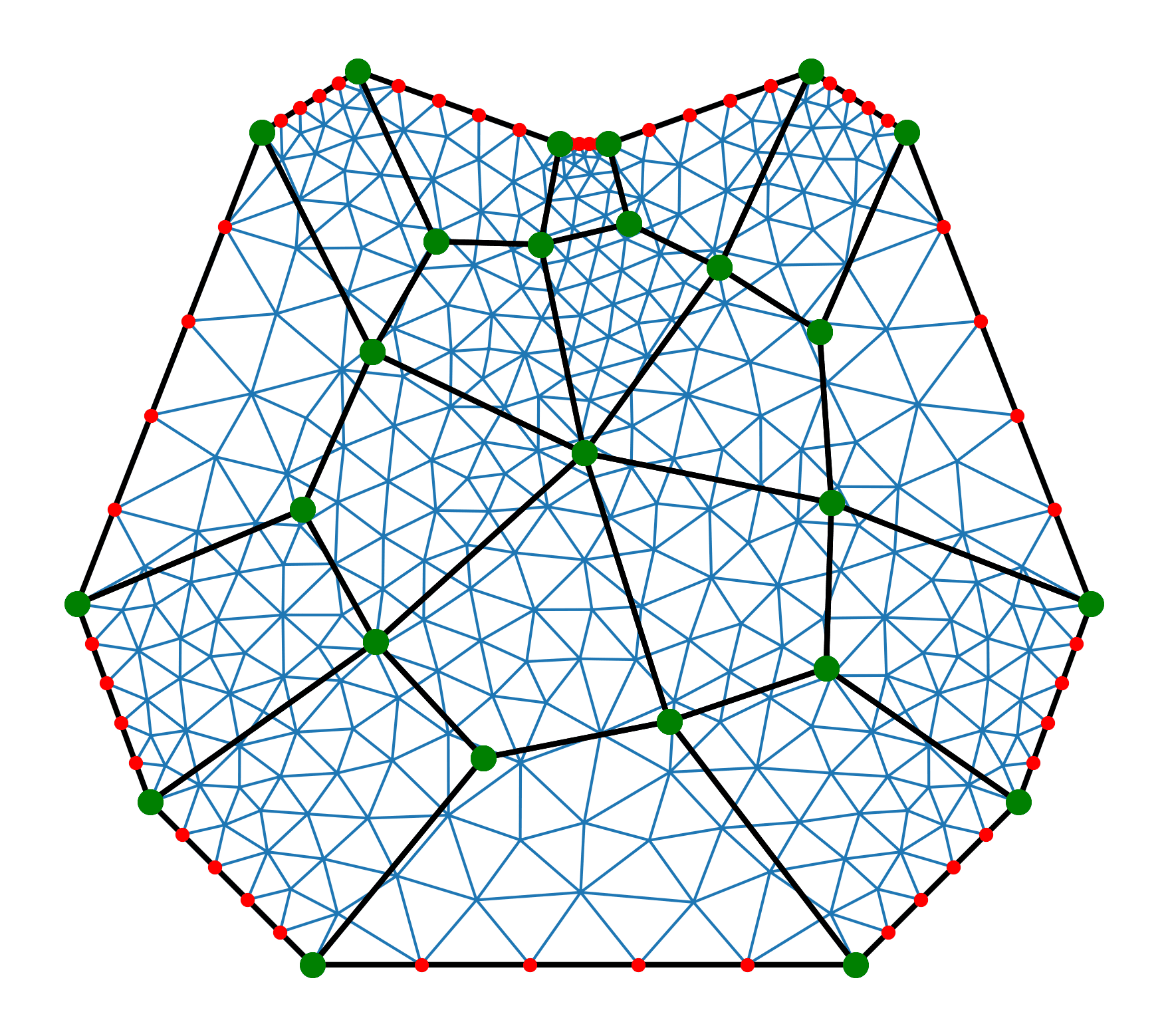}
    \end{subfigure}
    \begin{subfigure}[t]{0.3\textwidth}
        \centering
        \includegraphics[align=c, width=0.95\linewidth]{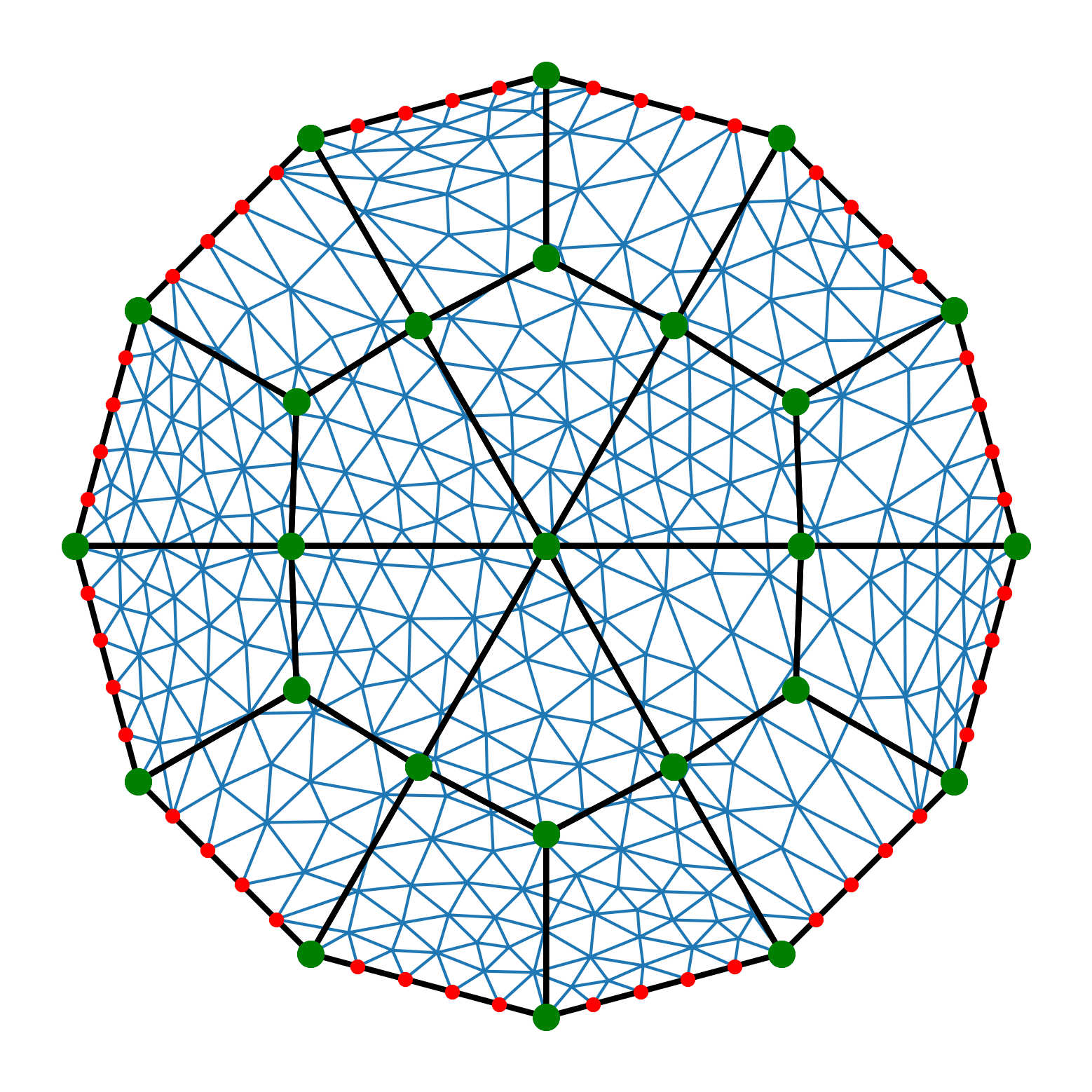}
    \end{subfigure}
\caption{Figure that visualises one step in the creation of the TSG using the techniques from Section \ref{sect:harmonic_maps} and \ref{sect:appendix_floater}. The skeleton (center) of the dense polygon corresponding to a face $F_i \in \mathcal{F}$ (left) is triangulated and the triangular mesh's vertices are mapped into the interior of template $T_i$ (right). Then, the template's vertices are mapped back into the skeleton in order to create a covering for it. The vertices that have been added in order to obtain a finer triangulation are highlighted in red. }
\label{fig:floater}
\end{figure}

\begin{figure}[h!]
\centering
\captionsetup[subfigure]{labelformat=empty}
    \begin{subfigure}[t]{0.4\textwidth}
        \centering
        \includegraphics[align=c, width=0.95\linewidth]{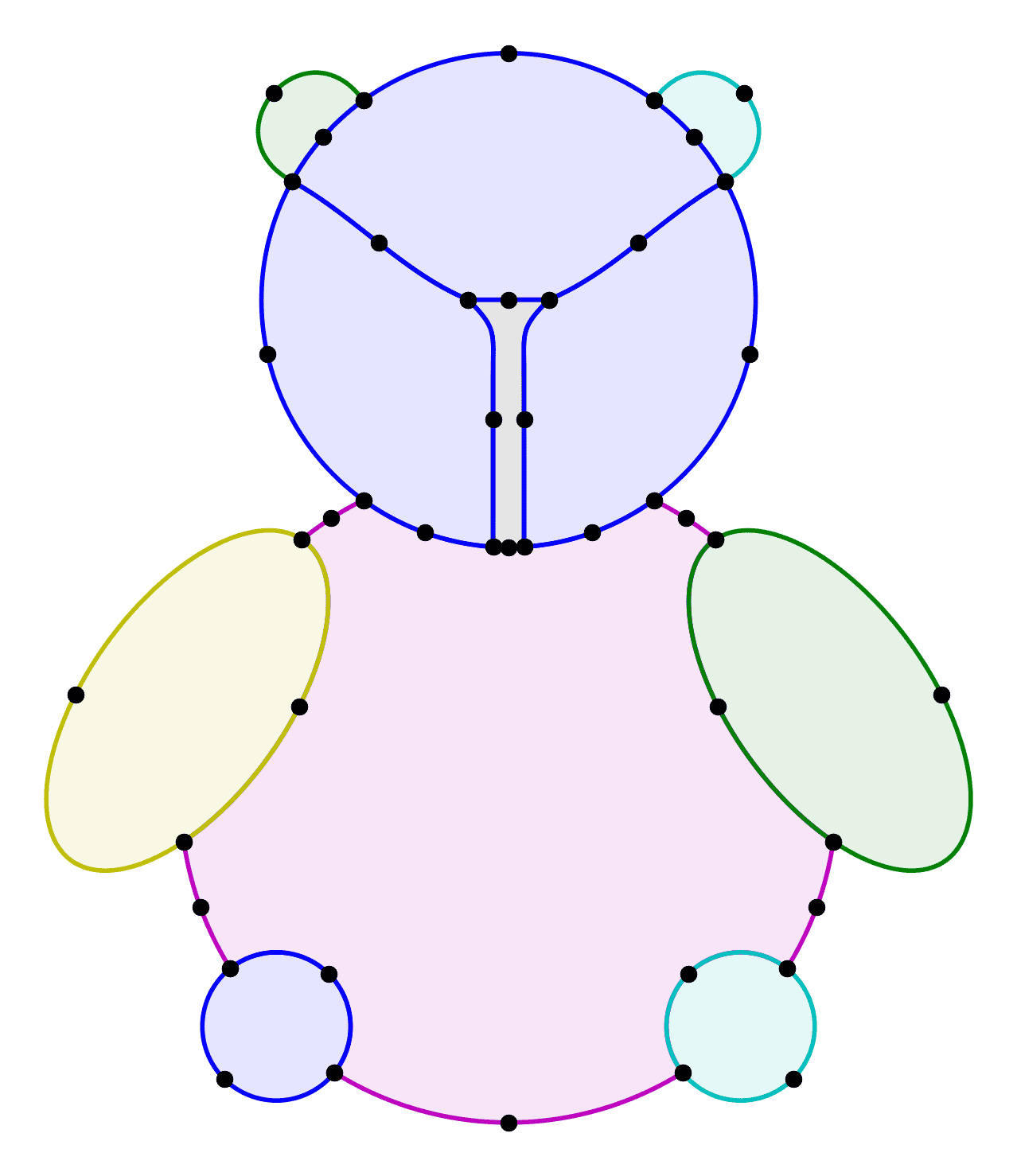}
    \end{subfigure} $\,$
    \begin{subfigure}[t]{0.4\textwidth}
        \centering
        \includegraphics[align=c, width=0.95\linewidth]{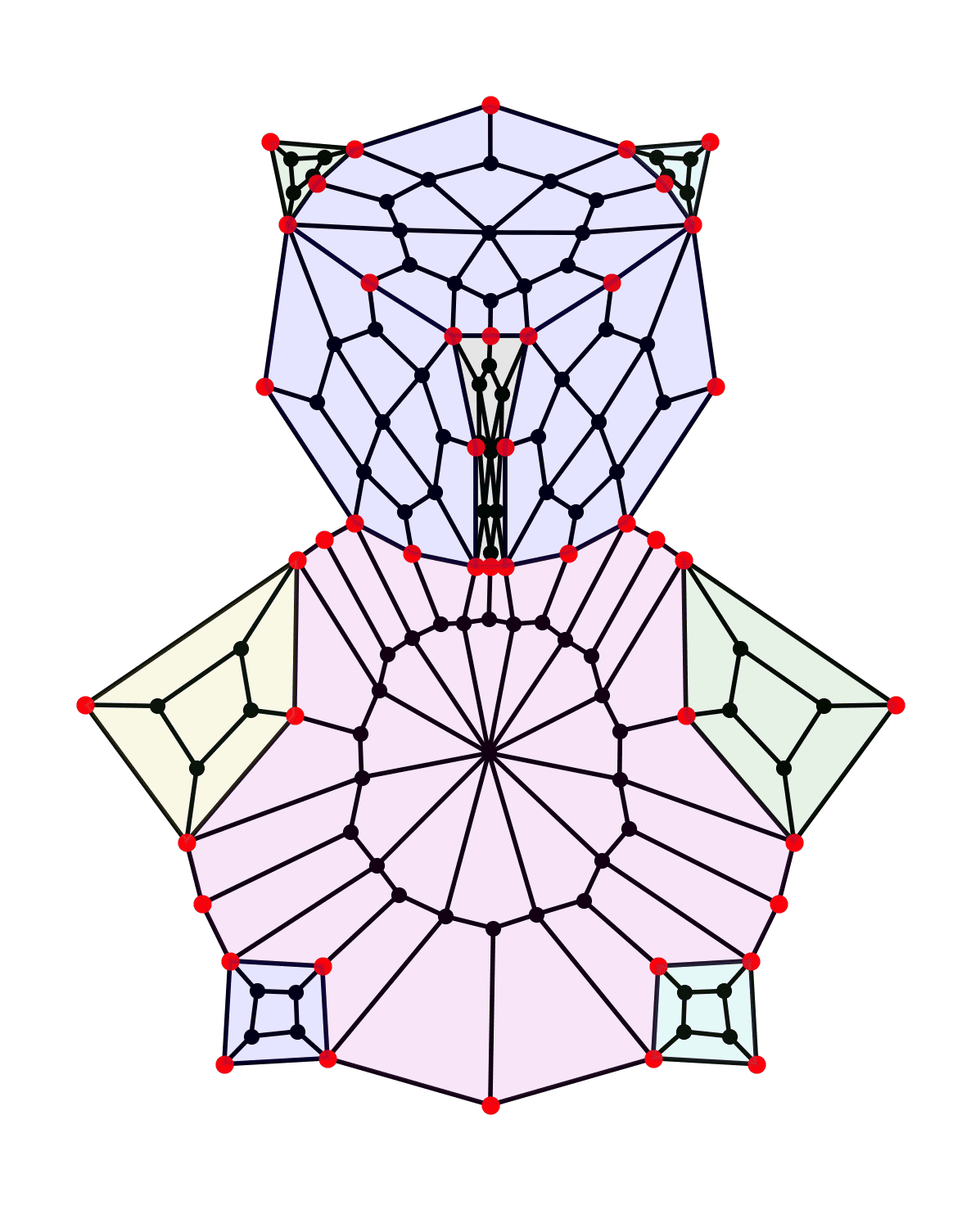}
    \end{subfigure}
\caption{An example of the TSG (right) of a plane graph with only even-sided faces whose associated point sets are shown (left). To each face $F_i \in F$, we assigned the RC $\vert F_i \vert$-leaf template. Here, vertices that follow from the $F_i \in \mathcal{F}$ are represented in red while nodes that correspond to the inner vertices of the $T_i$ are black.}
\label{fig:TSG}
\end{figure}

\noindent A partially templated plane graph can be conveniently visualised through its \textit{template skeleton graph} (TSG). The TSG takes the sub-graph comprised of templated faces $F_i \in \mathcal{F}$ and places additional vertices and edges in the interior of the even-sided polygons associated with the vertices $\mathbb{V}(F_i)$ in a way that mimics the quad layout that has been assigned (see Figure \ref{fig:TSG}). Given the pair $(F_i, T_i)$, we utilise Floater's algorithm (see Section \ref{sect:harmonic_maps}) to map the interior of the polygon with vertices $\mathbb{V}(F_i)$ onto the polygon $\hOm_i$ associated with $\mathbb{V}(\partial E_i)$. This creates a correspondence between mesh vertices in the physical and the parametric domain. The vertices $\hat{v} \in V_i$ of $T_i$ are then mapped into the physical domain using the mesh correspondence. This concept is depicted in Figure \ref{fig:floater}. While the two domains are generally not diffeomorphic, we may ignore this since the piecewise linear correspondence is evaluated in the interior only.  \\

We note that even though the TSG typically suffices for the purpose of visualisation, it may not be a plane graph since two of the sub-graph's linear edges $\{e_i, e_j\} \subset E$ may cross even though the point sets $\{w(e_i), w(e_j) \}$ do not. 

\begin{figure}[h!]
\centering
\includegraphics[align=c, width=0.6\linewidth]{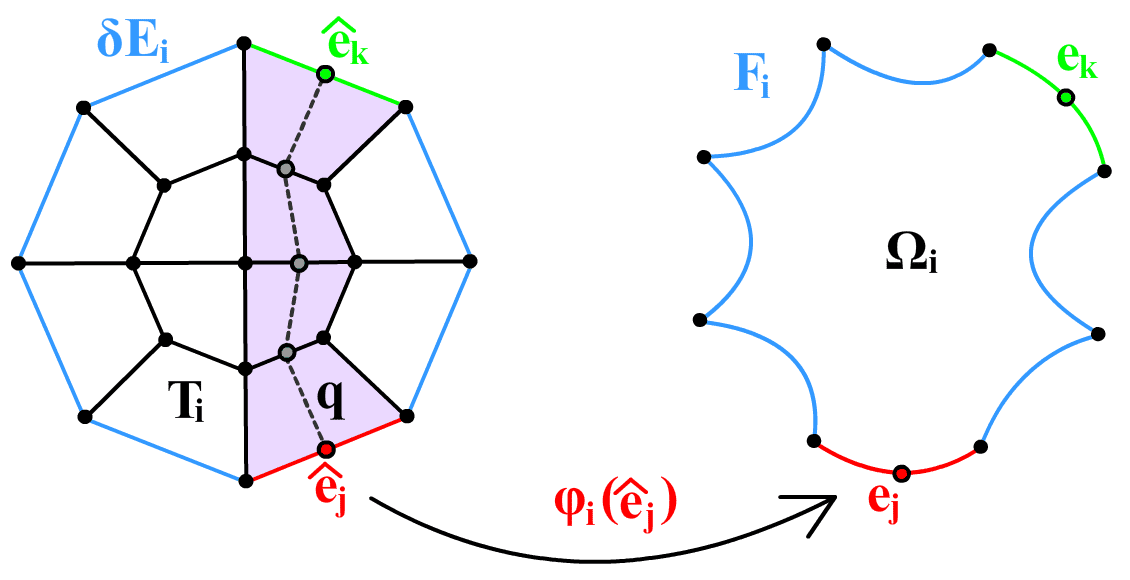}
\caption{Figure showing the concept of assigning a template $T_i = (V_i, E_i, \mathcal{Q}_i)$ to a face $F_i$ and how the assignment pairs edges $\{e_j, e_k\} \subseteq F_i$ under the correspondence $\phi_i: \partial E_i \rightarrow F_i$. The figure additionally shows the refinement path in $T_i$ originating from the root edge $\hat{e}_j = \phi_i^{-1}(e_j)$ which is contained in the root quad cell $q \in \mathcal{Q}_i$.}
\label{fig:template_refine}
\end{figure}

Having introduced the TSG, we are in the position to introduce the concept of template refinement, a critical aspect of point \textbf{B)} of this section's three main ingredients. Given an edge $e \in E$, let $F_e \subseteq \mathcal{F}$ be the set of faces $F_i \in \mathcal{F}$ with $\vert \{e, -e \} \cap F_i \vert \neq 0$ and let $\mathcal{F}_{\mathcal{T}} := \{F_i \in \mathcal{F} \, \vert \, T_i \in \mathcal{T} \}$ be the subset of templated faces. Clearly, we have $\vert F_e \vert \in \{1, 2 \}$.  Template refinement enables creating a new vertex by splitting an existing edge $e \in E$ while automatically performing operations on the graph $G = (V, E, \mathcal{F}, \mathcal{T})$ such that a suitable template structure is retained in the case $F_e \cap F_{\mathcal{T}} \neq \emptyset$. Let $e_j \in E$ be an edge with $| F_{e_j} | = 1$ that we wish to split. Let $F_i \in \mathcal{F}_{\mathcal{T}}$ be the face with $\pm e_j \in F_i$. The edge $e_j$, is associated with exactly one other $e_k \in F_i$ where the pairing follows from $T_i \in \mathcal{T}$ as follows: if the quad cell $q \in \mathcal{Q}_i$ that contains $\partial E_i \ni \hat{e}_j = \phi_i^{-1}(e_j)$ is refined in the direction orthogonal to $\hat{e}_j$, then $e_k = \phi(\hat{e}_k)$, where $\hat{e}_k \in \partial E_i$ is the other boundary edge that is split under the refinement (see Figure \ref{fig:template_refine}). Continuing the pairing on the (if applicable) next templated face, this leads to a sequence of templated edges $e_j, e_k, \ldots, e_\alpha, e_\beta$, with $|F_e| = 2$ for $e \in \{e_k, \ldots, e_\alpha\}$ while $|F_{e_\beta}| = 1$. Template refinement now splits all edges in the sequence $e_j, \ldots, e_\alpha$ while the $T_j \in \mathcal{T}$ associated with the affected faces $F_j \in F_{\mathcal{T}}$ are refined accordingly.

\noindent Let $e \in E$ be an edge from the sequence that is refined. Then, we split $w(e)$ at the point $p \in w(e)$ that minimises the length discrepancy of the resulting point sets $P^-$ and $P^+$. The point $p \in w(e)$ becomes a new vertex in $V$ and $e$ is split into $e^-$ and $e^+$ with $\iota(e^-) = (v_\alpha, p)$ and $\iota(e^+) = (p, v_\beta)$. \\
With the assumptions from Section \ref{sect:Introduction}, no new corner is created at the shared vertex of $(P^-, P^+)$. We finalise the operation by appropriately updating $G = (V, E, \mathcal{F}, \mathcal{T})$ (including the correspondences $\phi_i: \partial E_i \rightarrow F_i$) and the weight function $w(\, \cdot \,)$. 

\begin{figure}[h!]
\centering
    \begin{subfigure}[b]{0.32\textwidth}
        \centering
        \includegraphics[align=c, width=0.95\linewidth]{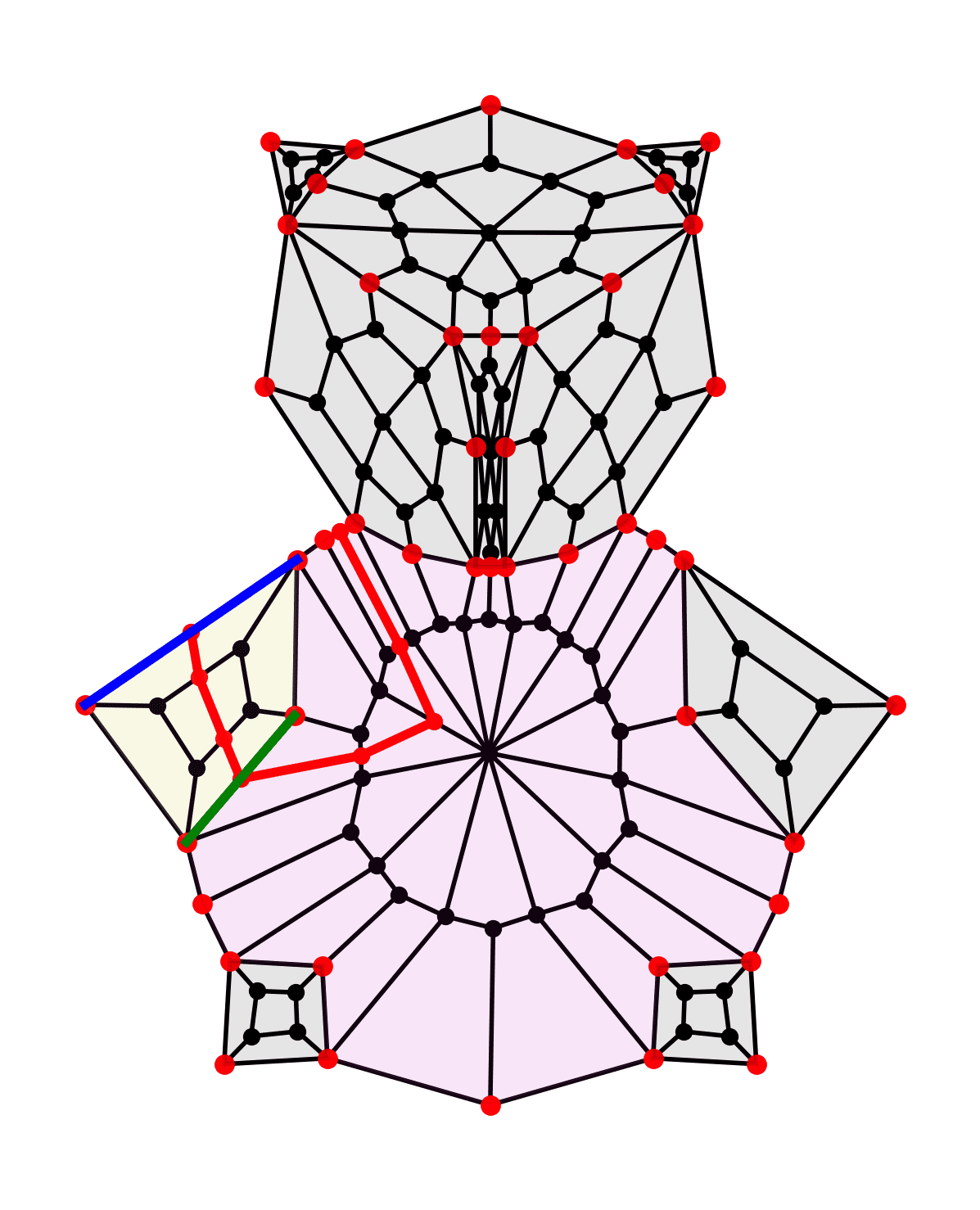}
        \caption{}
    \end{subfigure}
    \begin{subfigure}[b]{0.32\textwidth}
        \centering
        \includegraphics[align=c, width=0.95\linewidth]{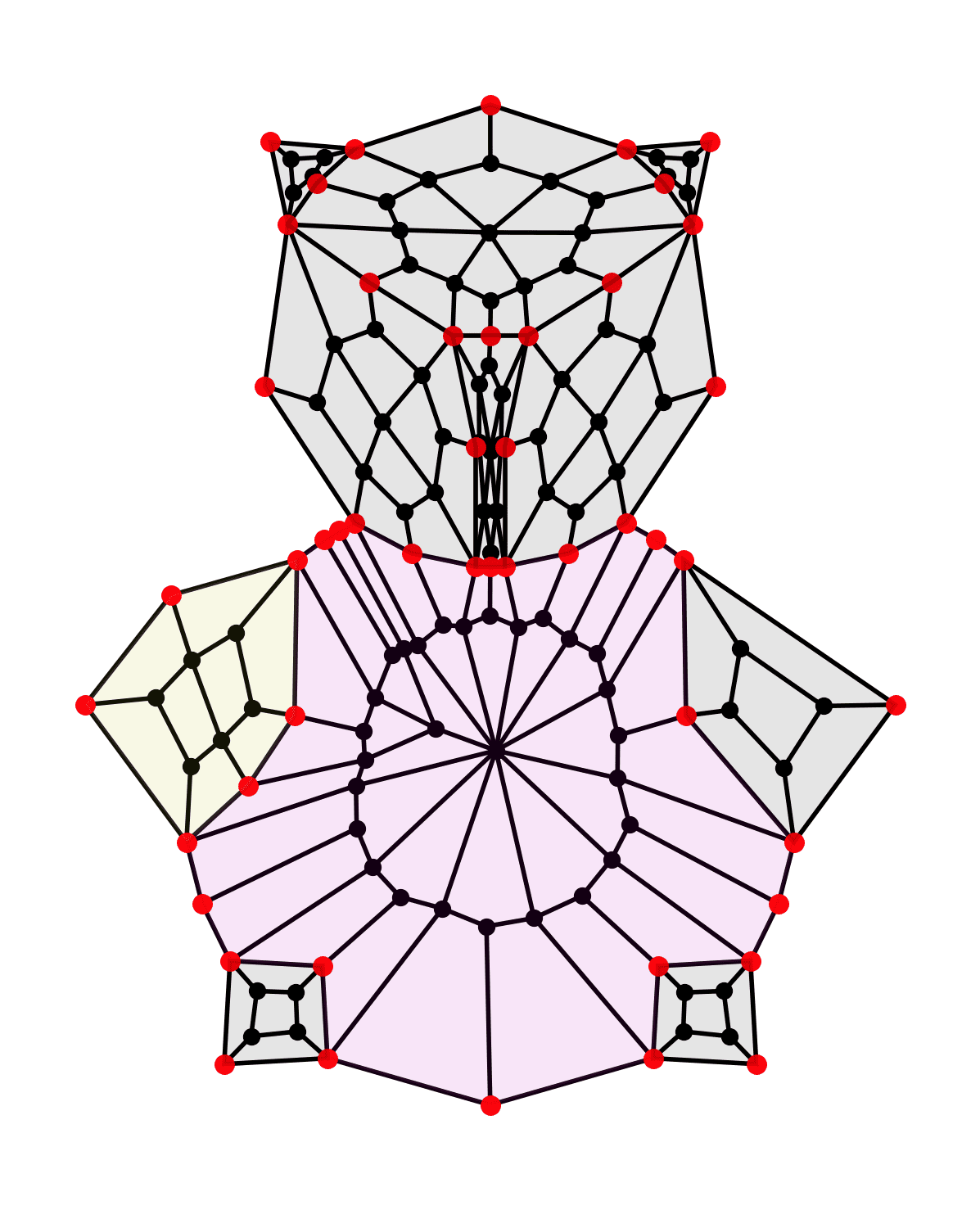}
        \caption{}
    \end{subfigure}
    \begin{subfigure}[b]{0.32\textwidth}
        \centering
        \includegraphics[align=c, width=0.95\linewidth]{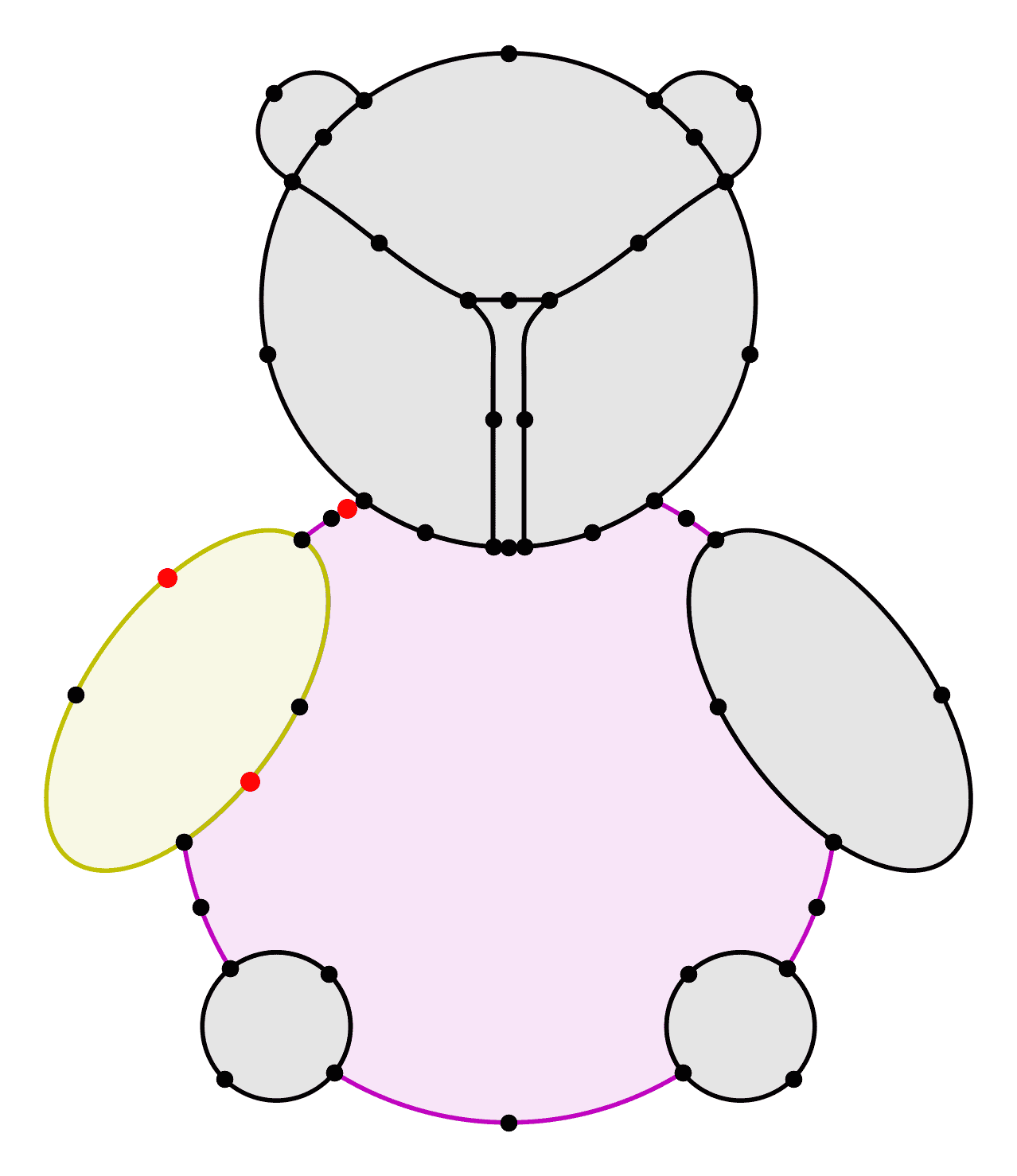}
        \caption{}
    \end{subfigure} \\
    \begin{subfigure}[b]{0.22\textwidth}
        \centering
        \includegraphics[align=c, width=0.95\linewidth]{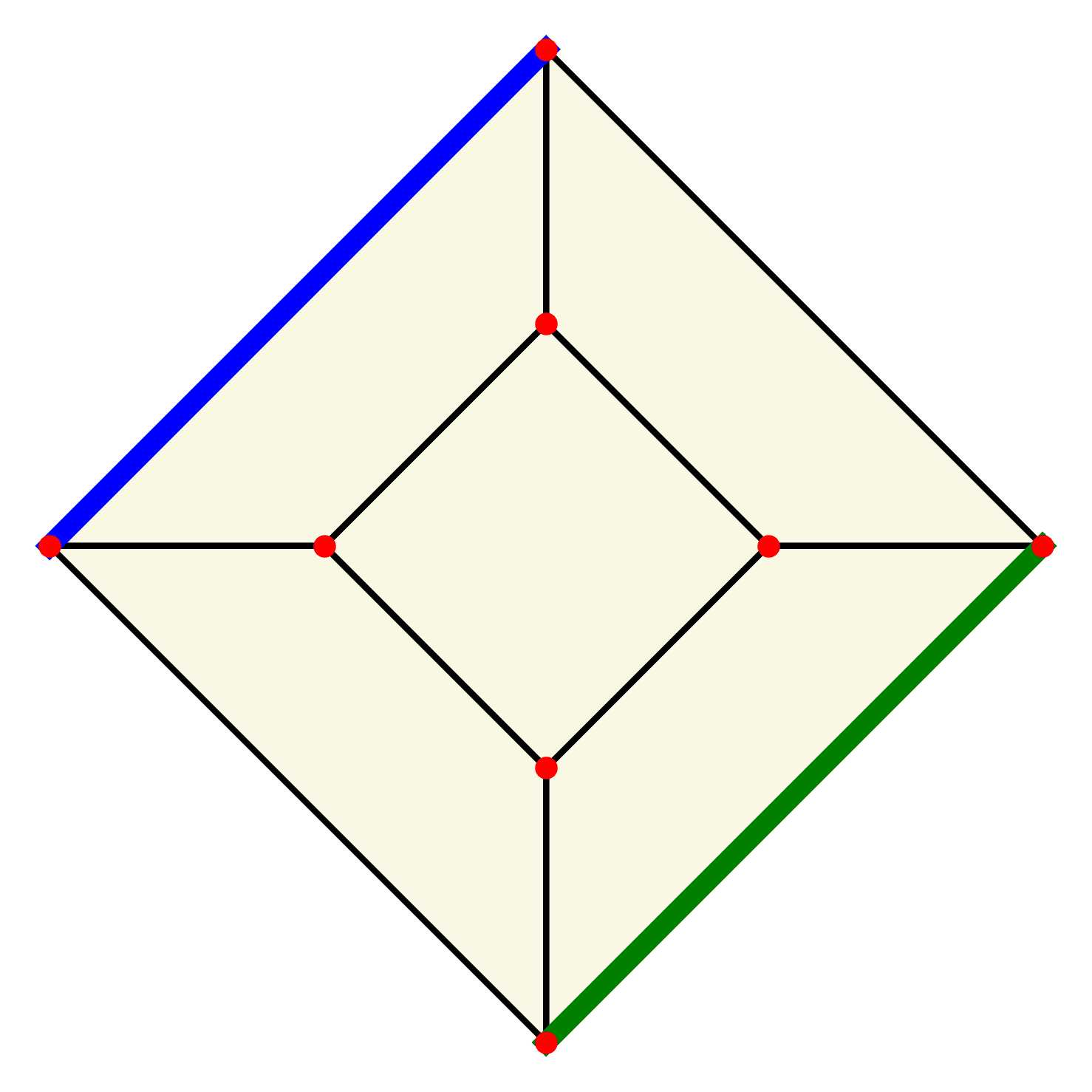}
        \caption{}
    \end{subfigure}
    \begin{subfigure}[b]{0.22\textwidth}
        \centering
        \includegraphics[align=c, width=0.95\linewidth]{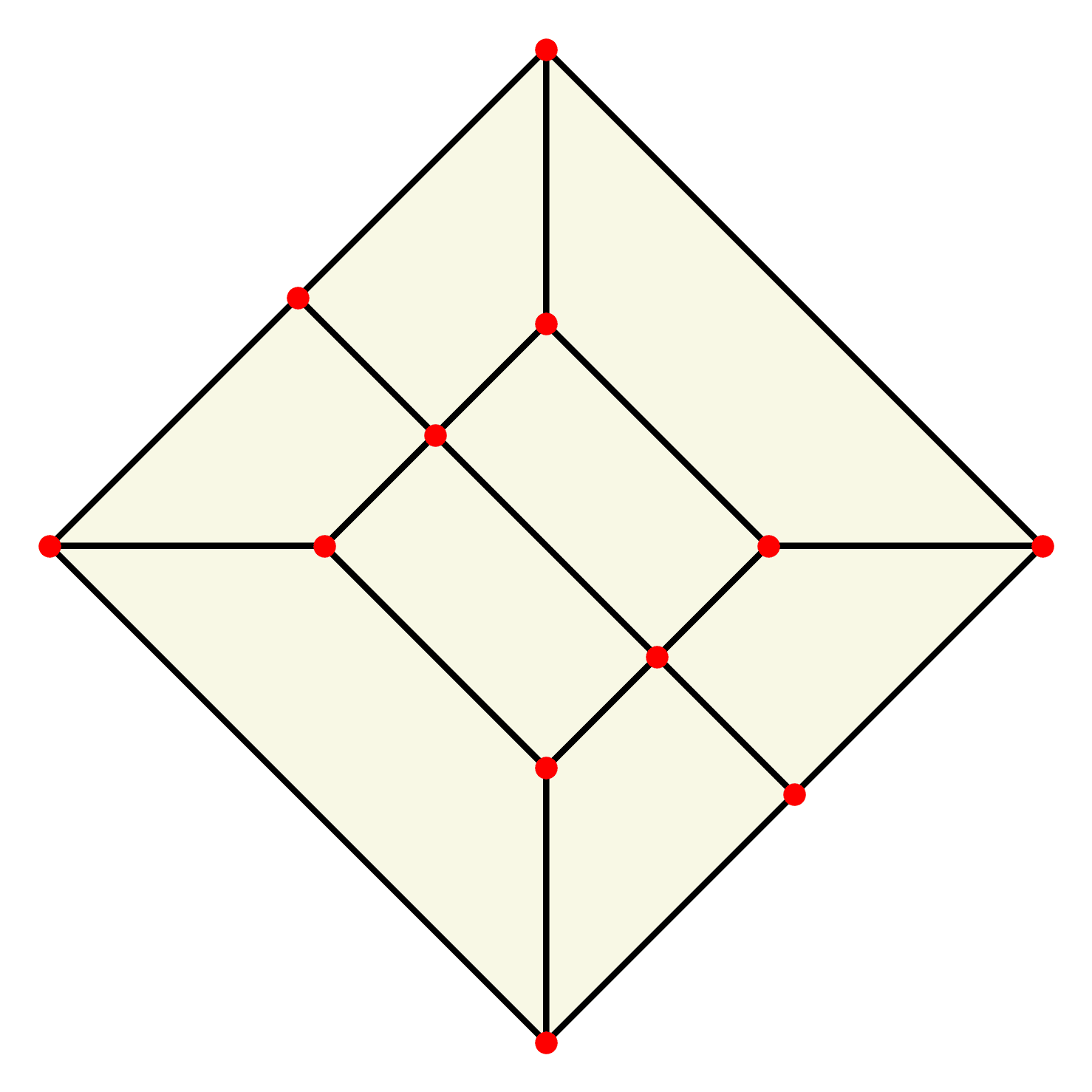}
        \caption{}
    \end{subfigure}
    \begin{subfigure}[b]{0.22\textwidth}
        \centering
        \includegraphics[align=c, width=0.95\linewidth]{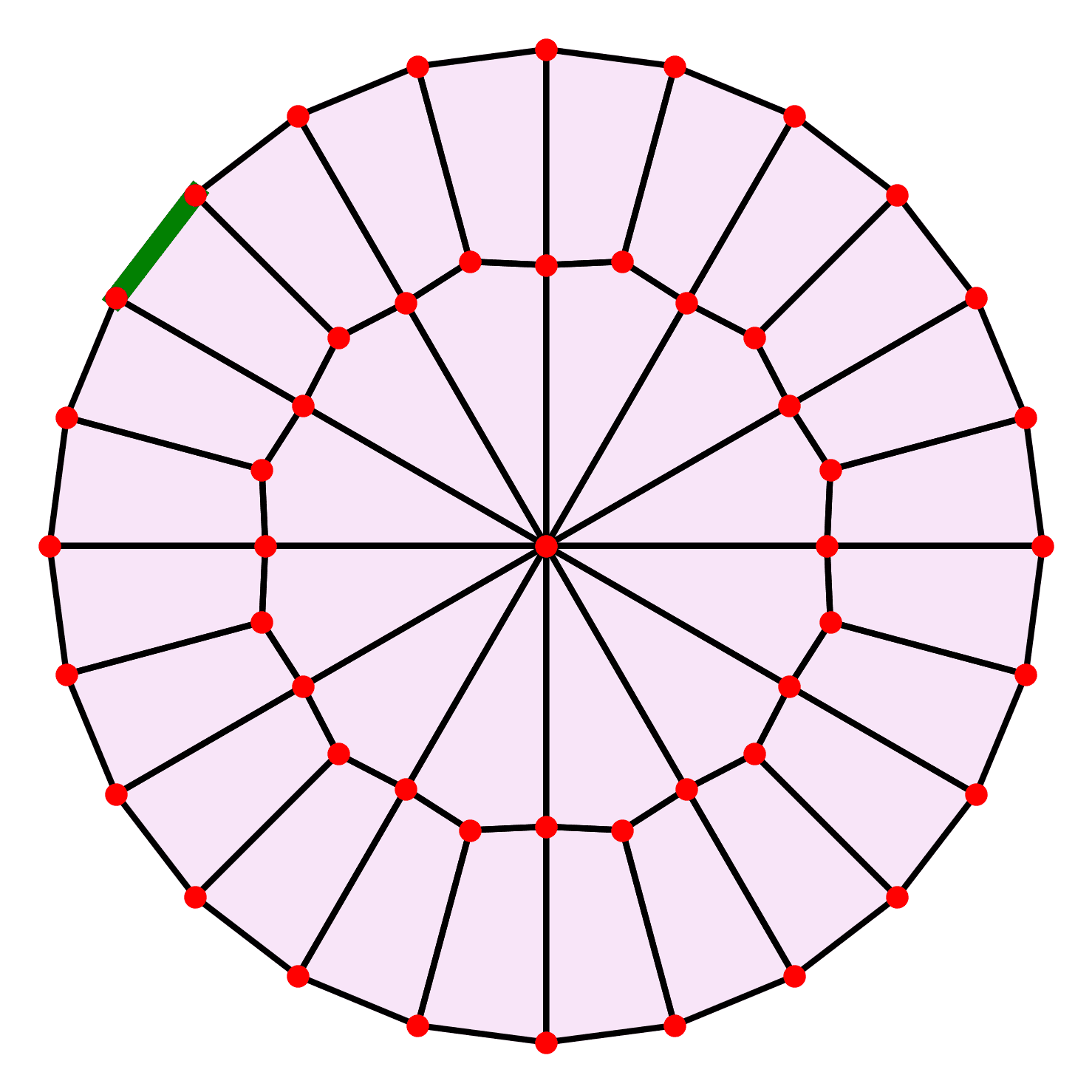}
        \caption{}
    \end{subfigure}
    \begin{subfigure}[b]{0.22\textwidth}
        \centering
        \includegraphics[align=c, width=0.95\linewidth]{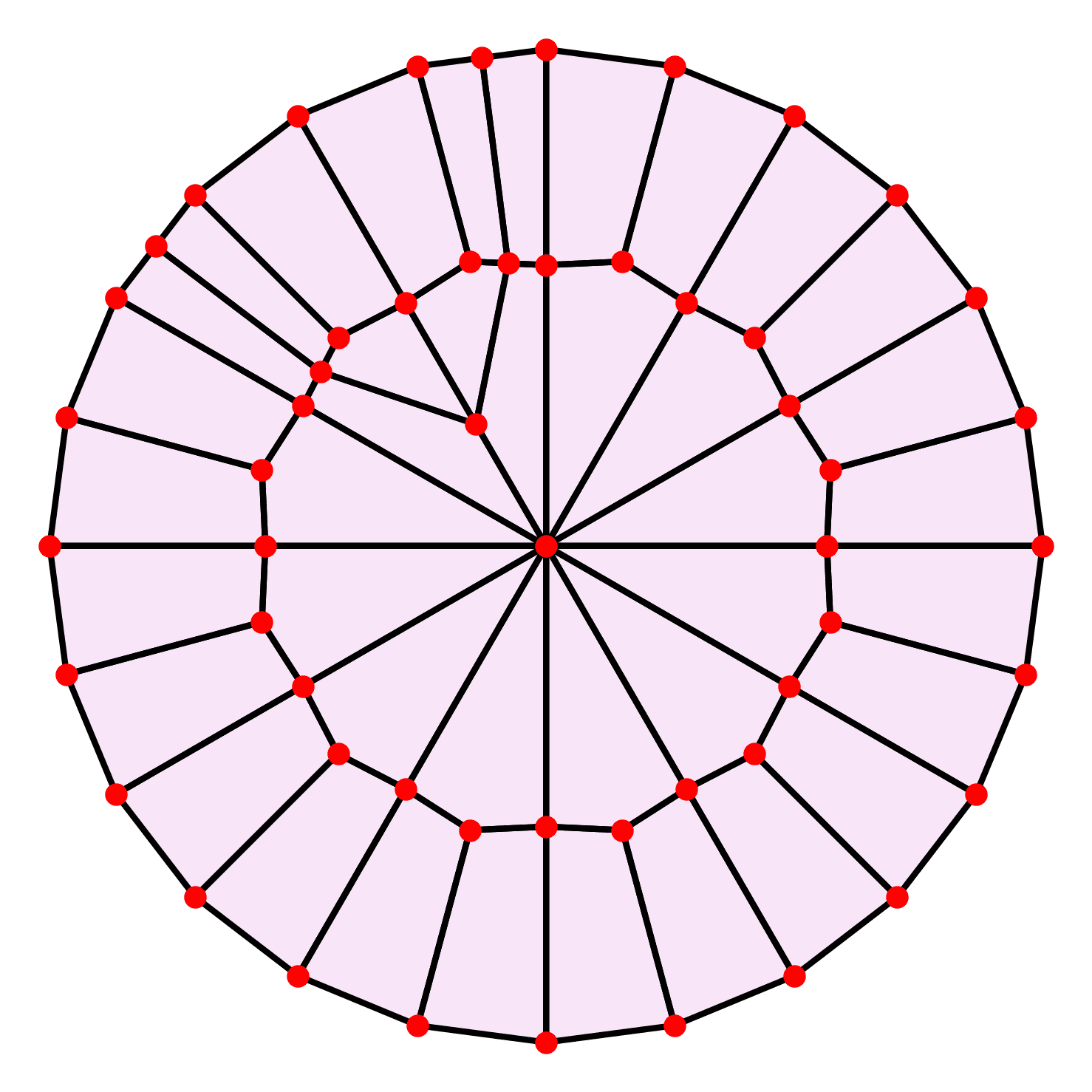}
        \caption{}
    \end{subfigure}
\caption{Figure depicting the concept of template refinement. We are given the TSG from figure \ref{fig:TSG} (right) and perform template refinement on a boundary edge. Figure (a) shows the edge highlighted in blue along with the refinement path along the TSG highlighted in red. The edge that connects the two affected faces along the refinement path is highlighted in green. Figure (b) shows the TSG after template refinement, which furthermore visualises the local refinements performed on the $T_i$. Figure (c) shows the updated point sets $w(e), \, e \in E$ along with three newly created vertices highlighted in red. Figures (d) to (g) depict the affected templates before and after refinement in the corresponding color. Figure (d) highlights the edge that corresponds to the highlighted edge from (a) in the same color, while highlighting the connecting edge between the two templates in green. The connecting edge is also highlighted in (f).}
\label{fig:TSG_trefine}
\end{figure}

This concept is visualised in Figure \ref{fig:TSG_trefine} which adopts the TSG from Figure \ref{fig:TSG} (right) and performs template refinement on the edge $e \in E$ highlighted in blue. The figure furthermore shows the refinement path and the updated TSG and highlights the newly created vertices $v \in V$ of the updated graph $G = (V, E, \mathcal{F}, \mathcal{T})$ while furthermore showing the two affected templates before and after refinement. Finally, the patchwise smooth controlmap $\mathbf{r}: \hOm_i \rightarrow \hOm^{\mathbf{r}}_i$ need not be updated because refinement does not change the corner pattern of the updated $F_i \in \mathcal{F}$ by assumption. \\

\noindent The concept of template refinement puts us in the position to develop concrete templatisation strategies. \\
The core of the iterative strategy, in words, is comprised of the following steps:
\begin{enumerate}
    \item Select an untemplated face $F \in \mathcal{F}$.
    \item If $\vert F \vert$ is odd or $\vert F \vert = 2$, select edge(s) $e \in F$ and perform (template) refinement on the edge(s).
    \item Select the best template $T$ from a large catalogue $\mathbb{T}$ and assign it to $F$.
    \item Repeat until all faces have been templated.
\end{enumerate}

Given the partially templated graph $G = (V, E, \mathcal{F}, \mathcal{T})$, let $E_{\mathcal{T}} \subset E$ denote the subset of edges associated with at least one template. The default face selection strategy (Point 1.), along with comments that explain each step, is given in Algorithm \ref{algo:select_face}. \\

\begin{algorithm}
    \caption{Select face}
    \label{algo:select_face}
    \begin{algorithmic}[1]	
        \Procedure{Select face}{Graph $G = (V, E, \mathcal{F}, \mathcal{T})$}
	   \Statex
          \State $\mathcal{F}_{\text{select}} \leftarrow \left \{ F \in \mathcal{F} \, \vert \, F \text{ is not templated} \right \}$ \Comment{Keep only untemplated faces.}
          \If{$\mathcal{F}_{\text{even}} := \left \{ F \in \mathcal{F}_{\text{select}} \, \vert \, \vert F \vert \geq 4 \text{ and } F \text{ is even-sided} \right \} \neq \emptyset$}
            \State $\mathcal{F}_{\text{select}} \leftarrow \mathcal{F}_{\text{even}}$ \Comment{If there are even-sided faces with $\vert F \vert \geq 4$, keep only them.}
          \EndIf
          \State $m \leftarrow \min \left \{\vert F \cap E_{\mathcal{T}} \vert \, \vert \, F \in \mathcal{F}_{\text{select}} \right \}$ \Comment{Get $m$, the minimum number of templated edges.}
          \State Update $\mathcal{F}_{\text{select}} \leftarrow \left \{ F \in \mathcal{F}_{\text{select}} \, \vert \, \vert F \cap E_{\mathcal{T}} \vert = m \right \}$ \Comment{Keep only faces with $m$ templated edges.}
          \Statex
          \Statex
        \Return{$\operatorname{argmin}_{\substack{F \in \mathcal{F}_{\text{select}}}} \vert F \vert$} \Comment{Of the remainder, return the face with the fewest number of edges.}
	\EndProcedure{}
    \end{algorithmic}
\end{algorithm}

Since the algorithm may perform operations that change the graph, we note that the greedy strategy from Algorithm \ref{algo:select_face} has to be called on the updated graph at the beginning of each iteration. \\

\noindent The next ingredient (see Point 2.) is a strategy for selecting a splitting edge in case $\vert F \vert$ is odd. The case $\vert F \vert = 2$ is treated separately. \\

Clearly, a possible, albeit crude, strategy is halving each edge $e \in E$ in a pre-step. This would ensure that all faces are even-sided with $\vert F \vert \geq 4$ before template assignment commences. However, we note that this significantly increases the number of required patches, an effect that may not be desirable. \\

By default, we favour longer edges for splitting. We favour untemplated edges $e \notin E_{\mathcal{T}}$ over templated edges $e \in E_{\mathcal{T}}$ and among the two groups, we favour boundary edges $e \in \partial E$ over interior edges. For this purpose, we introduce the penalty factors $\mu_{\mathcal{T}} \geq 1$ and $\mu_{\partial} \geq 1$.  Let $L(e)$, with $U(e) \in E$, denote the length of the piecewise linear curve that results from connecting the points $p \in w(e)$ by linear edges. Applying $L(\, \cdot \,)$ to some $F \in \mathcal{F}$, instead, performs the same operation element wise. We introduce the scaled length function
\begin{align}
\label{eq:scaled_edge_length}
    L^{\mu_{\mathcal{T}}, \mu_{\partial}}(e) := \left \{ \begin{array}{ll} 1 & e \in E_{\mathcal{T}} \\
                                                     \mu_{\mathcal{T}} & \text{else} \end{array} \right. \quad \times \quad \left \{ \begin{array}{ll} 1 & e \notin \partial E \\
                                                     \mu_{\partial} & \text{else} \end{array} \right. \quad \times \quad L(e).
\end{align}
By $L_{\text{max}}(F) := \max L(F)$ we denote the longest edge $e \in F$. Given some $0 < \varepsilon_L \leq 1$, the first step keeps only the edges $e \in F$ with 
\begin{align}
\label{eq:edge_length_filter}
    L(e) \geq \varepsilon_L L_{\text{max}}(F).
\end{align}
Note that after filtering, the set of eligible splitting edges $E_{\text{split}}(F)$ is nonempty. Of the remaining edges, we then select the one that maximises~\eqref{eq:scaled_edge_length}. In practice, we utilise $\mu_{\mathcal{T}} = 2$ and $\mu_{\partial} = 1.5$ while $\varepsilon_L$ is set to $\varepsilon_L = 0.3$. \\

\noindent For $\vert F \vert = 2$, we have the choice to either split both edges $e \in F$ or to split a single $e \in F$ twice. The selection strategy is similar to $\vert F \vert > 2$. However, we are more lenient in the first step by replacing $\varepsilon_L \rightarrow \tfrac{1}{2} \varepsilon_L$ in~\eqref{eq:edge_length_filter}. Next, we sort the remaining edges by~\eqref{eq:scaled_edge_length} in descending order. Then, if both edges $\{\hat{e}_1, \hat{e}_2 \} \subseteq F$ remain, we return either $\{\hat{e}_1 \}$ if $L^{\mu_{\mathcal{T}}, \mu_{\partial}}(\hat{e}_{1}) \geq 2 L^{\mu_{\mathcal{T}}, \mu_{\partial}}(\hat{e}_2)$ and $\{\hat{e}_1, \hat{e}_2 \}$ else. In the former case, the new vertices are set to the $p \in w(\hat{e}_1)$ that split the point set most evenly at $\tfrac{1}{3} L(e)$ and $\tfrac{2}{3} L(e)$, respectively while the edge $\hat{e}_{2}$ is not refined. \\

\noindent After (potentially) performing (template) refinement on the marked edge(s), we are in the position to select a template $T_i = (V_i, E_i, \mathcal{Q}_i)$ for the face $F \in \mathcal{F}$ from a large catalogue $\mathbb{T}$ of precomputed layouts. The non-optional pre-filter from Algorithm \ref{algo:template_pre_filter} is subject to the user-specified angle threshold $\mu_\angle$ which retains only templates with $\operatorname{val}(\hat{v}) \geq 3$ if some $\mathbb{V}(F) \ni v = \phi(\hat{v}), \, \hat{v} \in V_i$ satisfies $\angle(v) \geq \pi - \mu_\angle$. \\

\begin{algorithm}
    \caption{Template pre-filter}
    \label{algo:template_pre_filter}
    \begin{algorithmic}[1]	
        \Procedure{Template pre-filter}{Face $F \in \mathcal{F}$, catalogue $\mathbb{T}$, angle threshold $\mu_{\angle}$}
	   \Statex
          \State $\mathcal{T}_{\text{select}} \leftarrow \left \{ T_i \in \mathbb{T} \, \big \vert \, \vert \partial E_i \vert = \vert F \vert \right \}$
          \State Update $\mathcal{T}_{\text{select}} \leftarrow \left \{ T \in \mathcal{T}_{\text{select}} \, \big \vert \, \text{ if } \angle(v) \geq \pi - \mu_{\angle} \text{, then } \hat{v} := \phi^{-1}(v) \text{ satisfies } \operatorname{val}(\hat{v}) \geq 3  \right \}$
          \Statex
          \Statex
        \Return $\mathcal{T}_{\text{select}}$
	\EndProcedure{}
    \end{algorithmic}
\end{algorithm}

\begin{figure}[h!]
\centering
    \begin{subfigure}[b]{0.22\textwidth}
        \centering
        \includegraphics[align=c, width=0.95\linewidth]{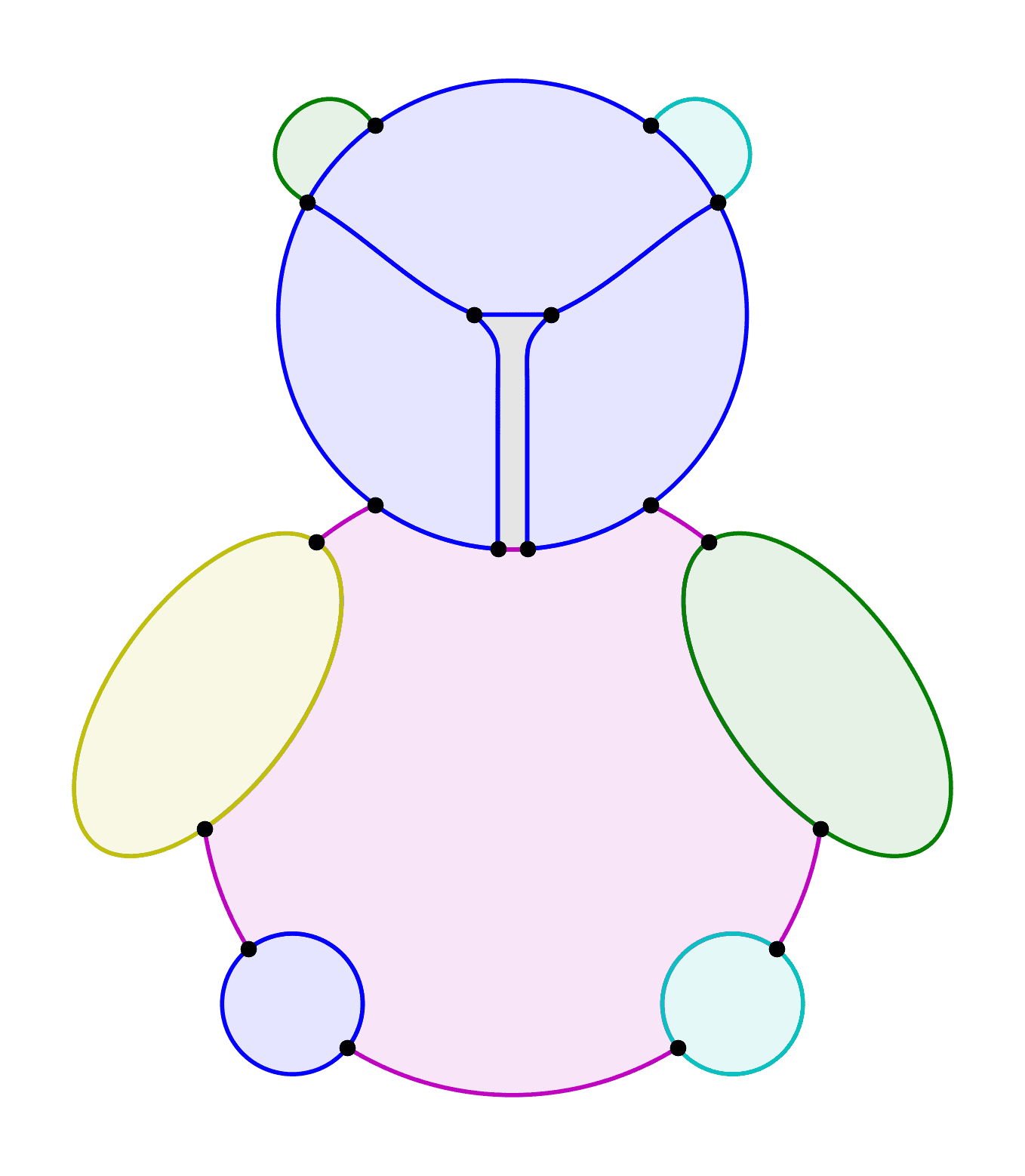}
        \caption{}
        \label{eq:minpatch_strategy_symmetric_input_graph}
    \end{subfigure}
    \begin{subfigure}[b]{0.22\textwidth}
        \centering
        \includegraphics[align=c, width=0.95\linewidth]{connie.png}
        \caption{}
    \end{subfigure}
    \begin{subfigure}[b]{0.22\textwidth}
        \centering
        \includegraphics[align=c, width=0.95\linewidth]{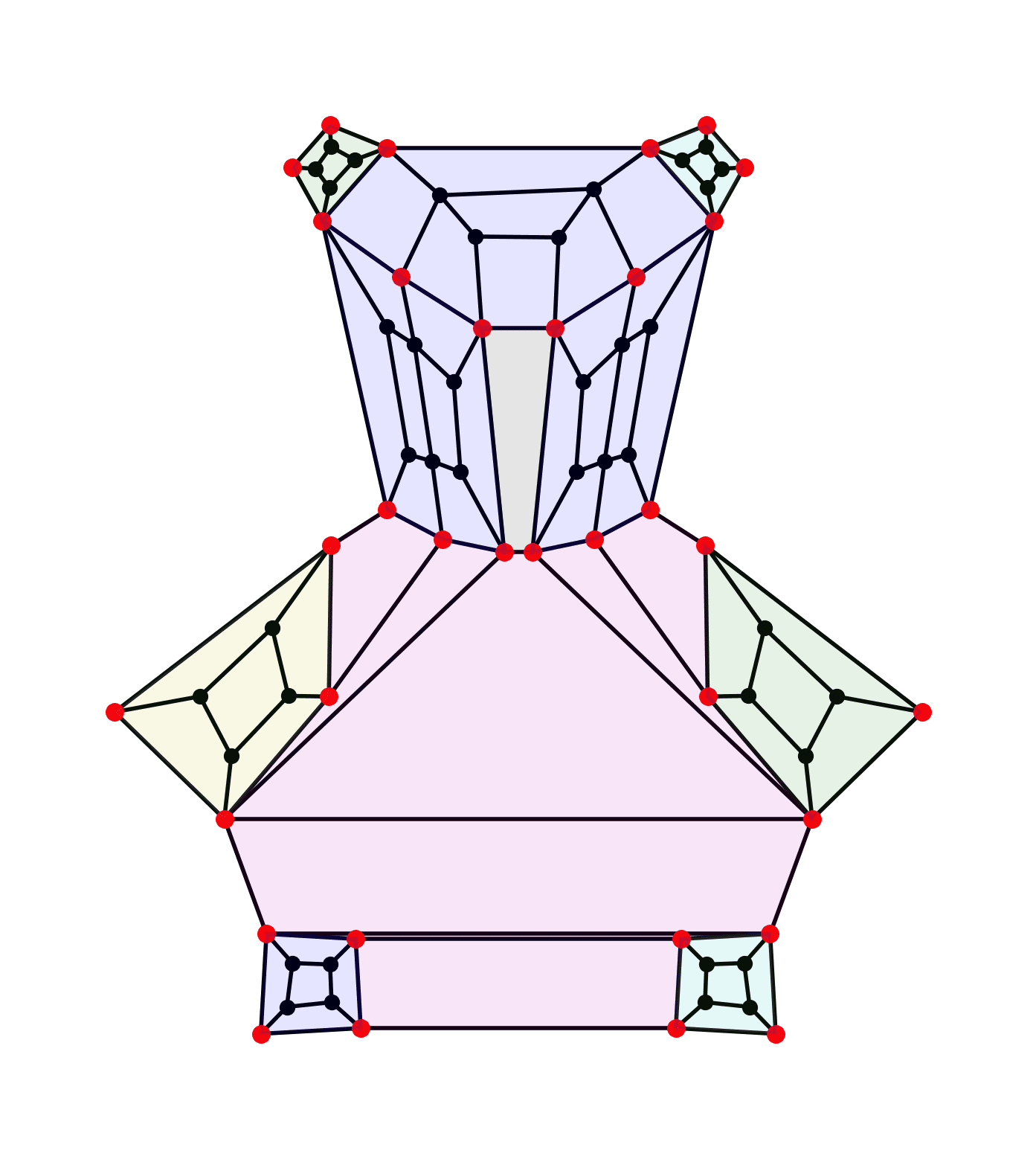}
        \caption{}
    \end{subfigure}
    \begin{subfigure}[b]{0.22\textwidth}
        \centering
        \includegraphics[align=c, width=0.95\linewidth]{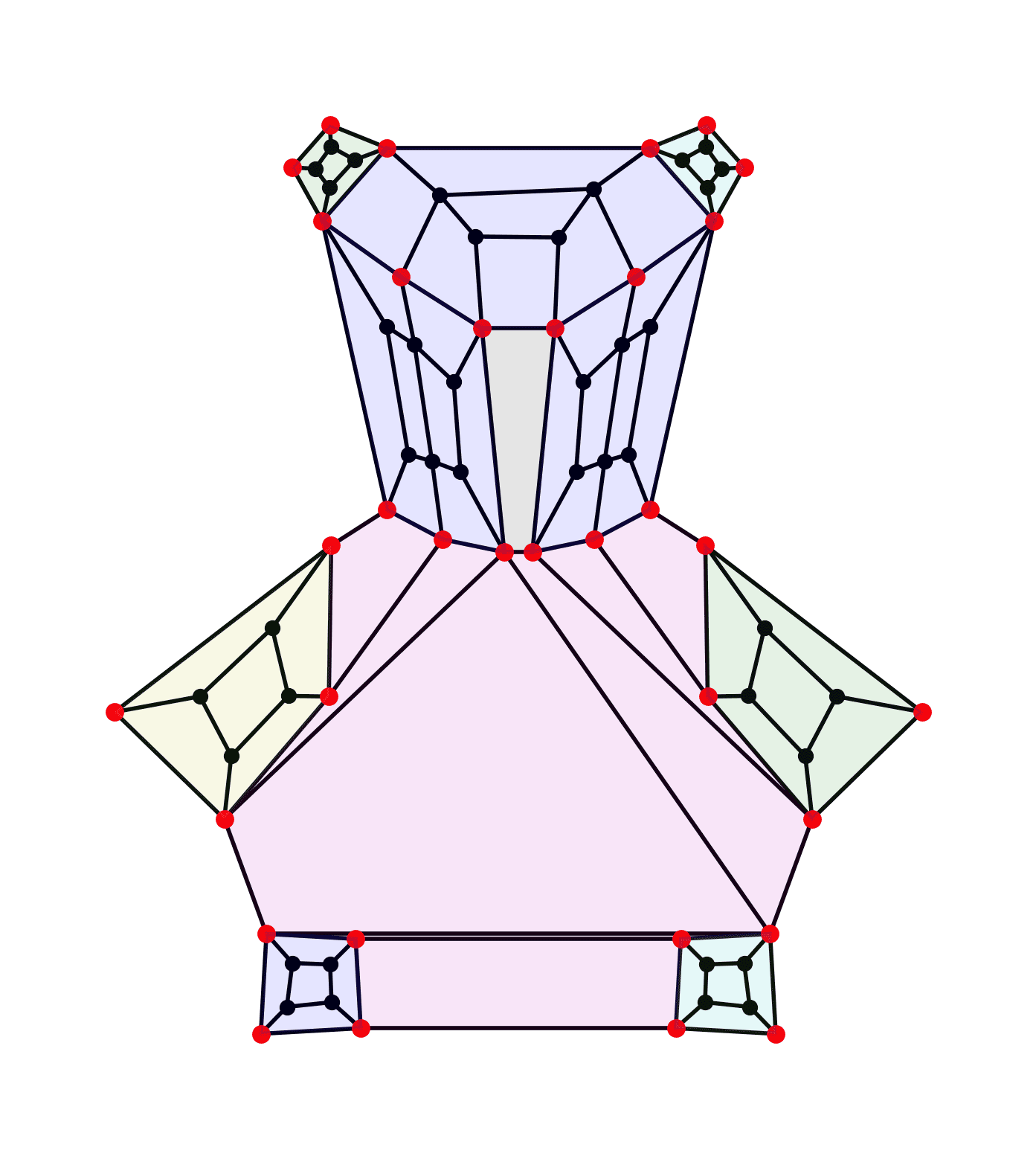}
        \caption{}
    \end{subfigure}
\caption{Figure showing the TSG of the input graph \textbf{(a)} using the min-patch strategy with respecting symmetries \textbf{(c)} and without \textbf{(d)}. In \textbf{(d)}, a template that cannot match the pink face's symmetry axes has been chosen. Figure \textbf{(b)} highlights the vertices that were added during the templatisation stage in red. The total number of patches may be further reduced by stronger penalising the splitting of templated edges via the parameter $\mu_{\mathcal{T}} \geq 1$.}
\label{fig:minpach_strategy_symmetric}
\end{figure}

After applying the necessary pre-filter, we apply an application-dependent post filter which returns the best template according to a specified quality criterion. Here we distinguish between two strategy types:
\begin{enumerate}[label=\textbf{Type }\arabic*., align=left]
    \item Heuristic filters;
    \item Deterministic filters.
\end{enumerate}
Type 1. filters are based on a set of heuristic, computationally lightweight selection criteria. For instance, a widely employed heuristic filter retains only the templates that match (if present) the face's symmetries and then selects the layout with the smallest number of patches. Figure \ref{fig:minpach_strategy_symmetric} depicts the TSG of the heuristic min-patch strategy with and without favouring symmetries. \\

\noindent Type 2. strategies make concrete predictions about the performance of some $\mathbb{T}_{\text{select}} \ni T_i = (V_i, E_i, \mathcal{Q}_i)$, where $\mathbb{T}_{\text{select}}$ is returned by Algorithm \ref{algo:template_pre_filter}, when employed in this paper's parameterisation method based on harmonic maps. Let $F \in \mathcal{F}$ be a face and let $\Om$ be the associated domain with $\hOm^\br$ the canonical control domain (see Section \ref{sect:harmonic_maps}). To enable testing a large number of templates, we utilise Floater's algorithm to approximate a harmonic map between the polygons $\Om_{h}$ and $\hOm^{\mathbf{r}}_{h}$, both (relatively) coarsely sampled from $\partial \Om$ and $\partial \hOm^{\mathbf{r}}(F)$, respectively. We sample $N_{\text{sample}}$ points from each edge of $\partial \Omega$ and $\partial \hOm^\br$. The polygon $\Om_h$ is triangulated using Gmsh \cite{geuzaine2009gmsh} and a surrogate (approximately harmonic) map $\mathbf{x}_{h}^{\mathbf{r}}: \hOm^\mathbf{r}_{h} \rightarrow \Om_{h}$ is constructed from the correspondence between $\Om_h$ and $\hOm^\br_h$ using the techniques from Section \ref{sect:harmonic_maps} and \ref{sect:appendix_floater}.

\begin{figure}[h!]
\centering
\includegraphics[align=c, width=0.7\linewidth]{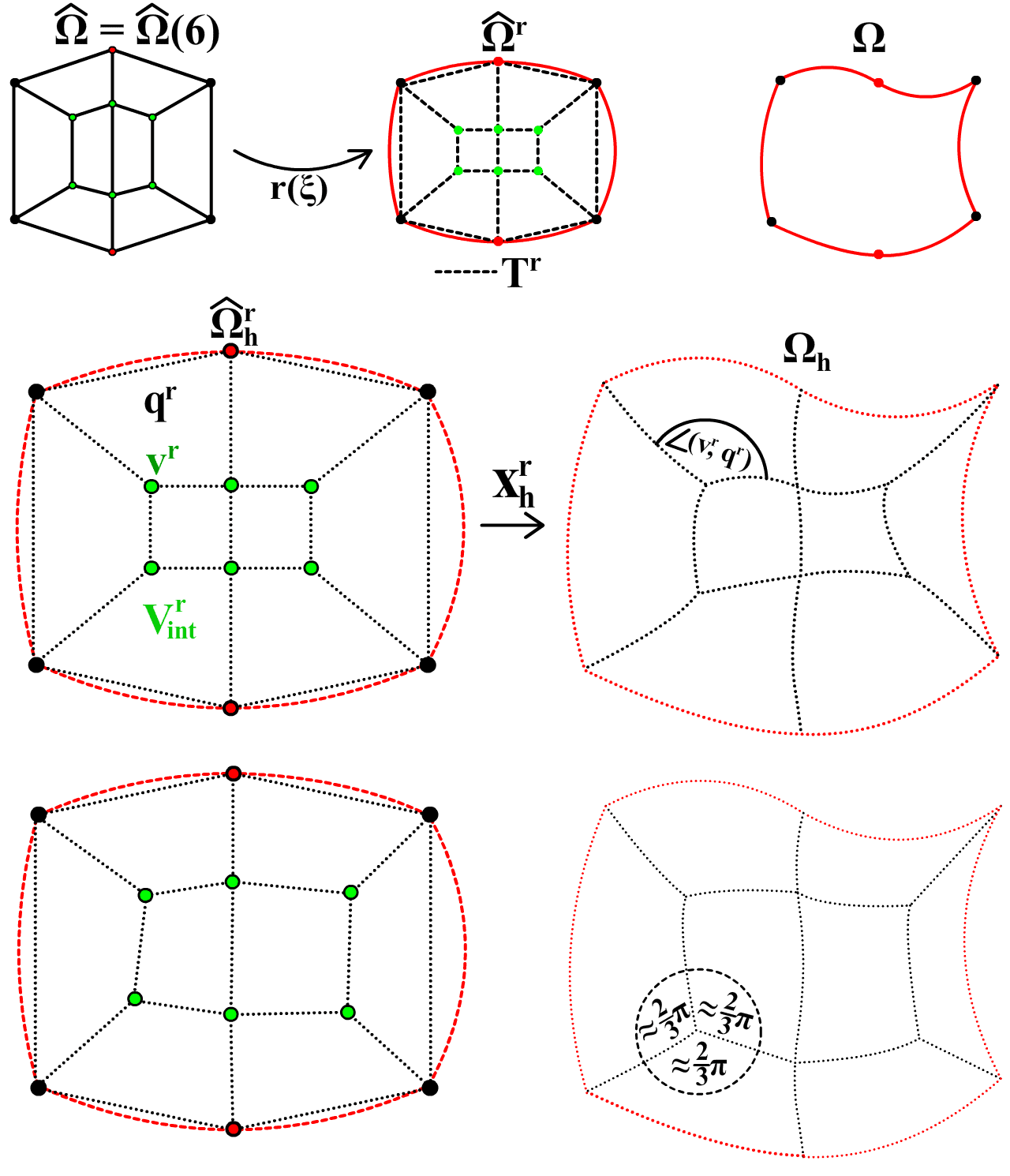}
\caption{Figure depicting the concept of utilising the surrogate map $\mathbf{x}_h^{\mathbf{r}}$ between the surrogate domains $\widehat{\Omega}_h^{\mathbf{r}} \approx \widehat{\Omega}^{\mathbf{r}}$ and $\Omega_h \approx \Omega$ to optimise the inner vertex positions (green) of the template $T^r = (V^r, E^r, \mathcal{Q}^r)$. Here, the inner vertex $v^r \in V_{\text{int}}^r$ positions are chosen such that the patch interface angles created by the quad layout vertices minimises the deviation from the preferred angle (see equation~\eqref{eq:preferred_angle_selection}). While the presented strategy has demonstrated effective results in practice, it is important to acknowledge that alternative quality metrics deserve further exploration.}
\label{fig:opt_angle}
\end{figure}

\noindent Given a candidate template $T_i \in \mathbb{T}_{\text{select}}$, we compute the canonical controlmap $\br:\hOm(\vert F \vert) \rightarrow \hOm^\br$ as a function of $T_i$ and $\hOm^\br$, see Section \ref{sect:harmonic_maps}. We use the controlmap to create a weight function $w^\bxi(\, \cdot \,)$ that assigns to an edge $\hat{e} \in E_i$ of $T_i$ a point set $p = w^\bxi(\hat{e})$ uniformly sampled over $\hat{e} \in E_i$ from the composition $\mathbf{x}_{h}^{\mathbf{r}} \circ \mathbf{r}: \hOm(\vert F \vert) \rightarrow \Om_{h}$.
Note that for this, the map $\mathbf{r}: \hOm(\vert F \vert) \rightarrow \hOm^\mathbf{r}$ needs to be only known on the edges $\hat{e} \in E_i$. By default, we assume that for boundary edges $\hat{e} \in \partial E_i$, the sampling density matches that of the sides of $\partial \hOm^{\mathbf{r}}_{h}$. With $w^\bxi(\, \cdot \,)$ at hand, we may now define a cost function on the pair $(T_i, w^\bxi)$ to gauge the quality of a template choice. \\

\noindent While more computationally demanding than type 1. strategies, type 2. strategies remain relatively inexpensive because the map $\mathbf{x}_{h}^{\mathbf{r}}: \hOm^\mathbf{r}_{h} \rightarrow \Om_{h}$ needs to be computed only once and can be reused to gauge the quality of a large number of templates $T_i \in \mathbb{T}$. Furthermore, the linear operator $A$ associated with the harmonic map employed in the computation of a controlmap $\mathbf{r}(T_i): \hOm(\vert F \vert) \rightarrow \hOm^{\mathbf{r}}$ (see \ref{sect:appendix_PAG}) can be reused across multiple faces $F_i$ that share the same parametric domain $\Om_i$. Here, differing control domains $\hOm^\br_i \neq \hOm^\br_j$ of faces $F_i \neq F_j$ with $\hOm_i = \hOm_j$ merely change the right-hand-sides in the linear problem $A x = f$ associated with the control map creation. \\

\noindent In what follows, we provide a cost function for selecting an appropriate template. Given some vertex $\hat{v} \in \mathbb{V}(q), q \in \mathcal{Q}_i$, we denote by $\angle(\hat{v}, q)$ the discrete angle created at $v := \mathbf{x}_{h}^{\br} \circ \br(\hat{v})$ by the two point sets $p_{1} = w^\bxi(\hat{e}_{1})$ and $p_{2} = w^\bxi(\hat{e}_{2})$ with $\{\hat{e}_{1}, \hat{e}_{2} \} \subset E_i$ the two edges incident to $\hat{v}$ in $q \in \mathcal{Q}_i$. Furthermore, by $\theta(\hat{v})$ we denote the \textit{preferred} angle created by two point sets incident to $v = \bx_h^\br \circ \br(\hat{v})$. It is defined as follows:
\begin{align}
\label{eq:preferred_angle_selection}
    \theta(\hat{v}) := \left \{ \begin{array}{ll} 
                    \min(\angle(v, \partial \Om_h), \frac{\pi}{2}) & v \text{ is a boundary vertex} \\
                    \frac{2 \pi}{\operatorname{val}(\hat{v})} & \text{else,} \end{array} \right.
\end{align}
where $\angle(v, \partial \Om_h)$ denotes the interior angle created at the vertex $v = \bx_h^\br \circ\br(\hat{v})$ on $\partial \Om_h$. We note that for boundary vertices $\hat{v} \in V_i$, $\theta(\, \cdot \,)$ tends to favour templates $T_i \in \mathbb{T}$ with $\operatorname{val}(\hat{v}) = 2$ unless the angle created at $v \in \partial \Om_h$ satisfies $\angle(v, \partial \Om_{h}) \geq \tfrac{3}{4} \pi$. As a next step, we define the \textit{angle deviation ratio}:
\begin{align}
    \overline{\angle}(\hat{v}, q) := \frac{\angle(\hat{v}, q)}{\theta(\hat{v})},
\end{align}
which measures the degree of over/undershoot of an interior patch angle from the preferred angle. \\
With $\overline{\angle}(\, \cdot \,, \, \cdot \,)$ at hand, we are in the position to define the following template quality measure (lower means better):
\begin{align}
\label{eq:general_template_selection_criterion}
    C_{T}(T_i) = \frac{\max \limits_{q \in \mathcal{Q}_i} \max \limits_{\hat{v} \in \mathbb{V}(q)} \overline{\angle}(\hat{v}, q)}{\min \limits_{q \in \mathcal{Q}_i} \min \limits_{\hat{v} \in \mathbb{V}(q)} \overline{\angle}(\hat{v}, q)} + \lambda_{\text{patch}} \vert \mathcal{Q}_i \vert.
\end{align}
In words, the criterion from~\eqref{eq:general_template_selection_criterion} favours templates that lead to a layout wherein the ratio between the maximum overshoot and undershoot from the preferred angle is close to unity. Furthermore, if the angle ratio between two templates is similar, the (optional) penalisation via $\lambda_{\text{patch}} \geq 0$ favours the layout with the fewest number of patches. In practice, we utilise $\lambda_{\text{patch}} = 0.5$. \\

\noindent A basic strategy selects the template $T_i \in \mathbb{T}_{\text{select}}$ that minimises~\eqref{eq:general_template_selection_criterion}. \\

\noindent A more sophisticated strategy adds a further optimisation stage to the controlmap $\br: \hOm \rightarrow \hOm^\br$. Let $T_i = (E_i, V_i, \mathcal{Q}_i)$ be an eligible template along with its canonical control template $T^\br = (V^\br, E^\br, \mathcal{Q}^\br)$, which follows from replacing $V_i \ni \hat{v} \rightarrow \br(\hat{v})$ (see \ref{sect:appendix_PAG}). Regarding the controlmap $\br: \hOm \rightarrow \hOm^\br$ as a function of the correspondence $V_i \ni \hat{v} \rightarrow v^\br \in V^\br$, we may further optimise the positions of the $v^\br \in V^\br$ in a way that reduces the value of $C_{T}(T_i)$, see Figure \ref{fig:opt_angle}. We define the preferred angle in a way analogous to~\eqref{eq:preferred_angle_selection}:
\begin{align}
    \theta^\br(v^\br) = \left \{ \begin{array}{ll} 
                    \frac{\angle(v, \partial \Om_h)}{\operatorname{val}(v^\br) - 1} & v^\br \text{ is a boundary vertex} \\
                    \frac{2 \pi}{\operatorname{val}(v^\br)} & \text{else,} \end{array} \right.
\end{align}
where $v = \bx_{h}^\br(v^\br)$ is the image vertex of $v^\br$ under $\bx_{h}^\br: \hOm^\br_{h} \rightarrow \Om_{h}$. By $\angle^\br(v^\br, q^\br)$, we denote the angle on $\Om_{h}$ created on the image of $q^\br \in \mathcal{Q}^\br$ at the vertex $v = \bx_{h}^\br(v^\br)$. Here, the angle is approximated by estimating the tangents associated with edges incident to $v^\br$ using finite differences. We would like to optimise the cost function
\begin{align}
    \widehat{C}_{T^\br}(T^\br) := \frac{\max \limits_{q^\br \in \mathcal{Q}^\br} \max \limits_{v^\br \in \mathbb{V}(q^\br)} \overline{\angle}^\br(v^\br, q^\br)}{\min \limits_{q^\br \in \mathcal{Q}^\br} \min \limits_{v^\br \in \mathbb{V}(q^\br)} \overline{\angle}^\br(v^\br, q^\br)},
\end{align}
where $\overline{\angle}^\br(v^\br, q^\br) := \angle(v^\br, q^\br) / \theta^\br(v^\br)$. Let $\mathbf{h}^\br(V^\br)$ be the vector of all $\overline{\angle}^\br(v^\br, q^\br), \, q^\br \in \mathcal{Q}^\br, \, v^\br \in \mathbb{V}(q^\br)$. To achieve better convergence behaviour, we replace the $\max$ and $\min$ operators by softmax and softmin operators and minimise the logarithm:
\begin{align}
\label{eq:sofmax_min_angle_costfunction}
    C_{T^\br}(T^\br) := \log{\left({\operatorname{softmax}_\beta} \left(\mathbf{h}^\br \right) \right)} - \log{\left({\operatorname{softmin}_\beta} \left(\mathbf{h}^\br \right) \right)},
\end{align}
where 
$$\operatorname{softmax}_\beta \left(\bx \right) := \frac{1}{\sum_{i} \exp(\beta x_i)} \sum \limits_{i} \exp(\beta x_i) x_i$$
and
$$\operatorname{softmin}_\beta \left(\bx \right) := \frac{1}{\sum_{i} \exp(-\beta x_i)} \sum \limits_{i} \exp(-\beta x_i) x_i.$$
In practice, we use $\beta = 6$. \\
The minimisation of~\eqref{eq:sofmax_min_angle_costfunction} is subjected to additional constraints that prevent the minimiser from folding. Let $\mathbf{g}(T^\br)$ be the vector of all cross products $\boldsymbol{\nu}(v^\br, q^\br)$, $q^\br \in \mathcal{Q}^\br$, $v^\br \in \mathbb{V}(q^\br)$ as defined by~\eqref{eq:controlmap_cross_product} and let $\mathbf{g}_0$ be the vector resulting from substituting the template $T^\br$ into $\mathbf{g}(\, \cdot \,)$ before the first iteration. By assumption, we have $\mathbf{g}_0 > \mathbf{0}$. To prevent the controlmap from folding, we optimise
\begin{align}
\label{eq:softmax_min_angle_optimisation}
     & C_{T^\br}(T^\br) \rightarrow \min \limits_{v^\br \in V^\br_{\text{int}}}, \nonumber \\
     \text{s.t.} \quad & \mathbf{g}(T^\br) \geq \mu_{\text{relax}} \mathbf{g}_0,
\end{align}
where $V^\br_{\text{int}}$ denotes the set of interior vertices of $T^\br$. The relaxation parameter $\mu_{\text{relax}} > 0$ tunes the degree of allowed deviation from the initial iterate. In practice, we utilise $\mu_{\text{relax}} = 0.2$. \\

\noindent As before, the optimisation converges in a fraction of a second thanks to the typically low number of inner vertices $v^\br \in V^\br_{\text{int}}$. We furthermore note that computing each internal angle $\angle^\br(v^\br, q^\br), q^\br \in \mathcal{Q}^\br, v^\br \in q^\br$ requires only four function evaluations of $\mathbf{x}_h^\br(\, \cdot \,)$ if the two incident edges are approximated through finite differences. \\

\noindent Since $\bx_h^\br: \hOm^\br_h \rightarrow \Om_h$ is piecewise linear, the gradient of $C_{T^\br}$ with respect to the $v^\br \in V^\br_{\text{int}}$ is defined \textit{almost everywhere} (with the exception of the triangulation's edges). In practice, convergence issues have been notably absent, likely due to the fact that $\bx_h: \hOm^\br_h \rightarrow \Om_h$ approximates a diffeomorphic map, along with the low number of optimisation variables $v^\br \in V^\br_{\text{int}}$.

\section{Face Parameterisation}
\label{sect:parameterisation}
Upon completion of the templatisation stage (see Section \ref{sect:prep_temp}), we proceed with the parameterisation stage. \\
This stage is comprised of two main steps:
\begin{enumerate}
    \item Fitting a spline to each point set $p = w(e)$, with $e \in E$. This results in the weight function $w^S(\, \cdot \,)$ that assigns a spline curve $s \in C^1([0, 1], \mathbb{R}^2)$ to each $e \in E$.
    \item Solving the parameterisation problem $\partial \Omega^S_i \rightarrow \Omega^S_i$, where the piecewise-smooth $\partial \Omega^S_i$ result from concatenating the $s_q = w(e_q)$ with $e_q \in F_i$.
\end{enumerate}
As for point 1., we assume to be in the possession of a (uniform) open base knot vector $\Xi^0$ with
\begin{align}
\label{eq:uniform_initial_knotvector}
    \Xi^0 = \left \{ \underbrace{0, \ldots, 0}_{p + 1 \text{ terms}}, \frac{1}{N_0}, \ldots, \frac{N_0 - 1}{N_0}, \underbrace{1, \ldots, 1}_{p + 1 \text{ terms}} \right \},
\end{align}
wherein $p \geq 2$ denotes the polynomial order while $N_0$ tunes the initial spacing of the knot spans. We note that, albeit uncommon, interior knots may also be repeated a maximum of $p - 1$ times in~\eqref{eq:uniform_initial_knotvector} without violating the assumption that each spline fit $s_q \in C^1([0, 1], \mathbb{R}^2)$. \\

\noindent Given $\Xi^0$ and some $p := w(e) = \{p_1, \ldots, p_{N} \}, \enskip e \in E$, we seek a parameterisation $s := w^S(e)$ by recursively solving the following quadratic optimisation problem:
\begin{align}
\label{eq:quadratic_regularised_least_squares_fit}
    & \min \limits_{s \in \mathcal{V}_\xi^2(\Xi^i)} \quad \frac{1}{N_i} \sum_j \left \| s(\xi_j) - p_j \right \|^2 + \lambda \int \limits_{(0, 1)} \left \| \frac{\partial^2 s}{\partial \xi^2} \right \|^2 \mathrm{d} \xi, \nonumber \\
    & \text{subject to constraints.}
\end{align}
Here, $\mathcal{V}_\xi(\Xi^i)$ denotes the univariate spline space generated by the knot vector $\Xi^i$ after the $i$-th fitting recursion (starting with the user-defined $\Xi^0$) while the parametric abscissae $\{ \xi_j \}_j$ follow from a chord-length parameterisation of the $p_j \in p = \{p_1, \ldots, p_N\}$. The constraints in~\eqref{eq:quadratic_regularised_least_squares_fit} require the spline fit's zeroth and first order derivatives to assume prescribed values in the end points $\xi = 0$ and $\xi = 1$, i.e., $s(0), s^{\prime}(0), s(1)$ and $s^{\prime}(1)$ are prescribed. By default, we impose $s(0) = p_1$ and $s(1) = p_{N}$ in order to avoid mismatches in neighboring edges' vertices. The decision to also impose the values of $s^\prime(0)$ and $s^\prime(1)$ is driven by the desire to leave the division of vertices into angle-forming and non angle-forming (according to the angle threshold $\mu_{\angle} \geq 0$, see Section \ref{sect:prep_temp}) unaffected by the fitting stage.  \\

\noindent Given an edge $e \in E$ with $w(e) = \{p_1, \ldots, p_{N} \}$, we require 
\begin{align}
\label{eq:strong_imposition_tangent_endpoints}
s^\prime(0) = \frac{p_2 - p_1}{\xi_2 - \xi_1} \quad \text{and} \quad s^{\prime}(1) = \frac{p_{N} - p_{N - 1}}{\xi_{N} - \xi_{N - 1}}.
\end{align}
We note that the imposition of~\eqref{eq:strong_imposition_tangent_endpoints} does not avoid the creation of an (albeit small) corner at a vertex $\mathbb{V}(F) \ni v \in V$ on the boundary $\partial \Om^S$ associated with the face $F \in \mathcal{F}$ even if a vertex is modeled as non angle-forming in $F \in \mathcal{F}$. However, if $\mu_\angle$ (see Section \ref{sect:prep_temp}) is small, the minuscule angle can usually be ignored in practice. \\
The regularisation parameter $\lambda \geq 0$ in~\eqref{eq:quadratic_regularised_least_squares_fit} avoids instabilities that may arise if the density of the $p_j \in p$ is (locally) too low for the least-squares fit to be well-conditioned. To assure scaling invariance of $\lambda \geq 0$, we always introduce an equal (in both coordinate directions) coordinate transformation that gauges the point set to unit length, an operation that is reversed upon finalisation of the fitting step(s). In practice, we choose $\lambda = 10^{-5}$. \\

\noindent The residual $r_j := \| s(\xi_j) - p_j \|$ serves as a local refinement criterion. Given $\Xi^i$, if some $r_j$ exceeds a threshold $\mu_{LS}$, we dyadically refine the knot span of $\Xi^i$ that contains $\xi_j$ by adding a new knot in the knot span's center. Doing this for all $r_j$ with $r_j \geq \mu_{LS}$ creates the new knot vector $\Xi^{i+1}$ which serves as input to the next recursion. Aforementioned steps are repeated until all $r_j$ are below a threshold. \\
As a rule of thumb, the regularisation $\lambda \geq 0$ should be chosen at least an order of magnitude smaller than the desired fitting accuracy $\mu_{LS} > 0$. As such, the aimed-for accuracy should satisfy $\mu_{LS} \geq 10^{-4}$ which leads to a fit that is sufficiently accurate in nearly all practical use cases. \\

\noindent Upon completion, we are in the possession of the weight function $w^S(\, \cdot \,)$ which assigns to each $e \in E$ a spline curve $s \in  C^1([0, 1], \mathbb{R}^2)$ generated from the knotvector $\Xi = w^{\Xi}(e)$ which, in turn, follows after the last fitting recursion associated with $p = w(e)$. If $s = w^S(e)$ with $e \in E$, we assume that $w^S(-e) = s(1 - \xi)$. Similarly, we assume that $w^\Xi(\, \cdot \,)$ assigns to $-e$ the coarsest knotvector $-\Xi := w^\Xi(-e)$ which generates the space that contains $s(1 - \xi) = w^S(-e)$. If the sorted unique knot values of $\Xi$ are given by $(0 = \xi_1, \ldots, \xi_n = 1)$, then the unique knot values $\{\xi_i^\prime \}_i$ of $-\Xi$ are given by
$$\xi_i^\prime = 1 - \xi_{\vert \Xi \vert + 1 - i}, \quad \text{where} \quad \vert \Xi \vert \text{ denotes the total number of unique knots.}$$
The collection of spline fits defines a spline face $\Om^S$ from $F \in \mathcal{F}$ as follows:
\begin{align}
    \partial \Om^S = \bigcup \limits_{e \in F} \overline{w^S(e)}.
\end{align}
Thanks to the strongly imposed data at $s(0)$ and $s(1)$, the domain $\hOm^S$ is simply connected. \\

\begin{figure}[h!]
\centering
    \includegraphics[align=c, width=0.95\linewidth]{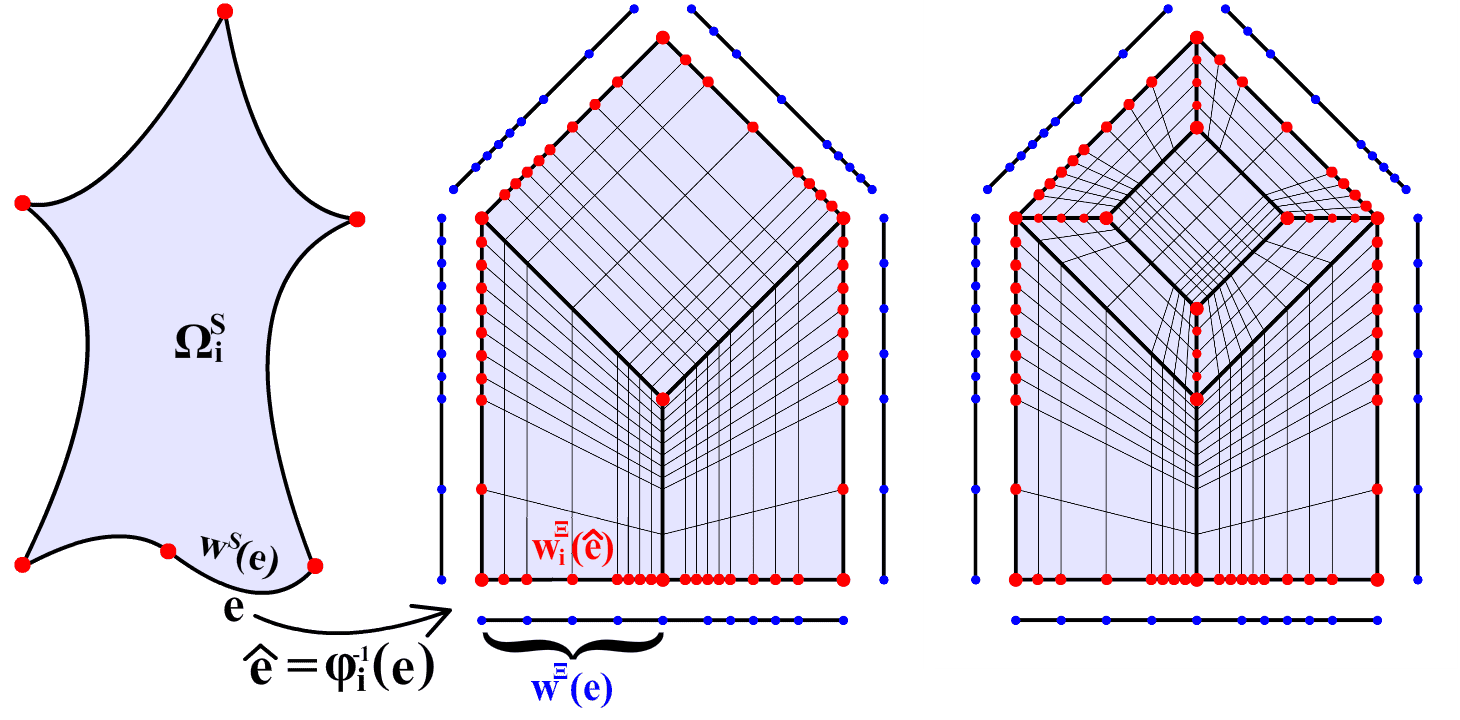}
    \caption{Figure showing how the knotvectors $\Xi = w^\Xi(e)$ (blue) are propagated from the global edge $e \in F_i$ to the associated local edge $\hat{e} = \phi_i^{-1}(e)$ while simultaneously taking the union of opposite knots to create a local weight function $w_i^\Xi(\, \cdot \,)$ (red). The result is a knotvector pair $(\Xi_\mu^j, \Xi_\nu^j)$ for each quad face $q_j \in \mathcal{Q}_i$ / patch $\hOm_j^\square \subset \hOm_i$ which conforms with neighboring cell's knotvectors (center). This results in the coarsest element segmentation that enables the creation of a basis $\mathcal{V}_{h, i} \subset H^1(\hOm_i)$ that is compatible with the boundary correspondence $\mathbf{f}_i: \partial \hOm_i \rightarrow \partial \Om^s_i$ in the sense that the set $\mathcal{U}^{\mathbf{f}_i} \neq \emptyset$. The right figure shows an example of a template $T_i$ where not all interior knotvectors follow from the boundary knotvectors, in which case a suitable knotvector is selected automatically.}
    \label{fig:example_unify_knotvectors}
\end{figure}

\noindent The finalised templatisation stage from Section \ref{sect:prep_temp} yields the graph $G = (V, E, \mathcal{F}, \mathcal{T})$ with $\vert \mathcal{T} \vert = \vert \mathcal{F} \vert$. The next step creates for all $T_i \in \mathcal{T}$ a weight function $w^{\Xi}_i(\, \cdot \,)$, with domain $E_i$ of $T_i = (V_i, E_i, \mathcal{Q}_i)$, that assigns to each $\hat{e} \in E_i$ a knotvector $\hat{\Xi}$. Let $\partial E_i = (\hat{e}_{j_1} \hat{e}_{j_2}, \ldots)$ and $F_i = (e_{j_1}, e_{j_2}, \ldots)$. Thanks to the canonical correspondence $F_i \ni e_{j_q} \rightarrow \hat{e}_{j_q} \in \partial E_i$, it is clear that $w_i^{\Xi}( \, \cdot \, )$ should assign to $\hat{e}_{j_q} \in \partial E_i$ a knotvector whose unique knots are a superset of the unique knots of $w^\Xi(e_{j_q})$. The correspondence $F_i \ni e \rightarrow \hat{e} \in \partial E_i$ furthermore induces a piecewise-smooth spline-based boundary correspondence $\mathbf{f}_i: \partial \hOm_i \rightarrow \partial \Om^S_i$ between each domain-template pair $(\hOm_i, T_i)$ and its corresponding spline domain $\Om^S_i$ via the weight function $w^S(\, \cdot \,)$. \\
Given the quadrangular faces $q \in \mathcal{Q}_i$, the unique knots assigned to the $F_i \ni e \rightarrow \hat{e} \in \partial E_i$ are propagated into the interior as in Figure \ref{fig:example_unify_knotvectors}. Depending on the layout, some interior edges' knotvectors may not follow from the boundary edges. In this case, we assign a uniform knotvector $\Xi^I$ (with the property that $\Xi^I = -\Xi^I$) to the isolated edges. The number of interior knots is chosen based on the average number of knots in all other knotvectors. For an example, see Figure \ref{fig:example_unify_knotvectors}, right.

The figure additionally illustrates how the use of dyadically refined knotvectors avoids the introduction of a dense element subdivision when taking the union of two or more knotvectors. Arbitrarily-placed knots may furthermore lead to minuscule knotspans upon taking the union of several knotvectors, which, in turn, leads to conditioning issues of matrices assembled over the associated spline space. While not immediately evident, we emphasise that potential conditioning issues make the use of dyadically refined knotvectors essentially mandatory in practice. \\

\noindent The result is the weight function $w^\Xi_i( \, \cdot \,)$ that assigns the same knotvector to opposite sides of the $q \in \mathcal{Q}_i$ and therefore a unique knotvector tuple $(\Xi_\mu^j, \Xi_\nu^j)$ to each patch $\hOm^{\square}_j \subset \hOm_i$. \\ 

Let $\mathcal{T} \ni T_i = (V_i, E_i, \mathcal{Q}_i)$ be a template with knotvectors $\Xi = w^\Xi(e), \, \pm e \in E_i$ and associated domains $(\hOm, \Om^S)$. The weight function $w^\Xi(\, \cdot \,)$ assigns a knotvector pair $(\Xi_\mu^j, \Xi_\nu^j)$ to each $q_j \in \mathcal{Q}_i$. Given the reference patch $\Om^{\square} = (0, 1)^2$, the pairs $(\Xi_\mu^j, \Xi_\nu^j)$ induce the canonical spline spaces $\widehat{\mathcal{V}}_{h, j} \subset H^2(\hOm^\square), \enskip j \in \{1, \ldots, \vert \mathcal{Q}_i \vert \}$ with canonical bases $\{ \widehat{\mathcal{V}}_{h, j} \} \subset H^2(\hOm^{\square})$ that follow from the default Cox-de Boor recursion formula \cite{piegl2012nurbs}. Then, we define the basis
\begin{align}
   \left \{ \mathcal{V}_{h}^{\text{disc}} \right \} := \bigcup \limits_{j \in \{1, \ldots, \vert \mathcal{Q}_i \vert \}} \left \{ v_h \circ \left( \mathbf{m}^j \right)^{-1} \, \vert \, v_h \in \{\widehat{\mathcal{V}}_{h, j} \} \right \} \subset L^2(\hOm)
\end{align}
along with
$$
\mathcal{V}_{h}^{\text{disc}} := \operatorname{span} \left \{\mathcal{V}_{h}^{\text{disc}} \right \},
$$
where $\mathbf{m}^j: \Om^\square \rightarrow \hOm_j^{\square}$ is the bilinear diffeomorphic map between $\Om^{\square} = (0, 1)^2$ and the $j$-th patch $\hOm^\square_j \subseteq \hOm$. We note that the canonical basis $\{ \mathcal{V}_{h}^{\text{disc}} \}$ of $\mathcal{V}_{h}^{\text{disc}} \subset L^2(\hOm_i)$ forms a partition of unity on $\hOm$. We may define the subspace $\mathcal{V}_{h} := \mathcal{V}_{h}^{\text{disc}} \cap C^0(\hOm)$, with $\mathcal{V}_{h} \subset H^1(\hOm)$, and the set $\mathcal{U}^\mathbf{F}_{h} := \{ \boldsymbol{\phi}_h \in \mathcal{V}_{h}^2 \, \, \vert \, \, \boldsymbol{\phi}_h = \mathbf{F} \text{ on } \partial \hOm \}$. By construction, the set $\mathcal{U}^{\mathbf{f}}$, that takes $\mathbf{F}$ as the spline-based boundary correspondence $\mathbf{f}: \partial \hOm \rightarrow \partial \Om^S$, satisfies $\mathcal{U}^{\mathbf{f}} \neq \emptyset$. \\

\noindent With the set $\mathcal{U}^{\mathbf{f}}_{h}$ at hand, we are in the position to compute a spline-based parameterisation for $\Om^S$. For simplicity, for the time being, we assume that the associated control surface satisfies $\hOm^\br = \hOm$, i.e. $\br(\bxi) = \bxi$ as this simplifies the discussion that follows. Generalisations to $\hOm^\br \neq \hOm$ are made afterwards. \\

Computing the discrete harmonic map $\mathbf{x}_{h} \in \mathcal{U}^{\mathbf{f}}_{h}$ leads to a nonlinear root-finding problem. Approach 1. (see Section \ref{sect:harmonic_maps}) follows the methodology outlined in \cite{hinz2024use} which approximates the inverse of a harmonic map between $\Om^S$ and $\hOm$ by seeking the root of the following discretised variational formulation of a quasilinear elliptic PDE in nondivergence form:
\begin{align}
\label{eq:C0_DG_discrete}
    \text{find } \mathbf{x}_{h} \in \mathcal{U}_{h}^{\mathbf{f}} \quad \text{s.t.} \quad \mathcal{L}_{\eta}^\mu(\mathbf{x}_h, \boldsymbol{\phi}_h) = 0 \quad \forall \boldsymbol{\phi}_h \in \mathcal{U}_{h}^{\mathbf{0}},
\end{align}
in which
\begin{align}
\label{eq:C0_DG_discrete_operator}
    \mathcal{L}_{\eta}^\mu(\mathbf{x}, \boldsymbol{\phi}) := \sum \limits_{\hOm^\square_j \subset \hOm} \int \limits_{\hOm^{\square}_j} \gamma^\bx \Delta \phi_{k} A^\mu \colon H(x_k) \mathrm{d} \bxi + \eta \sum \limits_{\gamma \in \Gamma^I} \frac{1}{h(\gamma)} \int \limits_{\gamma} [\![ \nabla x_k ]\!] \, \colon \, [\![ \nabla \phi_k ] \! ] \mathrm{d} \bxi,
\end{align}
and we sum over repeated indices. The matrix-valued function $A^\mu = A^\mu(\bx)$ satisfies $A^\mu(\bx) = \operatorname{cof}(G) + \mu \mathcal{I}^{2 \times 2}$ with $G(\bx)$ the metric tensor induced by the map $\bx: \hOm \rightarrow \Om^S$ and $\operatorname{cof}(C)$ the cofactor matrix of $C$. It introduces the nonlinearity in~\eqref{eq:C0_DG_discrete}. Meanwhile, $\gamma^\bx(\bx) := (A^\mu \colon A^\mu)^{-1} \operatorname{trace}(A^\mu)$ with $B: C$ the Frobenius inner product between two matrices while $H(u)$ denotes the Hessian matrix of $u: \hOm \rightarrow \mathbb{R}$ with $H(u)_{ij} = \partial_{\bxi_i} \partial_{\bxi_j} u$. The expression, $[ \! [ \, \boldsymbol{\phi} \, ] \! ]$ returns the (entry-wise) jump term of $\boldsymbol{\phi} \otimes \mathbf{n} \in L^2(\gamma, \mathbb{R}^{2 \times 2})$, with $\mathbf{n}$ the unit outer normal on $\gamma$ in arbitrary but fixed direction. Taking the stabilisation parameter $\mu > 0$ ensures that $A^\mu(\bx)$ is uniformly elliptic (a.e. in $\hOm$) even for intermediate iterates which may not satisfy $\det G(\bx) > 0$. In practice, we take $\mu = 10^{-5}$ and for achieving scaling invariance, both $\Om^S$ and $\hOm^\br$ are gauged to unit surface area through a coordinate transformation, an operation that is reversed upon finalisation of the spline map.\\
The $\hOm_j^\square \subset \hOm$ refer to the patches that follow from the quad faces $q \in \mathcal{Q}_i$ under the template $T_i \in \mathcal{T}$ while $\eta > 0$ is a penalty parameter. In practice, we choose $\eta = 10$. The set $\Gamma^I$ is defined as the collection of interior interfaces between the patches $\hOm_j^\square \subset \hOm$ and follows from the $\hat{e} \in (E_i \setminus \partial E_i)$. Finally, $h(\gamma)$ refers to the average diameter of all knotspans over the interior interface $\gamma \in \Gamma^I$. \\ 

\noindent The nonlinearity is tackled using a Newton algorithm with line search and the iterative algorithm is initialised with a map $\bx_{h}^0 \in \mathcal{U}^{\mathbf{f}}$ that is harmonic in $\hOm$ (rather than $\Om^S)$. For details, we refer to \cite{hinz2024use}. \\

\noindent Each template $T_i \in \mathcal{T}$ that was selected in the templatisation stage (see Section \ref{sect:prep_temp}) is accompanied by a control surface $\hOm^{\br}$ along with the corresponding controlmap $\br: \hOm \rightarrow \hOm^{\br}$ which allows a harmonic map $\mathbf{x}^{-1}: \hOm^S \rightarrow \hOm^{\br}$ to extend diffeomorphically to the closure $\overline{\Om}^S$ of $\Om^S$. By assumption, the map $\mathbf{r}: \hOm \rightarrow \hOm^{\br}$ is patchwise smooth with a bounded inverse. For $\hOm \neq \hOm^\br$ we may reinterpret~\eqref{eq:C0_DG_discrete} as a problem posed over $\hOm^\br$ and use basic differential geometry identities to reformulate the integrals over $\hOm$ through a pull-back. While the pullback changes the associated differential operators in $\hOm$, it does not change the nature of the equation and the same basic approach can be employed. For details, see \cite{hinz2024use}. \\

\begin{figure}[h!]
\centering
    \begin{subfigure}[b]{0.45\textwidth}
        \centering
        \includegraphics[align=c, width=0.95\linewidth]{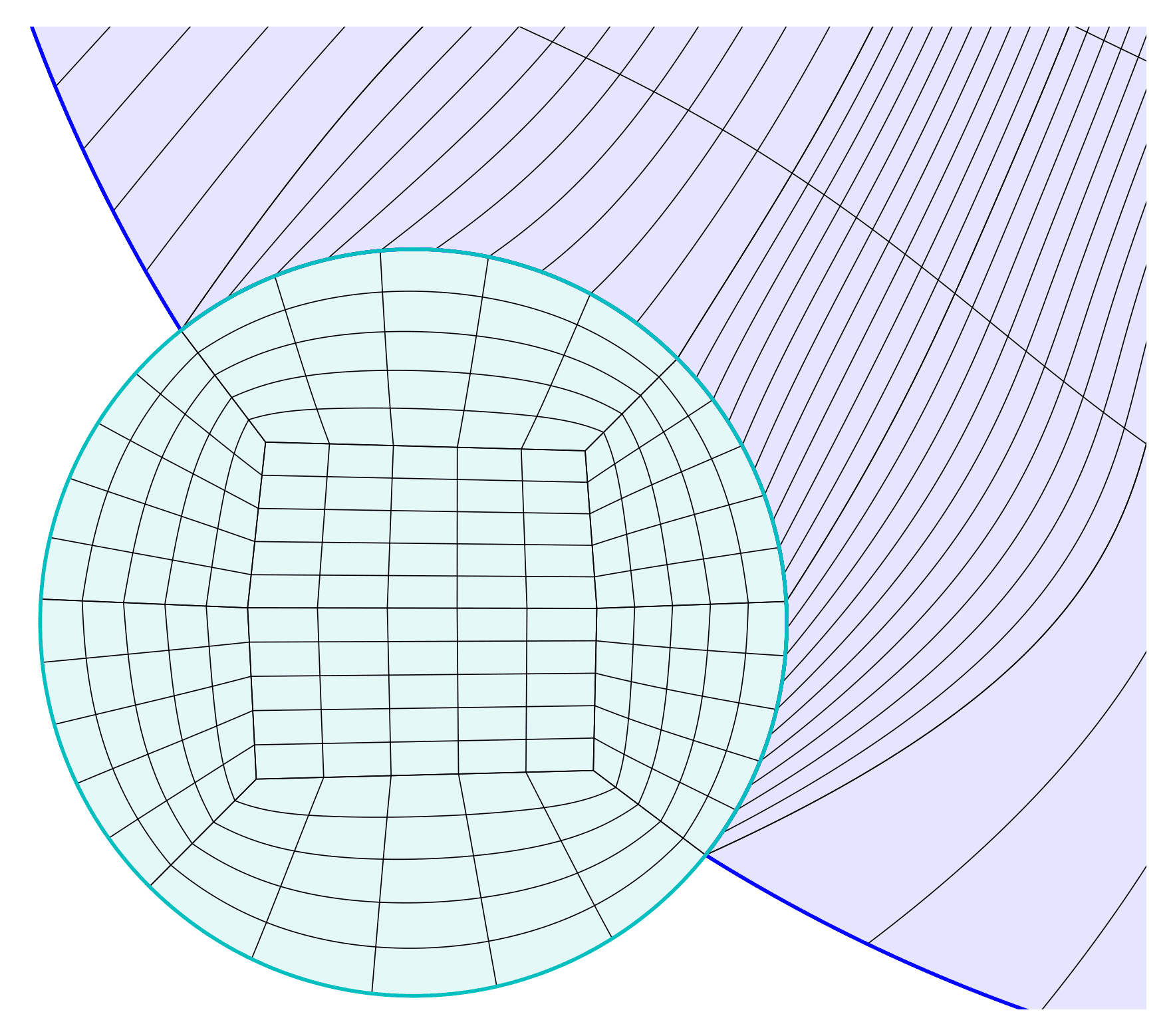}
        \caption{}
        \label{fig:expample_prolong_not_unified}
    \end{subfigure}
    \begin{subfigure}[b]{0.45\textwidth}
        \centering
        \includegraphics[align=c, width=0.95\linewidth]{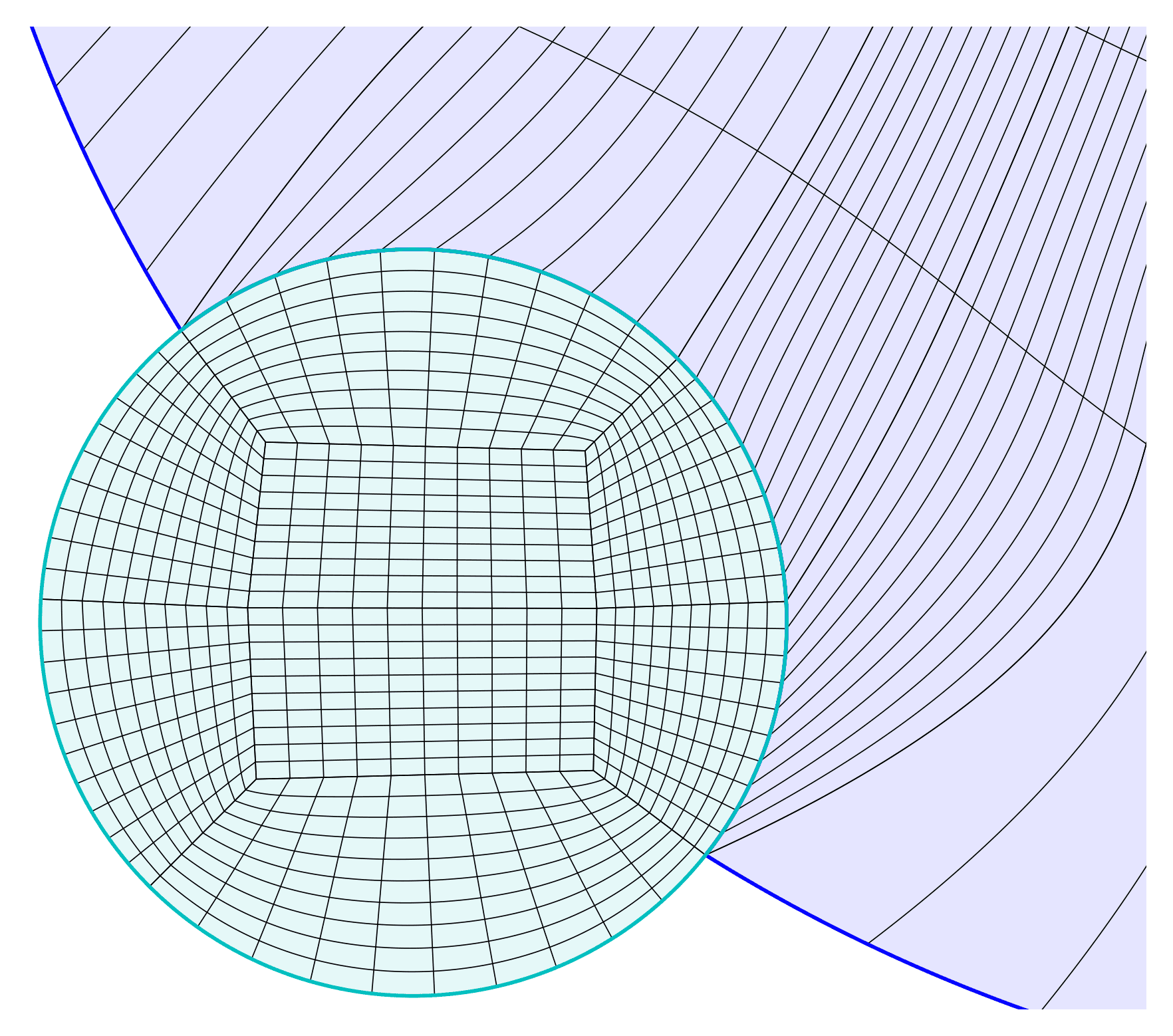}
        \caption{}
        \label{fig:expample_prolong_unified}
    \end{subfigure}
\caption{Figure showing the interface between two faces in differing colors. While the interface is watertight and parameterically conforming in \textbf{(a)}, the splines are represented in differing bases on either sides of the interface which can be seen from the differing knot span densities shown in the plot. In \textbf{(b)} the coarser face has been prolonged to match the density of its neighboring face without changing the parameterisation itself.}
\label{fig:example_prolong}
\end{figure}

\noindent Upon solving~\eqref{eq:C0_DG_discrete} for each face $F_i \in \mathcal{F}$, we are in the possession of $\vert \mathcal{F} \vert$ maps $\mathbf{x}_{h, i}: \hOm_i \rightarrow \Om^S_i$. The union of the maps parameterises the entire plane graph. The interfaces
\begin{align}
    \Gamma^{\bx} := \left \{\gamma_{ij}^\bx \, \, \vert \, \, \gamma_{ij}^\bx := \operatorname{Int} \left(\overline{\Om}^S_i \cap \overline{\Om}^S_j \right) \text{ is a (collection of) open line segment(s)} \right \},
\end{align}
are given by the union of the images of the $s = w^S(e), \, e \in E$. As such, they are fully conforming, both in the sense of watertightness as well as parameterically. We note, however, that if two faces $F_i$ and $F_j$ have shared edges, the $s = w^S(e)$, with $\pm e \in F_i$ and $\mp e \in F_k$, may be expressed in differing supersets of the canonical knotvector $\Xi = w^\Xi(e)$ of $s$ under the two maps $\mathbf{x}_{h,i}$ and $\mathbf{x}_{h, j}$ (see Figure \ref{fig:expample_prolong_not_unified}). More precisely, let $e \in F_i$ with $-e \in F_j$ be a shared edge. Furthermore, let $\hat{e}_1$ and $-\hat{e}_2$ be the associated local edges under the correspondences $\phi_i^{-1}: F_i \rightarrow \partial E_i$ and $\phi_j^{-1}: F_j \rightarrow \partial E_j$, respectively. Then, generally, $w^\Xi_j(\hat{e}_1) \neq w^\Xi_k(-\hat{e}_2)$. If desired, this can be computationally inexpensively remedied by prolonging neighboring faces to a unified element segmentation, as in Figure \ref{fig:expample_prolong_unified}. 

\section{Results}
\label{sect:results}
In this section, we apply the methodology developed in Sections \ref{sect:prep_temp} and \ref{sect:parameterisation} to select geometries. In all cases, the input plane graphs are preprocessed using the techniques from Section \ref{sect:prep_temp} after which we utilise the techniques from Section \ref{sect:parameterisation} to parameterise each individual face. In case a face's parameterisation is singular due to the truncation associated with~\eqref{eq:C0_DG_discrete}, we remedy this by employing the untangling routines detailed in Section \ref{sect:post_processing}. \\

\noindent In what follows, strategy \textbf{1} refers to the template selection based on the heuristic minpatch approach which, after the mandatory pre-filter from Algorithm \ref{algo:template_pre_filter}, retains layouts that (if applicable) align with the face's symmetries and then selects a template with the lowest number of required patches at random. Meanwhile, the deterministic selection based on~\eqref{eq:general_template_selection_criterion} without and with a posteriori optimisation based on~\eqref{eq:softmax_min_angle_optimisation} is referred to as strategy \textbf{2} and \textbf{3}, respectively. \\

\noindent As a basic first example, we are considering the geometry from Figure \ref{fig:minpach_strategy_symmetric} \textbf{(a)}. 

\begin{figure}[h!]
\centering
    \begin{subfigure}[b]{0.45\textwidth}
        \begin{tikzpicture}
        \node[anchor=south west,inner sep=0] (image) at (0,0) {\includegraphics[align=c, width=0.95\linewidth]{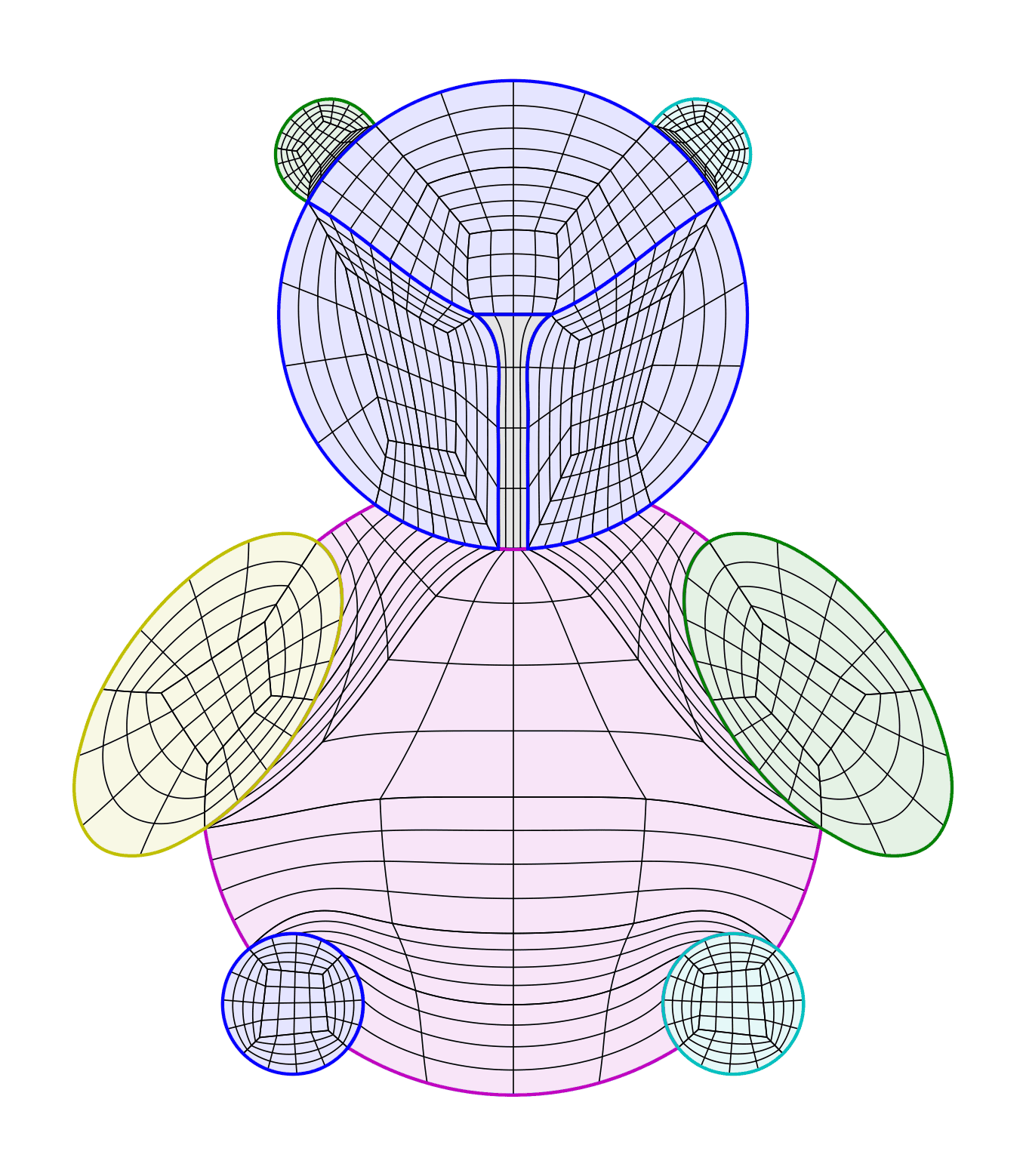}};
        \begin{scope}[x={(image.south east)},y={(image.north west)}]
        \node[anchor=center] at (0.52, 0.31) {\color{blue}\textbf{A}}; 
        \end{scope}
        \end{tikzpicture}
        \caption{}
    \end{subfigure}
    \begin{subfigure}[b]{0.45\textwidth}
        \begin{tikzpicture}
        \node[anchor=south west,inner sep=0] (image) at (0,0) {\includegraphics[align=c, width=0.95\linewidth]{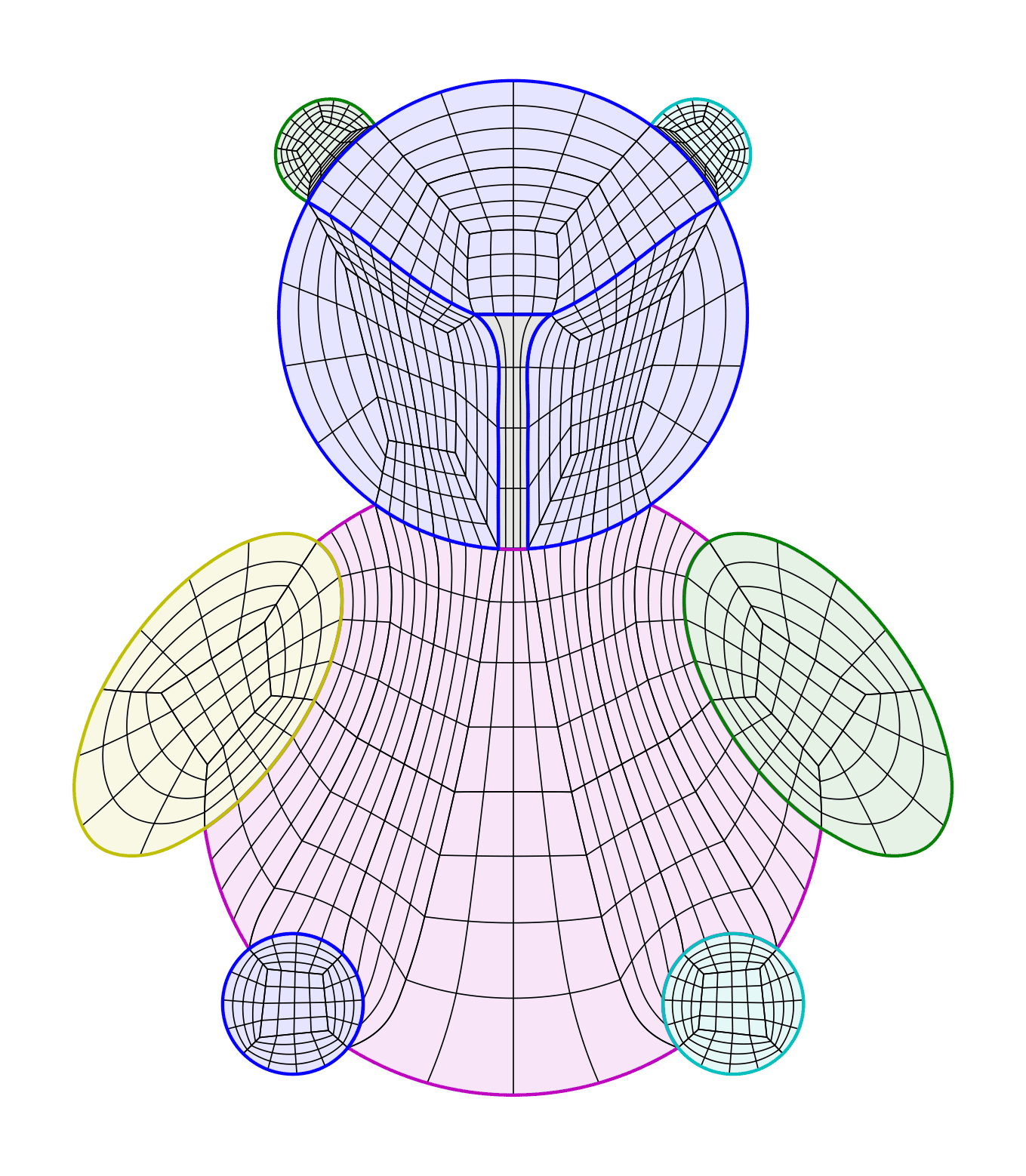}};
        \begin{scope}[x={(image.south east)},y={(image.north west)}]
        \node[anchor=center] at (0.52, 0.3) {\color{blue}\textbf{A}}; 
        \end{scope}
        \end{tikzpicture}
        \caption{}
    \end{subfigure}
\caption{The result of applying strategies $\textbf{1}$ and $\textbf{2}$ to an elementary input plane graph.}
\label{fig:results_connie}
\end{figure}

Applying the heuristic strategy \textbf{1}, we obtain the TSG depicted in Figure \ref{fig:minpach_strategy_symmetric} \textbf{(c)} and the finalised spline-based parameterisation of the input graph in Figure \ref{fig:results_connie} \textbf{(a)}. Figure \ref{fig:results_connie} \textbf{(b)} shows the finalised parameterisation after selection based on strategy \textbf{2}. Both strategies yield similar results with the exception of the face marked by the letter \textit{A}, for which strategy \textbf{2} selects a template that results in a parameterisation with notably improved uniformity. However, both spline parameterisations' faces are all bijective without the need for untangling. \\

\begin{figure}[h!]
\centering
\begin{subfigure}[b]{0.45\textwidth}
    \includegraphics[align=c, width=0.9\linewidth]{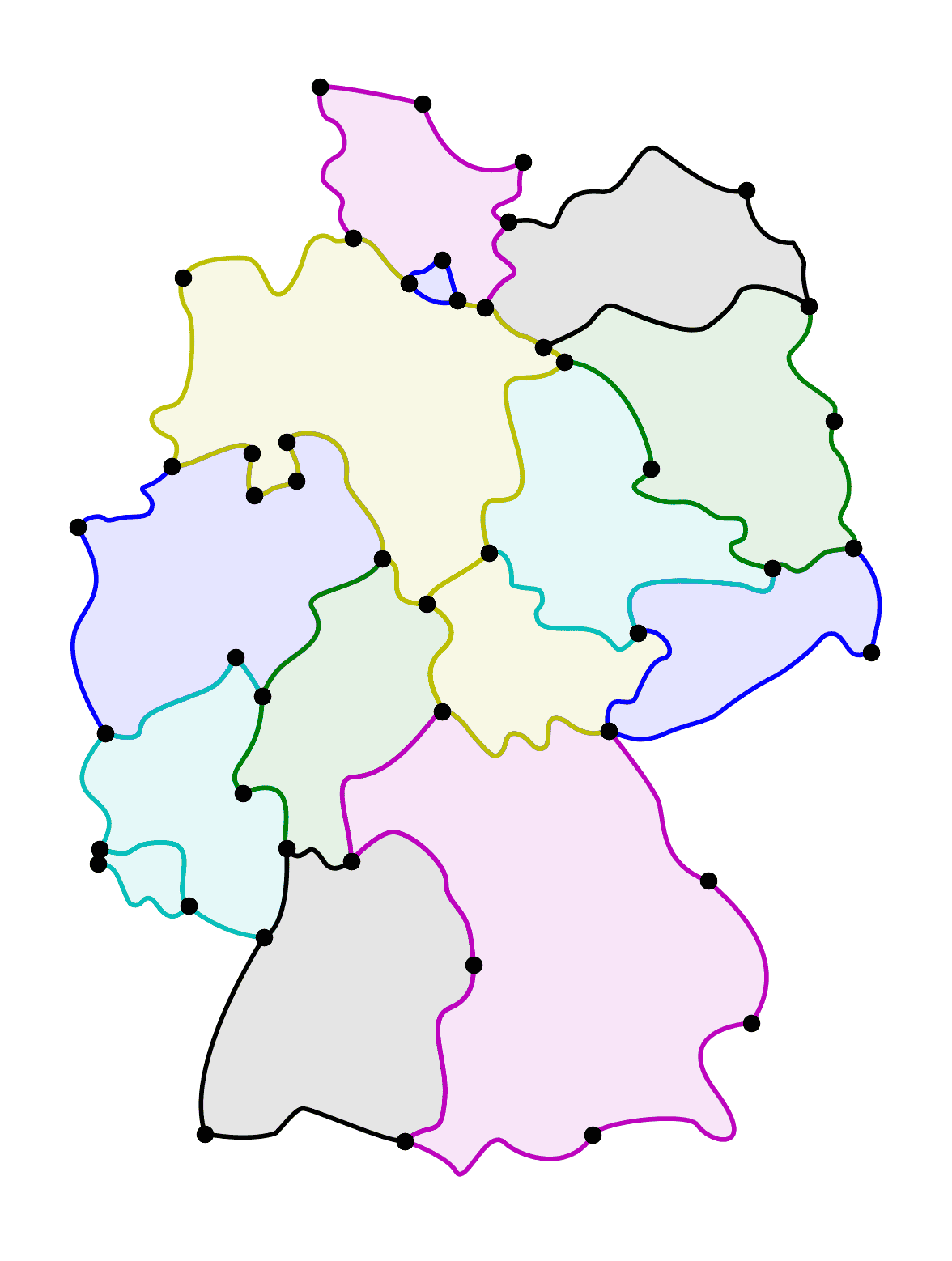}
    \caption{}
\end{subfigure}
\begin{subfigure}[b]{0.45\textwidth}
        \begin{tikzpicture}
            \node[anchor=south west,inner sep=0] (image) at (0,0) {\includegraphics[align=c, width=0.9\linewidth]{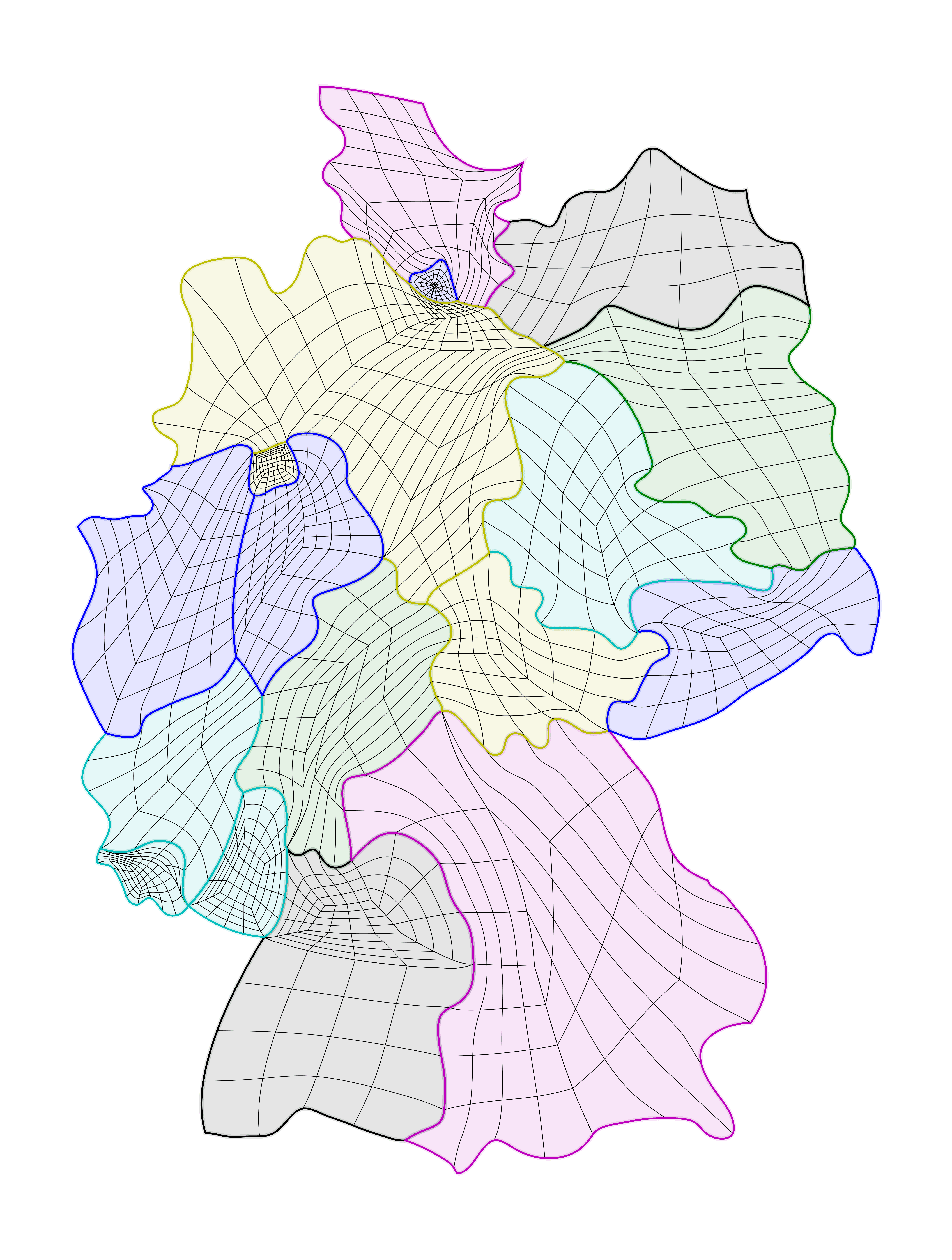}};
            \begin{scope}[x={(image.south east)},y={(image.north west)}]
                \node[anchor=center] at (0.55, 0.29) {\color{blue}\textbf{A}}; 
                \node[anchor=center] at (0.445, 0.85) {\color{blue}\textbf{B}}; 
            \end{scope}
        \end{tikzpicture}
    \caption{}
\end{subfigure} \\
\begin{subfigure}[b]{0.45\textwidth}
    \includegraphics[align=c, width=0.9\linewidth]{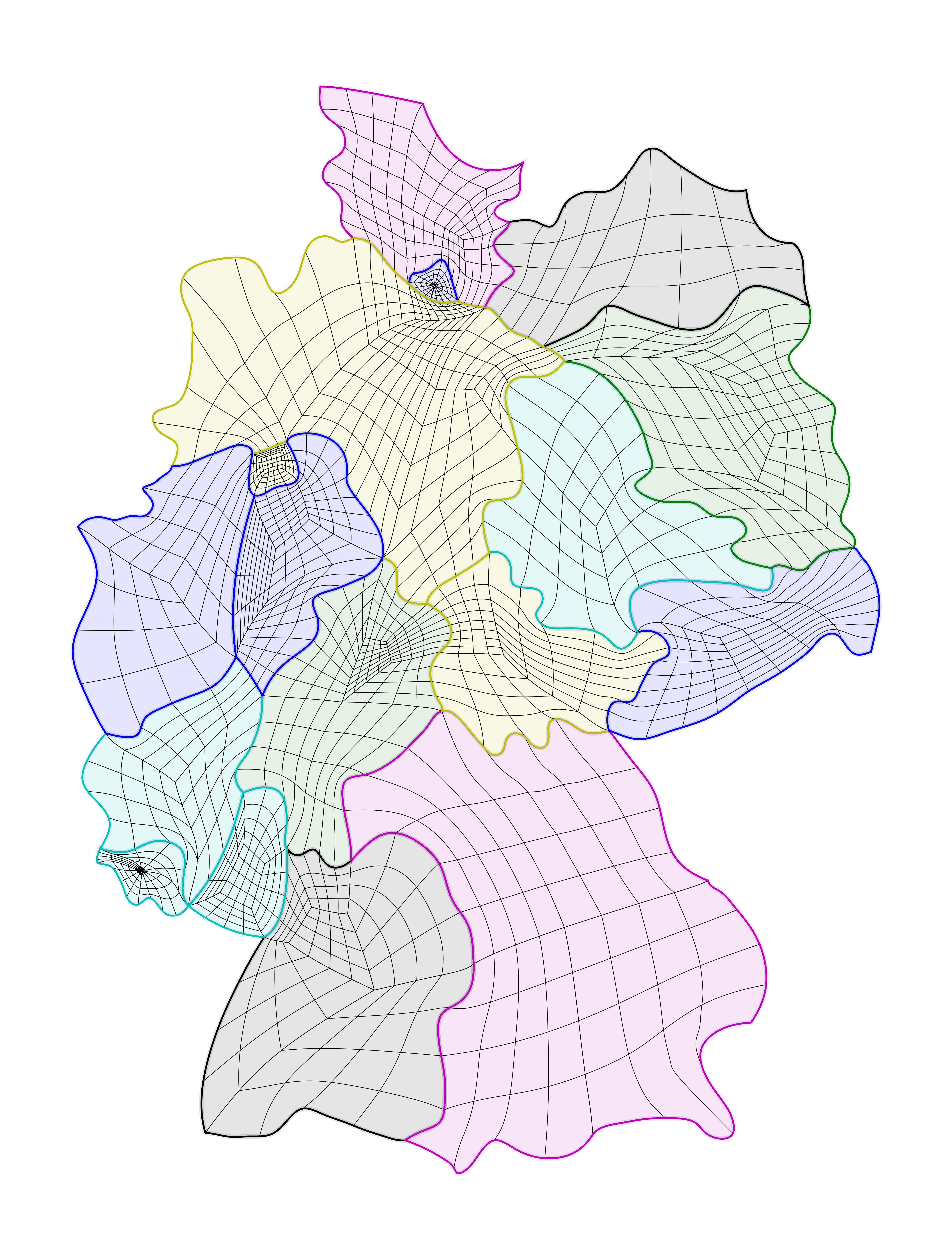}
    \caption{}
\end{subfigure}
\begin{subfigure}[b]{0.45\textwidth}
    \includegraphics[align=c, width=0.9\linewidth]{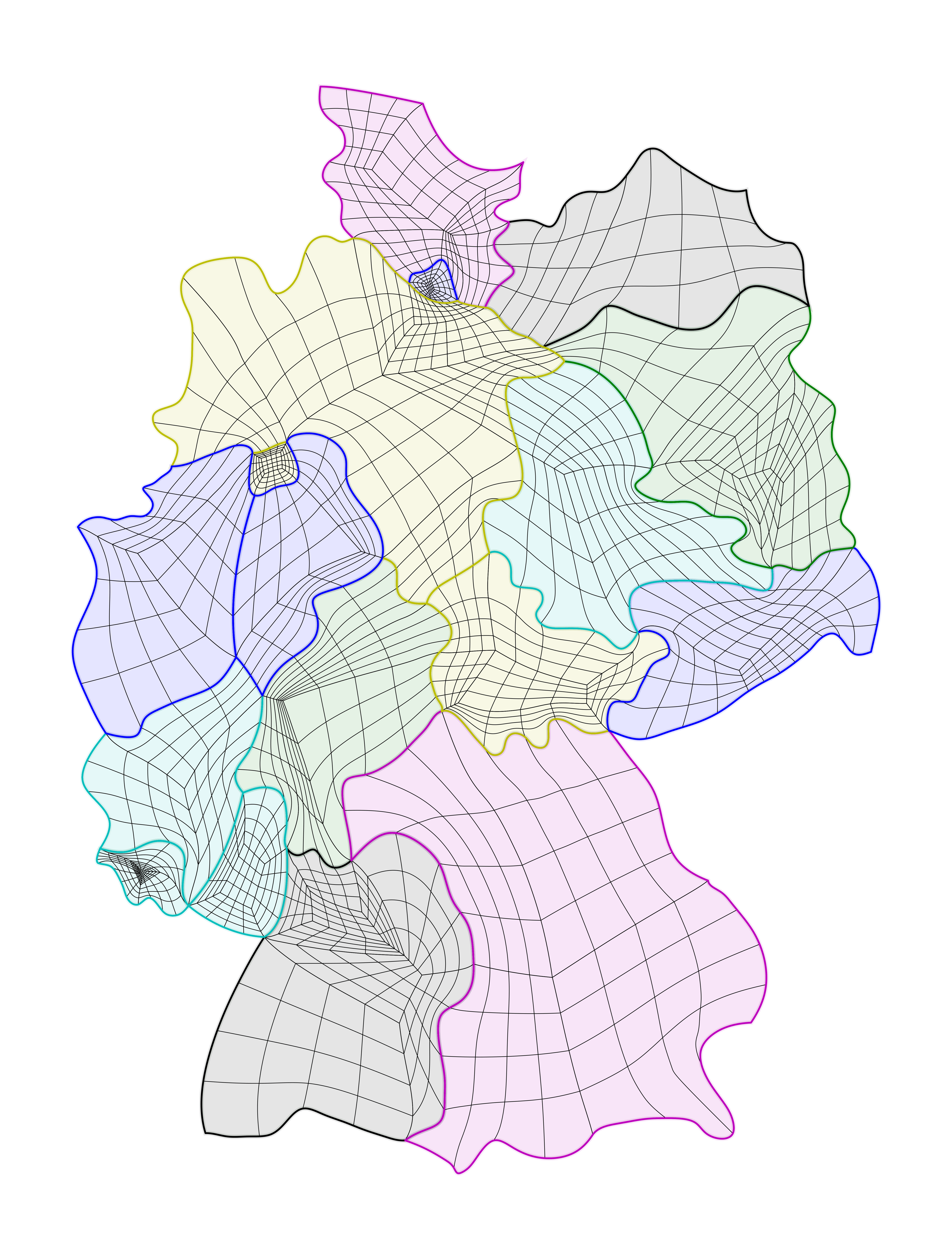}
    \caption{}
\end{subfigure}
\caption{Spline-based parameterisation of all German provinces (excluding enclaves) using template selection strategies \textbf{1} - \textbf{3}.}
\label{fig:results_germany}
\end{figure}

\noindent Generally, the heuristic minpatch strategy \textbf{1} from Figure \ref{fig:results_connie} \textbf{(a)} is an inexpensive method that often suffices for simple plane graphs. However, the same strategy may not align with the complex requirements of more advanced input graphs. As an example, we are considering the plane graph from Figure \ref{fig:results_germany} \textbf{(a)} whose faces constitute a rough approximation of all German provinces (with the exception of enclaves, since multiply-connected domains are beyond this paper's scope). \\

\noindent Figures \ref{fig:results_germany} \textbf{(b)}-\textbf{(d)} show the finalised spline-based parameterisations of the input plane graph for strategies \textbf{1}, \textbf{2} and \textbf{3}. All approaches lead to a spline parameterisation in which all faces are bijective. However, the minpatch strategy \textbf{1} required untangling on the faces marked by the letters \textit{'A'} and \textit{'B'} (see Section \ref{sect:post_processing}). \\

\begin{figure}[h!]
\centering
\begin{subfigure}[b]{0.32\textwidth}
    \includegraphics[width=0.85\linewidth]{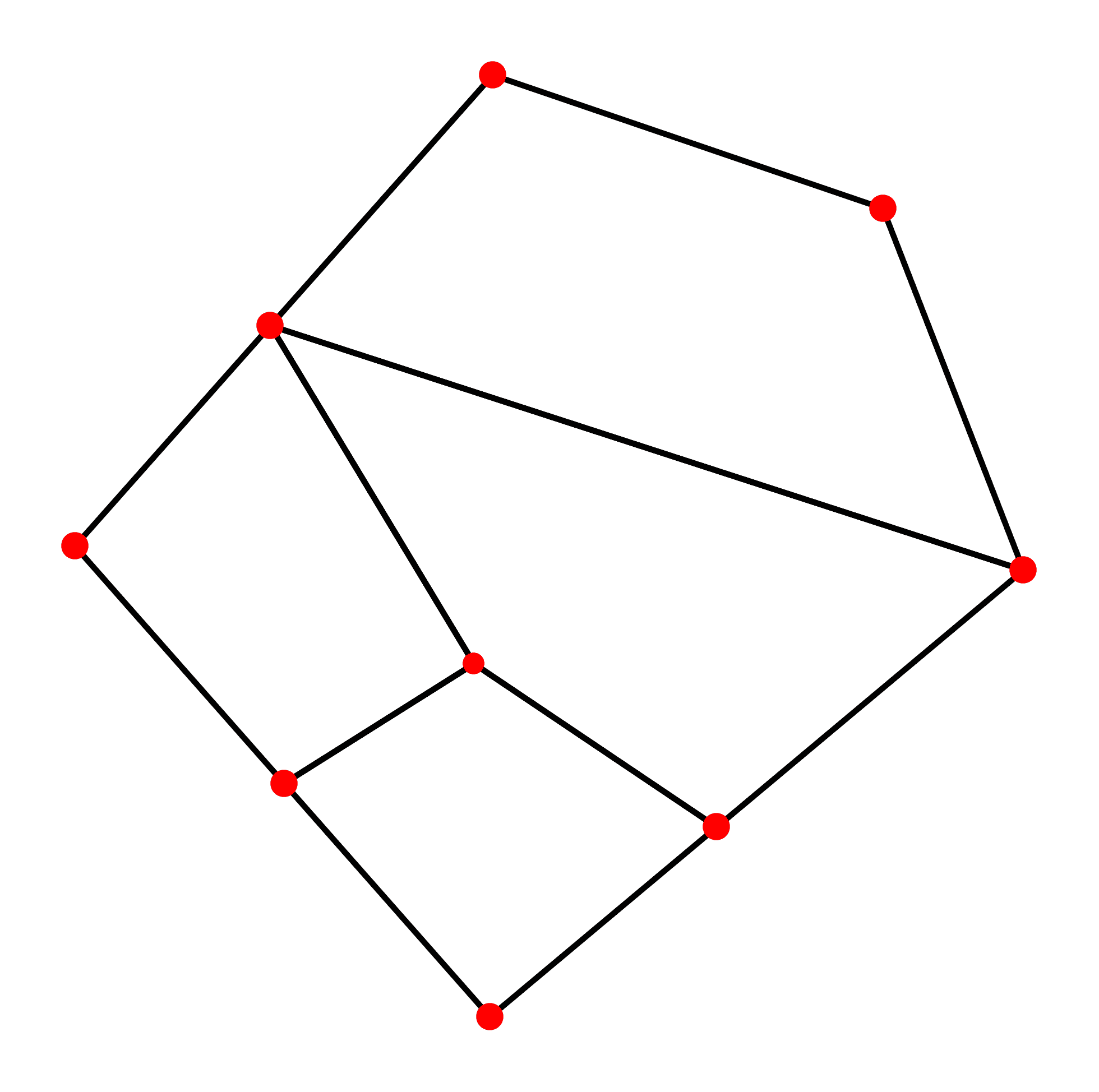}
    \caption{}
\end{subfigure}
\begin{subfigure}[b]{0.32\textwidth}
    \includegraphics[width=0.85\linewidth]{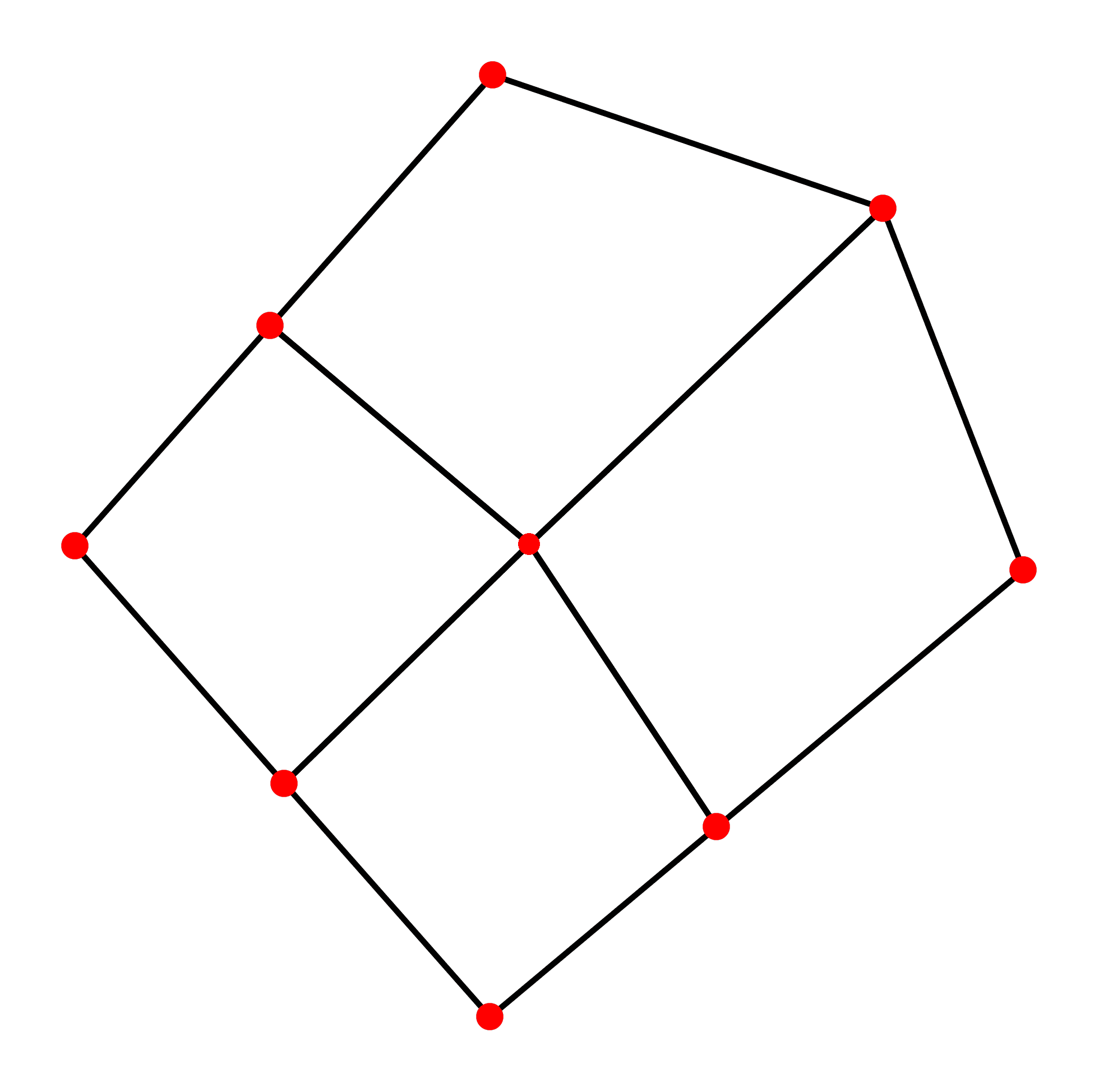}
    \caption{}
\end{subfigure}
\begin{subfigure}[b]{0.32\textwidth}
    \includegraphics[width=0.85\linewidth]{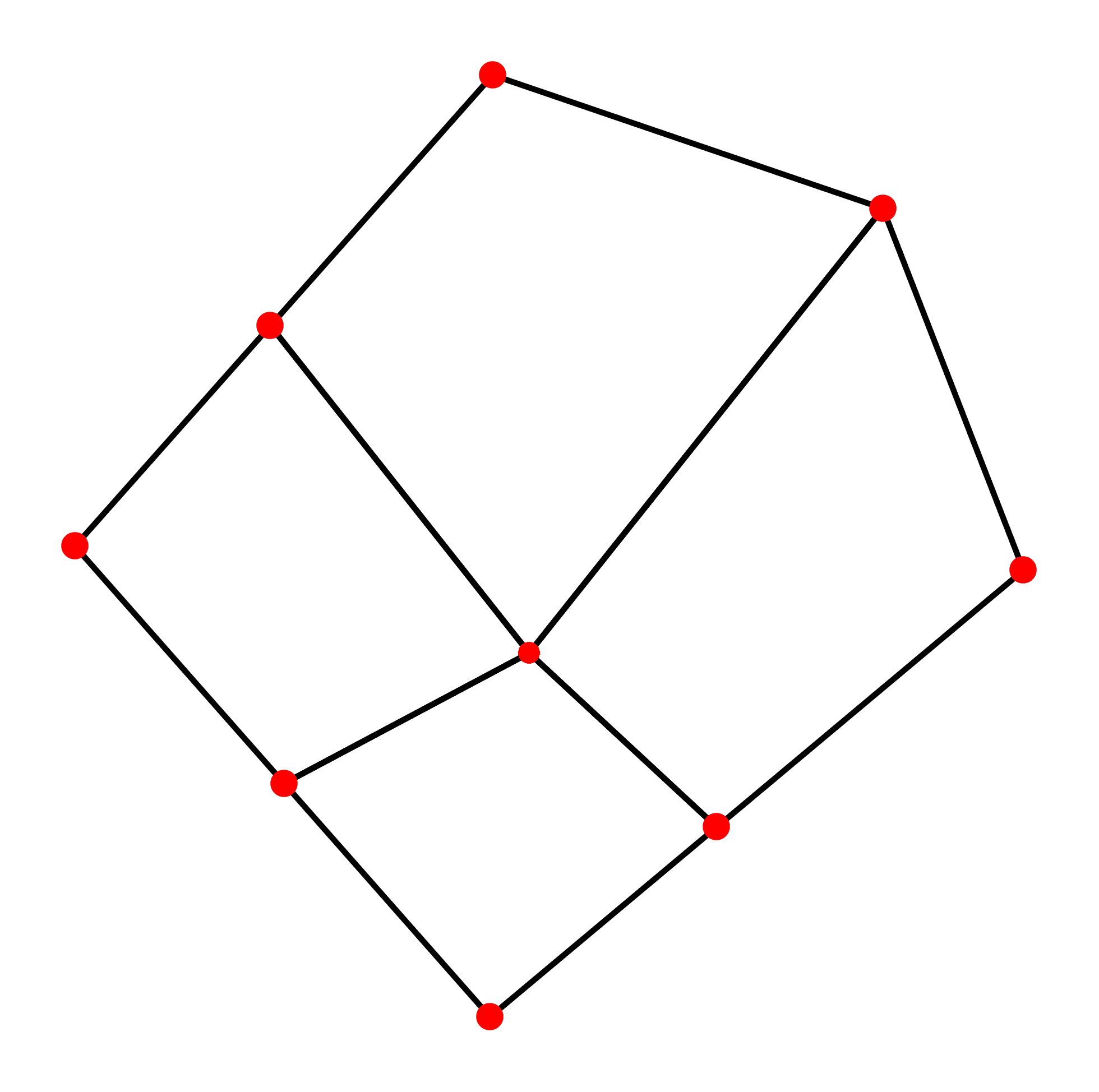}
    \caption{}
\end{subfigure}
\caption{Figure showing the multipatch coverings chosen by strategies \textbf{1} - \textbf{3} of the control domain $\hOm^\br$ associated with the face marked by the letter \textit{'A'} in Figure \ref{fig:results_germany} \textbf{(b)}.}
\label{fig:germany_contromap_5}
\end{figure}

\noindent Besides that, Figures \ref{fig:results_germany} \textbf{(c)} and \textbf{(d)} both tend to select templates which lead to face parameterisations with a notably improved uniformity. Figure \ref{fig:germany_contromap_5} depicts the control domain $\hOm^\br$ along with its covering selected by strategies \textbf{1} - \textbf{3} for the face marked by \textit{'A'}. All approaches select a template with the face's minimum number of four required patches. However, both strategy \textbf{2} and \textbf{3} select a template whose structure aligns more seamlessly with the contour of face \textit{'A'}. Since strategy \textbf{1}, selects a template with the lowest number of required patches at random, the choice's appropriateness, unlike strategies \textbf{2} and \textbf{3}, is contingent on a stochastic process. Compared to strategy \textbf{2}, approach \textbf{3} additionally optimises the position of the template's inner vertex in an effort to minimise~\eqref{eq:softmax_min_angle_optimisation}. \\

\noindent We conclude that besides potentially improved parametric quality, strategies \textbf{2} and \textbf{3} furthermore require less a posteriori untangling in practice. On the other hand, the heuristic minpatch strategy \textbf{1} is algorithmically more lightweight, leading to speedups by the factors $\sim 8$ and $\sim 12$ during the template selection stage with respect to strategies \textbf{2} and \textbf{3}, respectively. \\

\begin{figure}[h!]
\begin{subfigure}[b]{0.45\textwidth}
        \begin{tikzpicture}
            \node[anchor=south west,inner sep=0] (image) at (0,0) {\includegraphics[align=c, width=0.95\linewidth]{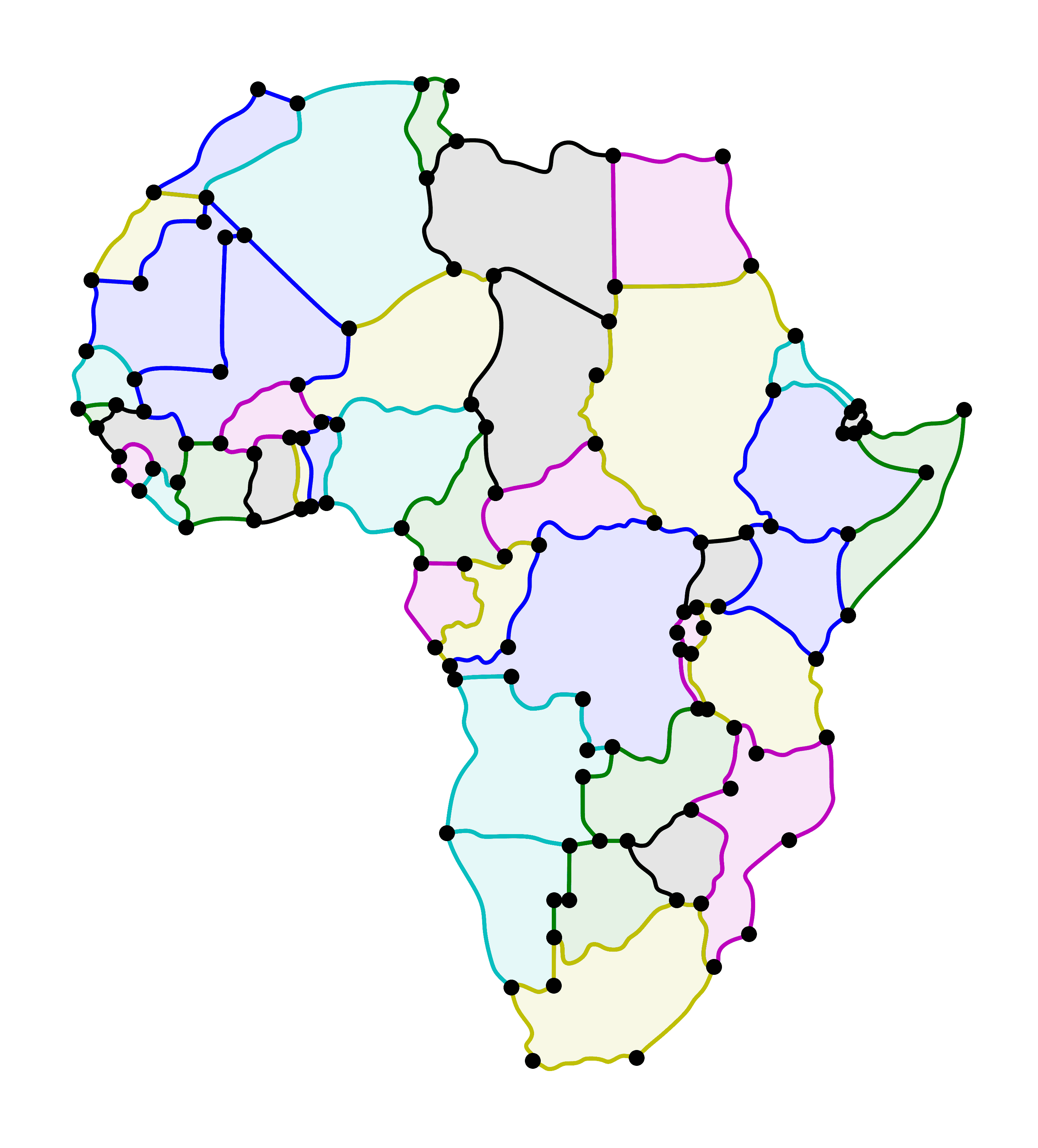}};
            \begin{scope}[x={(image.south east)},y={(image.north west)}]
                \node[anchor=center] at (0.48, 0.48) {\color{blue}\textbf{A}}; 
            \end{scope}
        \end{tikzpicture}
    \caption{}
\end{subfigure}
\begin{subfigure}[b]{0.45\textwidth}
    \includegraphics[width=0.95\linewidth]{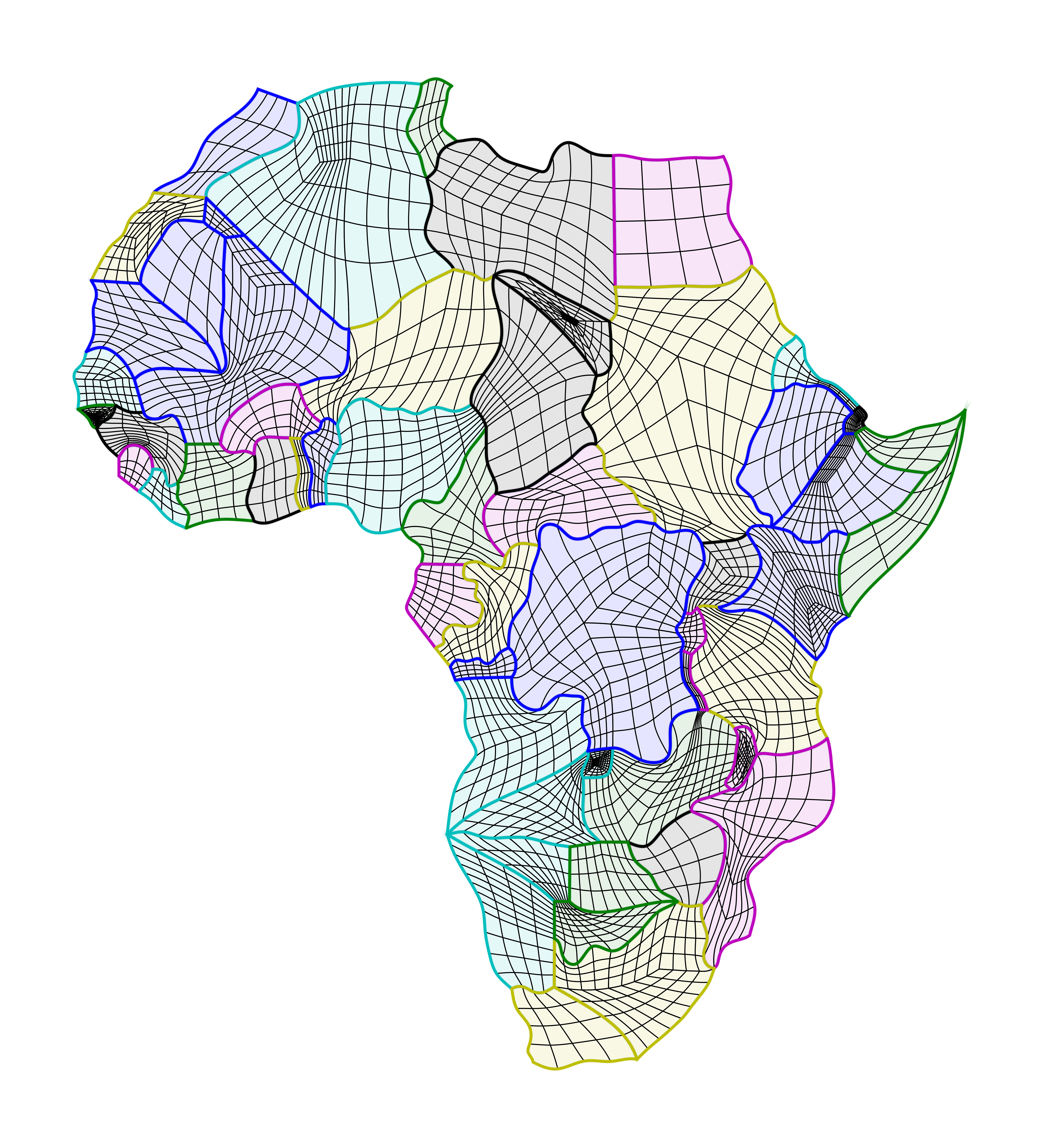}
    \caption{}
\end{subfigure}
\caption{Spline-based parameterisation of a plane graph representing the African continent with template selection based on~\eqref{eq:general_template_selection_criterion} and template optimisation based on~\eqref{eq:softmax_min_angle_optimisation}.}
\label{fig:results_africa}
\end{figure}

\noindent As a final example, Figure \ref{fig:results_africa} shows a plane graph $G = (V, E, \mathcal{F})$ representing the African continent. Figure \ref{fig:results_africa} \textbf{(b)} shows the spline-based parameterisation resulting from applying template selection strategy \textbf{3}. The algorithm successfully removed all concave corners by splitting faces using a total of $13$ Hermite curves.
Of the resulting $55$ faces, all lead to a bijective outcome whereby three faces required additional untangling using the techniques from Section \ref{sect:post_processing}.
\section{Post-Processing}
\label{sect:post_processing}
The purpose of this stage is two-fold:
\begin{enumerate}
    \item Untangling folded parameterisations;
    \item (Optionally) further optimising the parametric properties.
\end{enumerate}
Due to RKC theory, Theorem \ref{thrm:RKC}, folded parameterisation are a direct consequence of the truncation error associated with the scheme from~\eqref{eq:C0_DG_discrete}. Clearly if a face $F \in \mathcal{F}$ leads to a map $\mathbf{x}: \hOm \rightarrow \Om^S$ that is folded, recomputing $\mathbf{x}: \hOm \rightarrow \Om^S$ under a refined spline space may provide relief. However, since this paper employs structured spline spaces, refinements performed on a face $F \in \mathcal{F}$ have to be extended to neighboring faces in order to maintain equal knotvectors on shared edges leading to a DOF proliferation that may not be desirable. \\
As an alternative to refinement, we present an untangling strategy based on the regularised Winslow function \cite{turner2018curvilinear, wang2021smooth}. As before, we assume $\br(\bxi) = \bxi$ for convenience. Generalisations to $\hOm \neq \hOm^\br$ follow as before. A harmonic map is the unique minimiser of the so-called Winslow function
$$W(\mathbf{x}) = \int_{\hOm} \frac{\operatorname{trace} G}{\det J} \mathrm{d} \bxi,$$
where $J = \partial_{\bxi} \bx$ denotes the map's Jacobian matrix while $G = J^T J$ is the map's metric tensor. The appearance of $\det J$ in the integrand's denominator prohibits the substitution of a folded $\bx: \hOm \rightarrow \Om^S$. As such, we introduce the \textit{regularised} Winslow function
\begin{align}
\label{eq:regularised_Winslow}
    W_{\varepsilon}(\bx) = \int_{\hOm} \frac{\operatorname{trace} G}{\mathcal{R}_{\varepsilon}(\det J)} \mathrm{d} \bxi,
\end{align}
where
$$\mathcal{R}_{\varepsilon}(x) := \frac{x + \sqrt{4 \varepsilon^2 + x^2}}{2}$$
is a regulariser that ensures $W_\varepsilon(\bx) < \infty$. The functional~\eqref{eq:regularised_Winslow} has a barrier property wherein the regularisation parameter $\varepsilon > 0$ tunes the degree of penalisation for maps $\bx: \hOm \rightarrow \Om^s$ with $\det J(\bx) < 0$. We may therefore optimise~\eqref{eq:regularised_Winslow} to eliminate folds by exploiting the integrand's barrier property:
\begin{align}
\label{eq:regularised_Winslow_opt}
    W_{\varepsilon}(\bx_h) \rightarrow \min \limits_{\bx_h \in \mathcal{U}^{\mathbf{f}}_h}.
\end{align}
In practice, we utilise $\varepsilon = 10^{-4}$. Despite regularisation, the strong barrier property gives rise to a small radius of convergence. As such, an initial iterate that is close to the minimiser of~\eqref{eq:regularised_Winslow_opt}, such as the solution of~\eqref{eq:C0_DG_discrete}, is mandatory. \\
As an alternative to~\eqref{eq:regularised_Winslow_opt} we may adopt the regularised weak form discretisation introduced in \cite{hinz2024use} which approximates a harmonic map by seeking the root
\begin{align}
\label{eq:inverse_elliptic_weak_pullback}
    \text{find } \bx_h \in \mathcal{U}^{\mathbf{f}}_h \quad \text{s.t.} \quad \mathcal{L}^{\text{W}}_{\varepsilon}(\bx_h, \boldsymbol{\phi}_h) = 0, \quad \forall \boldsymbol{\phi}_h \in \mathcal{U}^{\mathbf{0}}_h,
\end{align}
where
\begin{align}
\label{eq:inverse_elliptic_weak_pullback_operator}
    \mathcal{L}^{\text{W}}_{\varepsilon}(\bx, \boldsymbol{\phi}) := \int \limits_{\hOm}\frac{\partial_{\bxi} \boldsymbol{\phi} \, \colon \, A(\partial_{\bxi} \bx)}{\mathcal{R}_\varepsilon(\det J)} \mathrm{d} \bxi
\end{align}
and $A(\partial_\bxi \bx) := \operatorname{cof}(G(\partial_\bxi \bx))$ the metric tensor's cofactor matrix (c.f.~\eqref{eq:C0_DG_discrete_operator}). As in~\eqref{eq:regularised_Winslow}, we utilise $\varepsilon = 10^{-4}$. As before, the resulting root-finding problem has a strong barrier property, thus leading to a discrete root that is rarely folded (even when~\eqref{eq:C0_DG_discrete} results in a folded map). To find the discrete problem's root, we utilise Newton's method. As before, initialisation by the discrete root $\bx_h^\circ: \hOm \rightarrow \Om^S$ of~\eqref{eq:C0_DG_discrete} is mandatory because the operator's barrier property results in a small radius of convergence. \\
In \cite{hinz2024use}, it is shown that~\eqref{eq:inverse_elliptic_weak_pullback} has a faster rate of convergence than~\eqref{eq:C0_DG_discrete}, which generally leads to more satisfactory outcomes from a parameterisation quality perspective. As such, we may choose to solve~\eqref{eq:inverse_elliptic_weak_pullback} initialised by the discrete root $\bx_h^\circ: \hOm \rightarrow \Om^S$ of~\eqref{eq:C0_DG_discrete} even when $\bx_h^\circ: \hOm \rightarrow \Om^S$ is not folded. \\

\noindent In the rare event of $\bx_h^\circ: \hOm \rightarrow \Om^s$ being folded, we first untangle using~\eqref{eq:regularised_Winslow_opt} and initialise~\eqref{eq:inverse_elliptic_weak_pullback} by the untangled map. For untangling, a total of $\sim 25$ gradient-based iterations suffice in the vast majority of cases. The steps of solving the algorithmically favorable root-finding problem~\eqref{eq:C0_DG_discrete} and initialisation of~\eqref{eq:inverse_elliptic_weak_pullback} with the proxy solution $\bx_h^\circ$, potentially in conjunction with untangling based on~\eqref{eq:regularised_Winslow_opt}, forms a highly effective strategy that consistently achieves a near-perfect success rate in practical applications. \\

\begin{figure}[h!]
\centering
\begin{subfigure}[b]{0.45\textwidth}
    \includegraphics[width=0.95\linewidth]{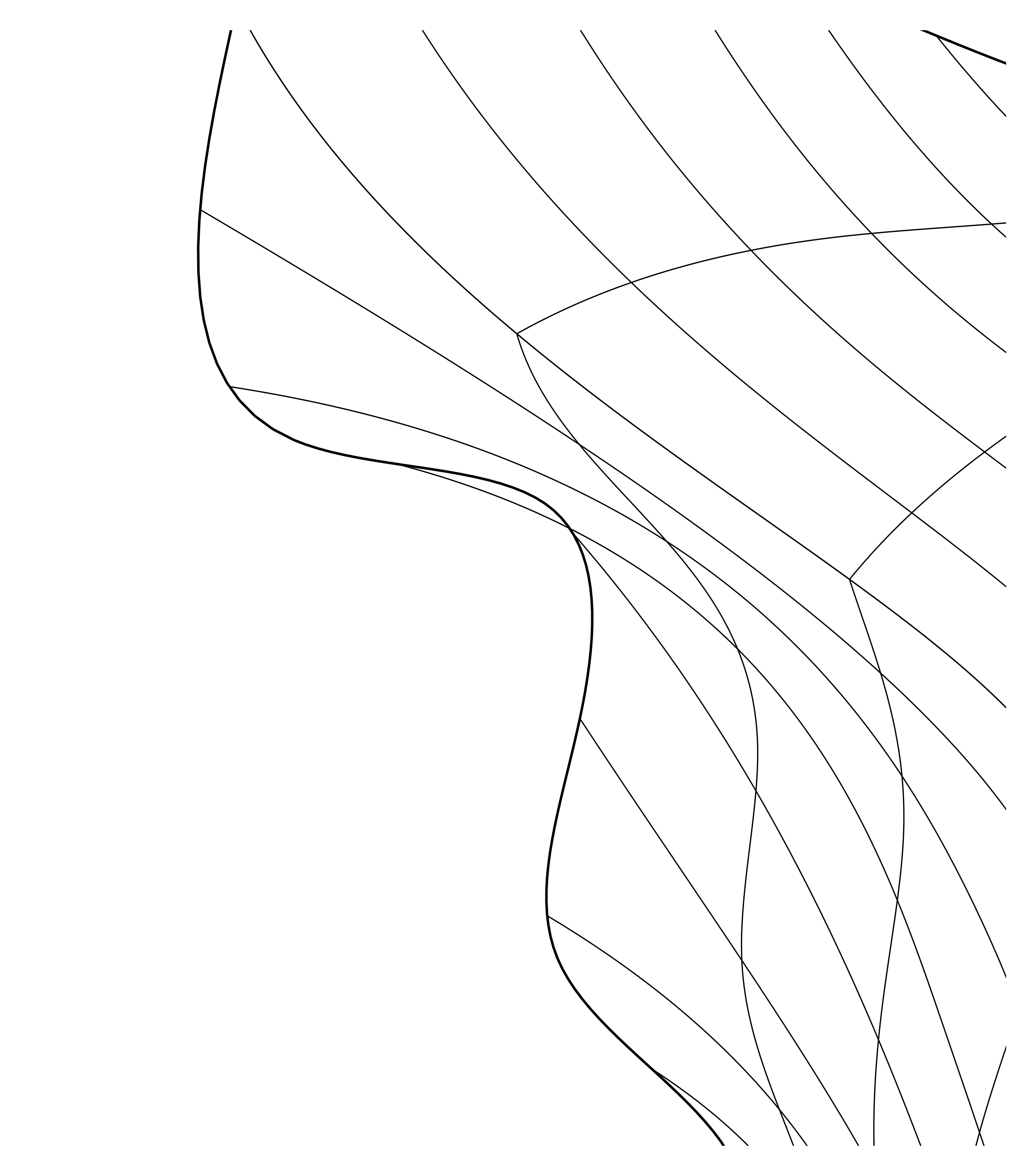}
    \caption{}
\end{subfigure}
\begin{subfigure}[b]{0.45\textwidth}
    \includegraphics[width=0.95\linewidth]{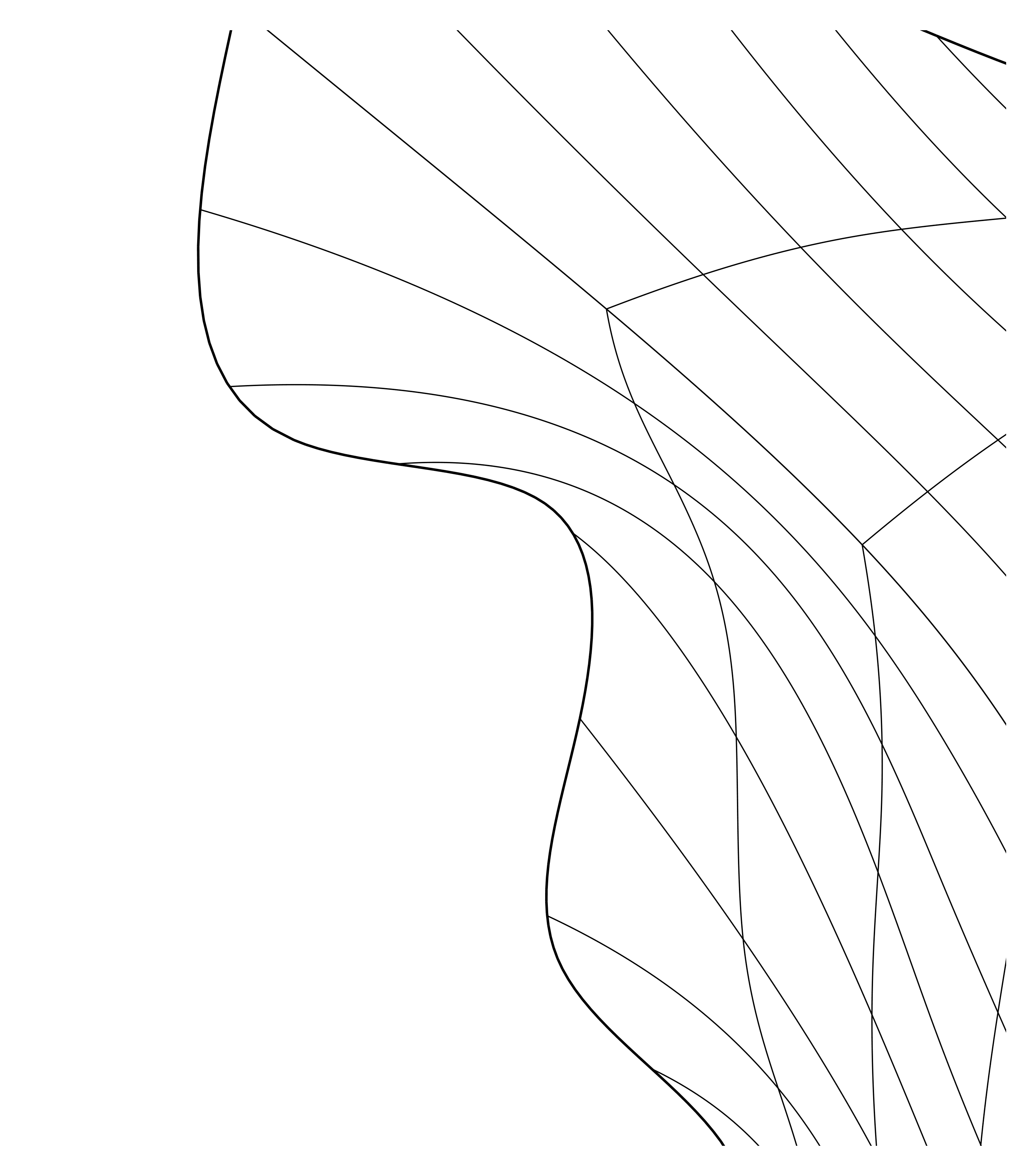}
    \caption{}
\end{subfigure}
\caption{Figure showing an example of removing a foldover using this section's untangling strategies.}
\label{fig:post_processing_untangle}
\end{figure}

\noindent Figure \ref{fig:post_processing_untangle} shows the result of removing a slight foldover near a significant inward concavity of $\partial \Om^S$ using this section's untangling strategies. This example is taken from the spline parameterisation of the face marked by the letter \textit{'A'} in Figure \ref{fig:results_africa} \textbf{(a)} which is one of the three faces that required a posteriori untangling in the parameterisation of the entire plane graph.

\begin{figure}[h!]
\centering
\includegraphics[align=c, width=0.22\textwidth]{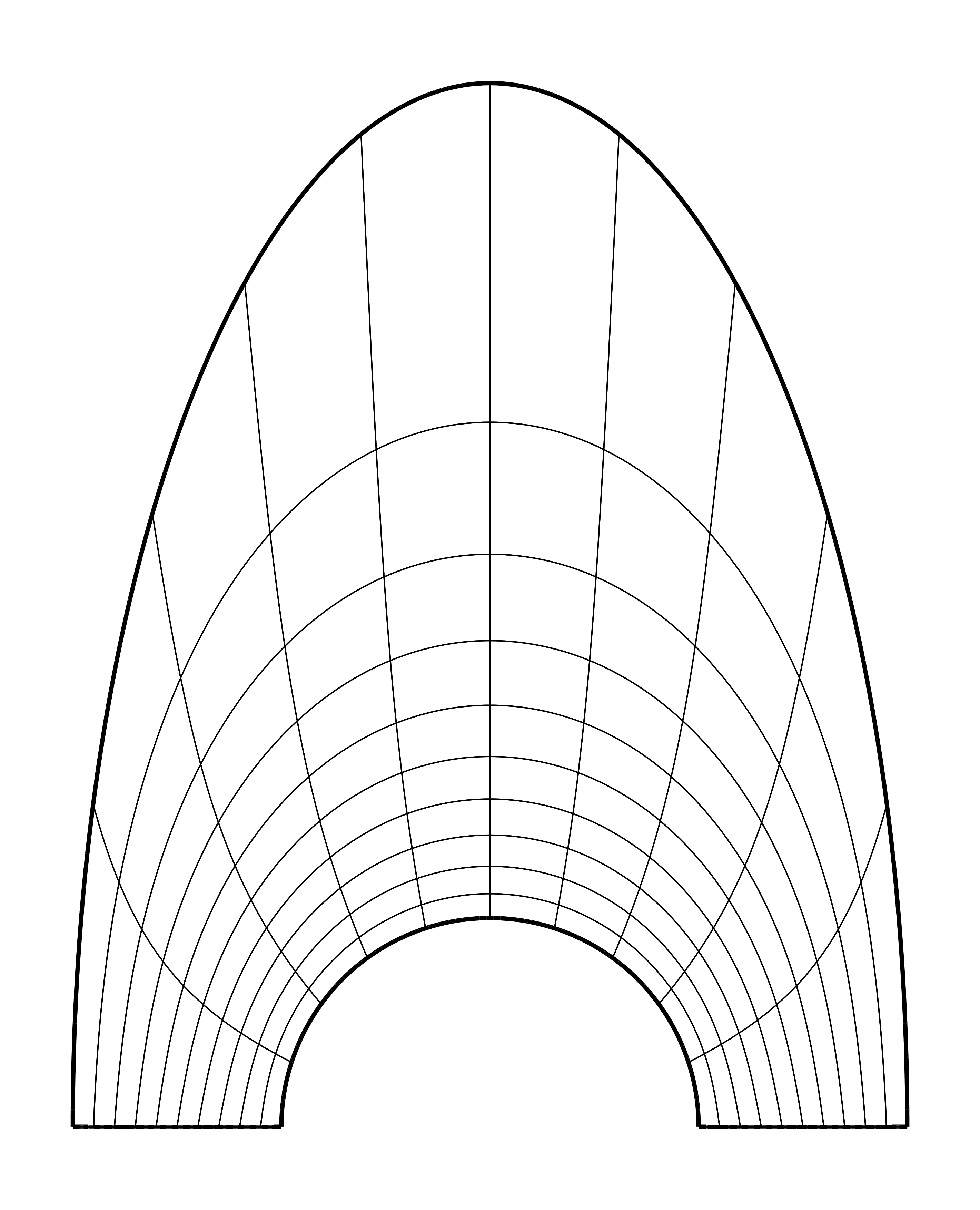}
\includegraphics[align=c, width=0.22\textwidth]{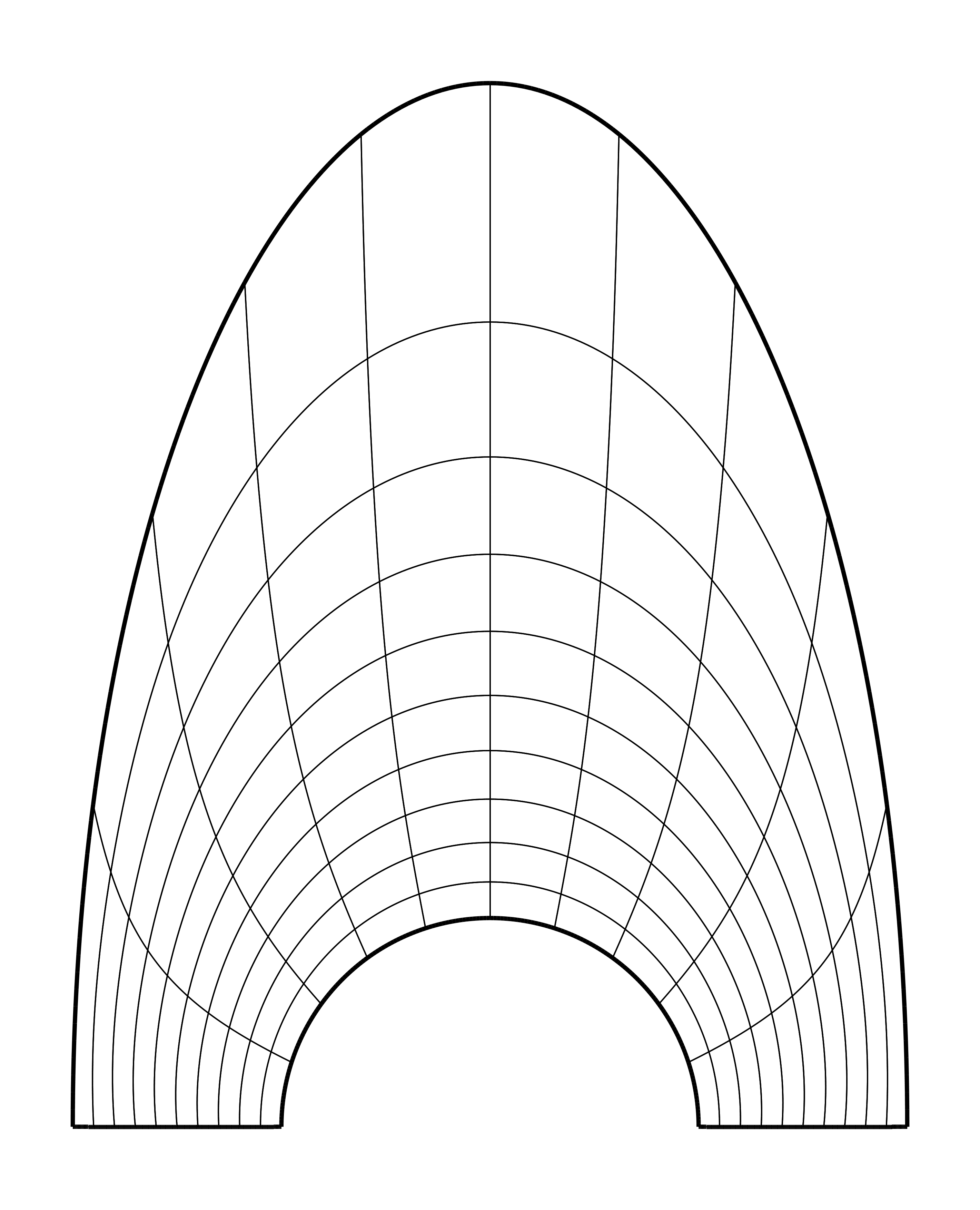}
\includegraphics[align=c, width=0.22\textwidth]{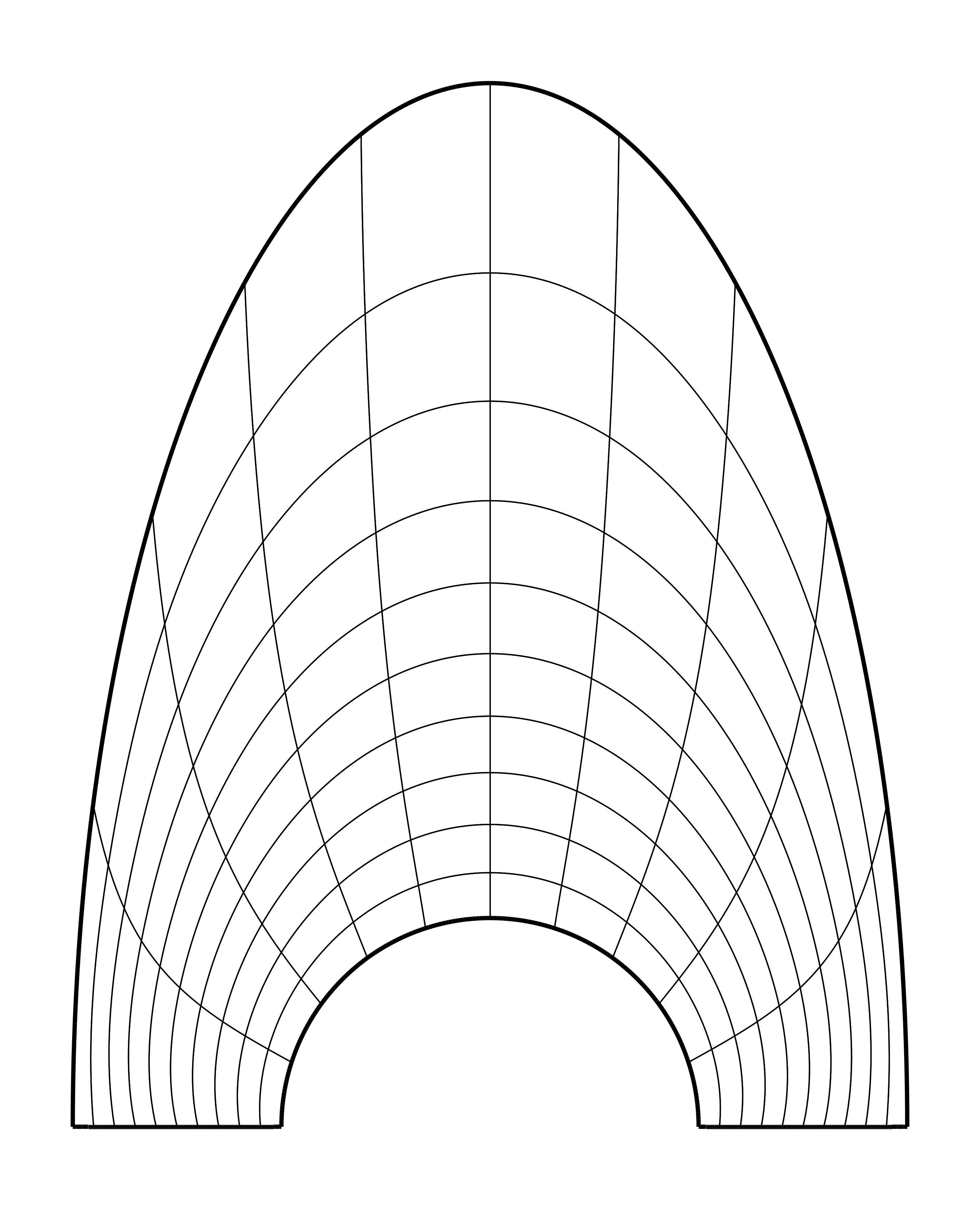}
\includegraphics[align=c, width=0.22\textwidth]{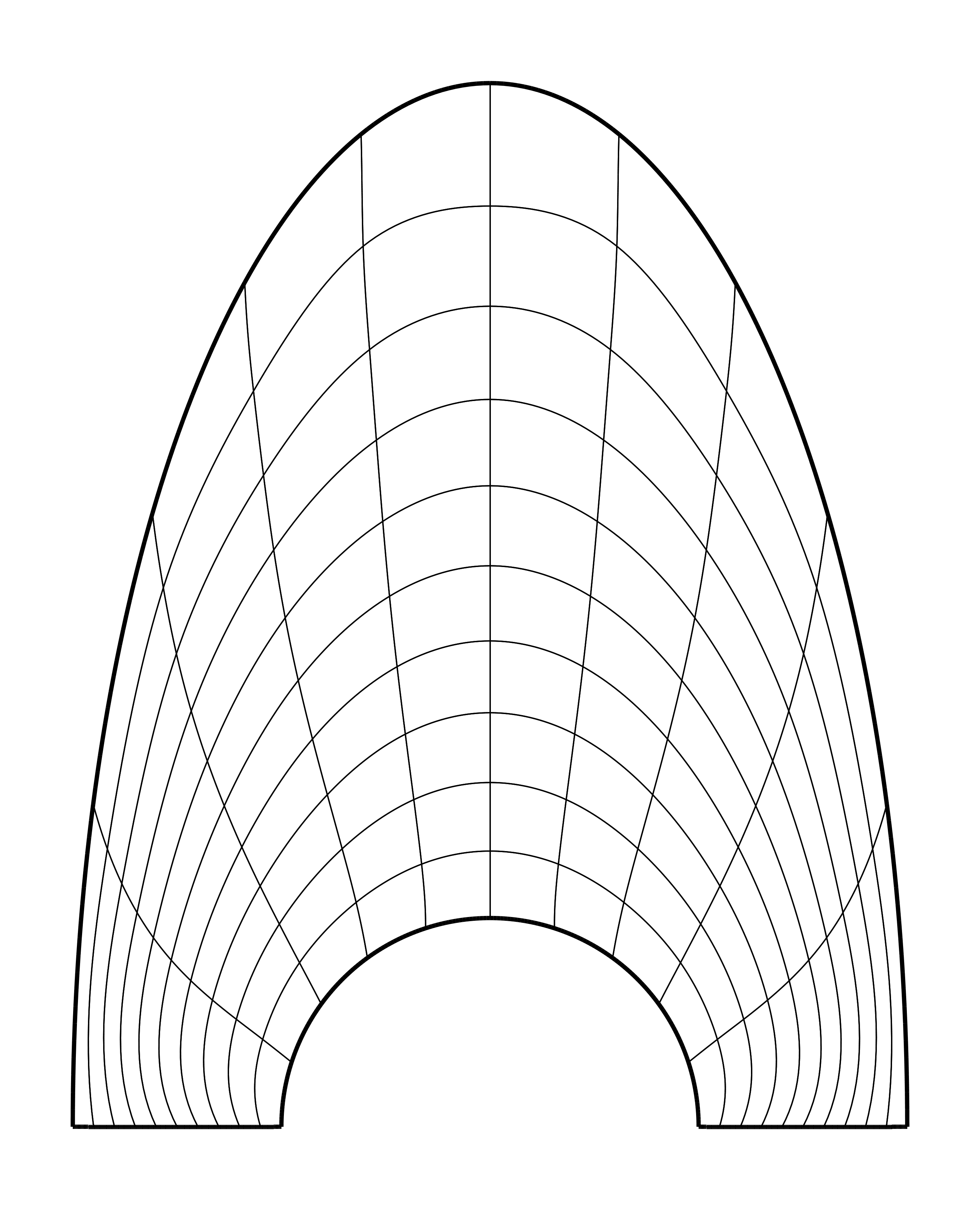}
\caption{The result of performing cell size homogenisation of a horse shoe domain using~\eqref{eq:inverse_elliptic_weak_D} in combination with~\eqref{eq:homogeneity_diffusivity} for $k \in \{0, 0.5, 1, 3\}$. For $k = 0$, the horse shoe domain's protrusion leads to significant parametric nonuniformity which decreases for increasing values of $k$.}
\label{fig:cell_size_hom_horse_shoe}
\end{figure}

\begin{figure}[h!]
\centering
\begin{subfigure}[b]{0.48\textwidth}
\begin{tikzpicture}
            \node[anchor=south west,inner sep=0] (image) at (0,0) {\includegraphics[align=c, width=0.95\linewidth]{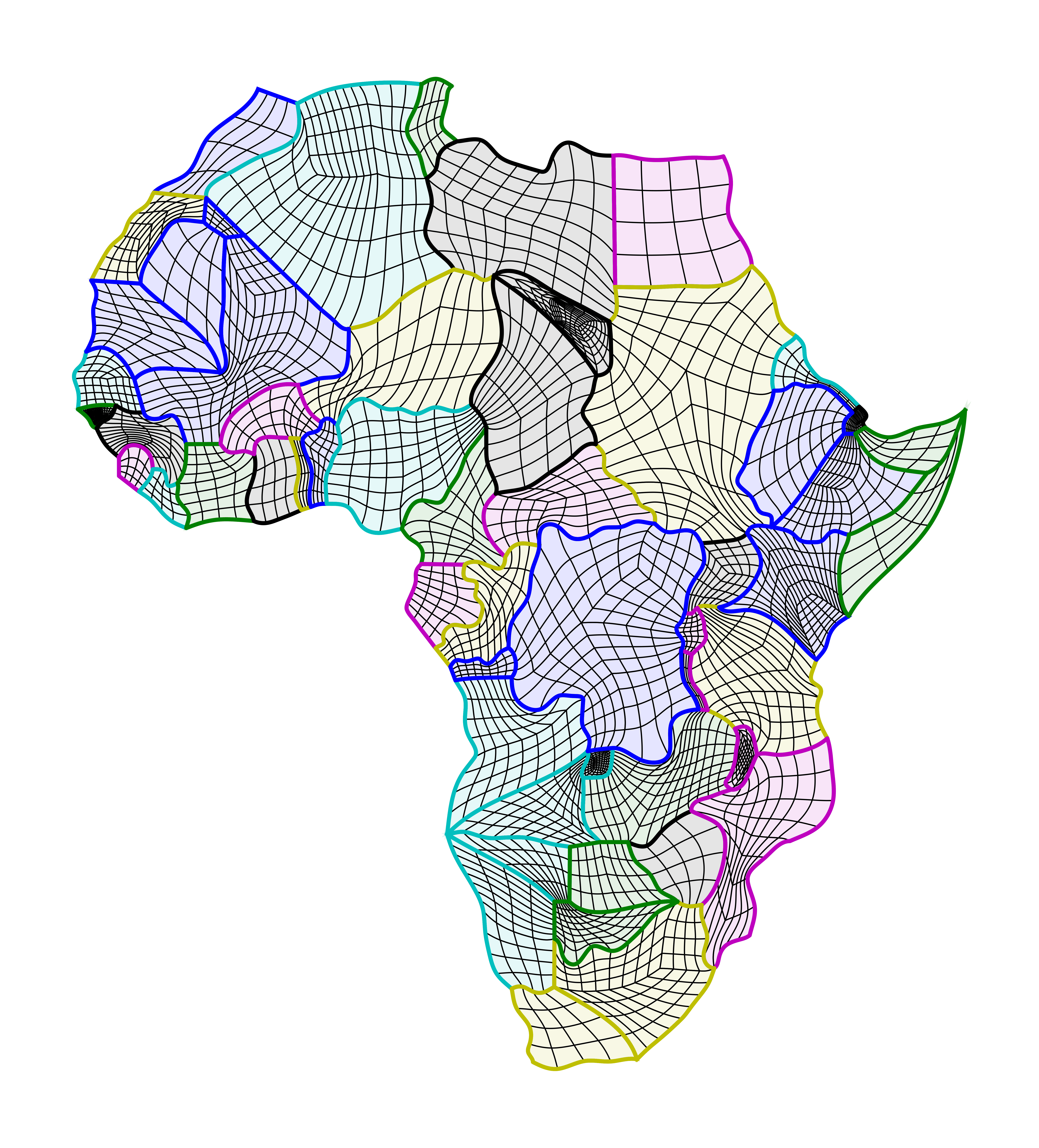}};
            \begin{scope}[x={(image.south east)},y={(image.north west)}]
                \node[anchor=center] at (0.56, 0.46) {\color{red}\textbf{A}}; 
            \end{scope}
        \end{tikzpicture}
        \caption{}
\end{subfigure}
\begin{subfigure}[b]{0.48\textwidth}
    \includegraphics[width=0.95\linewidth]{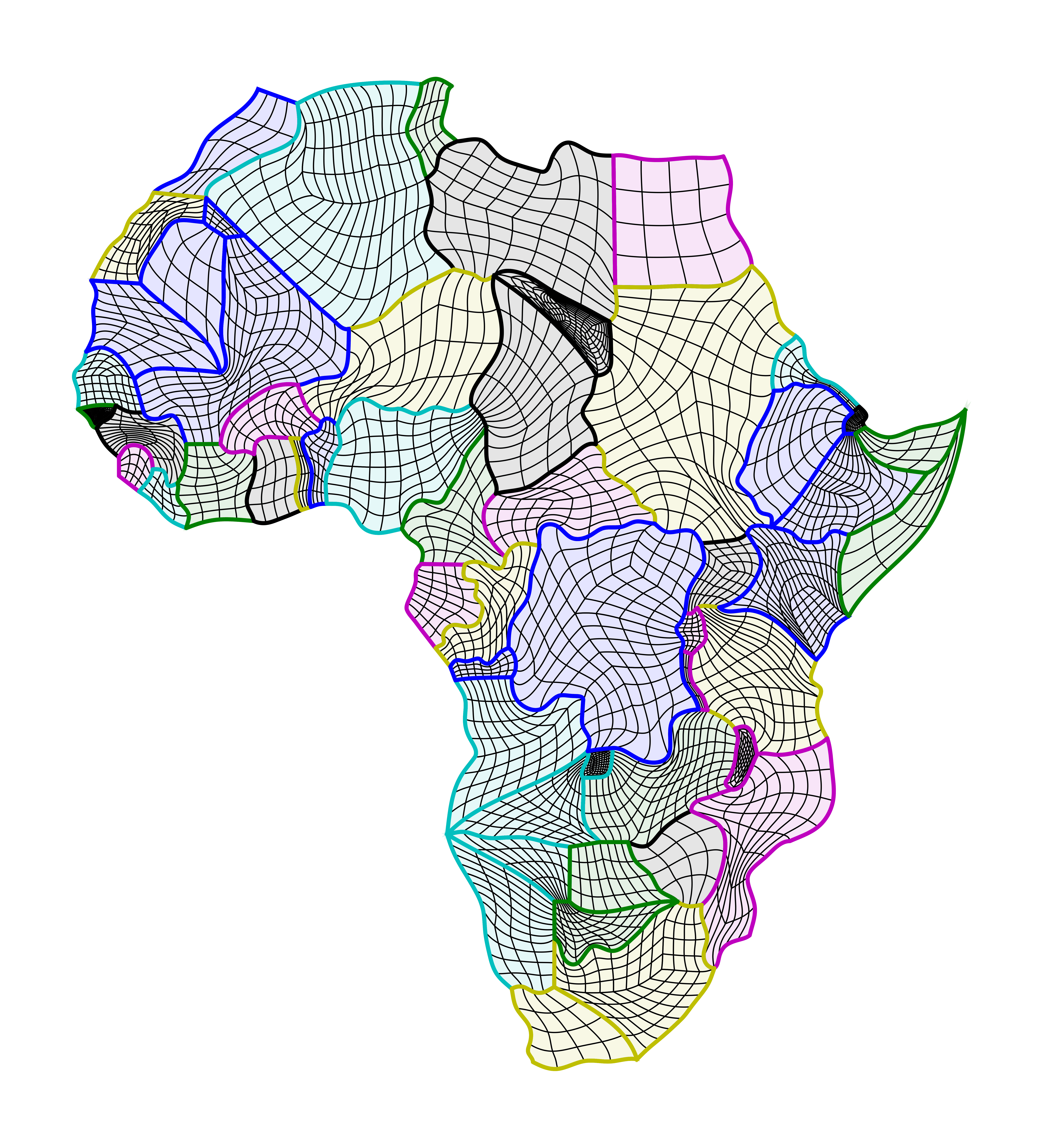}
    \caption{}
\end{subfigure}
\caption{Figure showing an example of performing cell size homogenisation on the plane graph from Figure \ref{fig:results_africa} \textbf{(b)} using $k = 1$ \textbf{(a)} and $k=4$ \textbf{(b)}.}
\label{fig:cell_size_hom_africa}
\end{figure}

Upon finalisation of the (if applicable) untangling stage, we may utilise the techniques presented in \cite{hinz2024use} to further optimise each face's $F \in \mathcal{F}$ spline parameterisation. We focus on the desired parametric features of:
\begin{enumerate}
    \item cell size homogeneity;
    \item patch interface removal.
\end{enumerate}
A well-documented pathology of harmonic maps is their tendency to exhibit significant parametric nonuniformity in the vicinity of protruded boundary segments (see Figure \ref{fig:cell_size_hom_horse_shoe} \textbf{(a)}). A remedy is the introduction of a nonconstant diffusivity $D \in \operatorname{SPD}^{2 \times 2}$ which requires the map to be the solution of the following inverted PDE problem
\begin{align}
   i \in \{1, 2\}: \quad \nabla_\bx \cdot (D \nabla_\bx \bxi_i) = 0, \quad \text{s.t. } \bx \vert_{\partial \hOm} = \partial \Om^S,
\end{align}
which reduces to the ordinary harmonic map problem for $D = I^{2 \times 2}$.  The associated semi-linear operator is a generalisation of~\eqref{eq:inverse_elliptic_weak_pullback_operator} and we utilise the same symbol for convenience:
\begin{align}
\label{eq:inverse_elliptic_weak_D}
    \mathcal{L}^{\text{W}}_{\varepsilon}(\bx, \boldsymbol{\phi}, D) = \int \limits_{\hOm}\frac{\partial_\bxi \boldsymbol{\phi} \colon (C^T D C)}{\mathcal{R}_\varepsilon(\det J)} \mathrm{d} \bxi, \quad \text{where} \quad C(\partial_{\bxi} \bx) := \operatorname{cof}(\partial_\bxi \bx)^T
\end{align}
is the Jacobian's transposed cofactor matrix. A generalisation to $\hOm^\br \neq \hOm$ is found in \cite{hinz2024use}. For point 1. (cell size homogeneity), we introduce the family of diffusivities $D^k \in \operatorname{SPD}^{2 \times 2}$ as proposed in \cite{hinz2024use} with
\begin{align}
\label{eq:homogeneity_diffusivity}
    D^{k} := \kappa \left( \det \partial_\bmu \bx \right)^{-k} I^{2 \times 2} \quad \text{on} \quad \hOm_i \subseteq \hOm, \quad \text{where} \quad \bmu = (\mu_1, \mu_2)^T
\end{align}
denotes the free coordinate functions in $\Om^\square = (0, 1)^2$ and $\bx = \bx(\bmu)$ is regarded as a function of $\bmu$ (rather than $\bxi$) on each individual $\hOm_i \subseteq \hOm$. The parameter $k \geq 0$ tunes the degree of homogenisation, with more drastic outcomes for larger values of $k$, while $\kappa = \kappa(k)$ is a normalisation constant. We note that~\eqref{eq:inverse_elliptic_weak_D} retains the barrier property, which consistently yields bijective maps, even for extreme diffusive anisotropy.
Figure \ref{fig:cell_size_hom_africa} shows the plane graph from Figure \ref{fig:results_africa} \textbf{(b)} after homogenisation of each face using~\eqref{eq:inverse_elliptic_weak_D} and~\eqref{eq:homogeneity_diffusivity} for $k \in \{1, 4\}$. \\

\begin{figure}[h!]
\centering
\begin{subfigure}[b]{0.48\textwidth}
    \centering
    \includegraphics[width=0.95\linewidth,  angle=200]{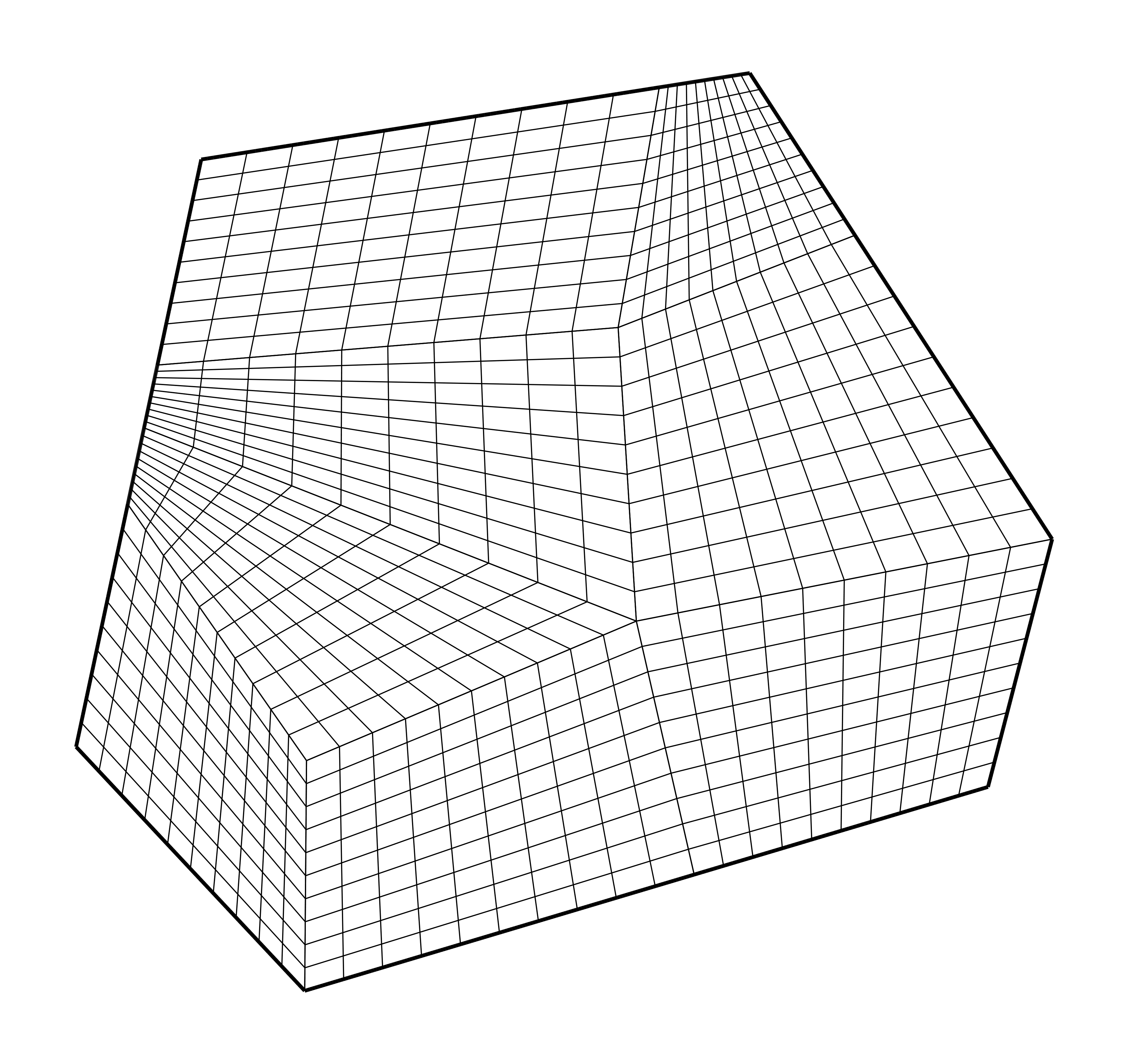}
    \caption{}
\end{subfigure}
\begin{subfigure}[b]{0.48\textwidth}
    \includegraphics[width=0.95\linewidth]{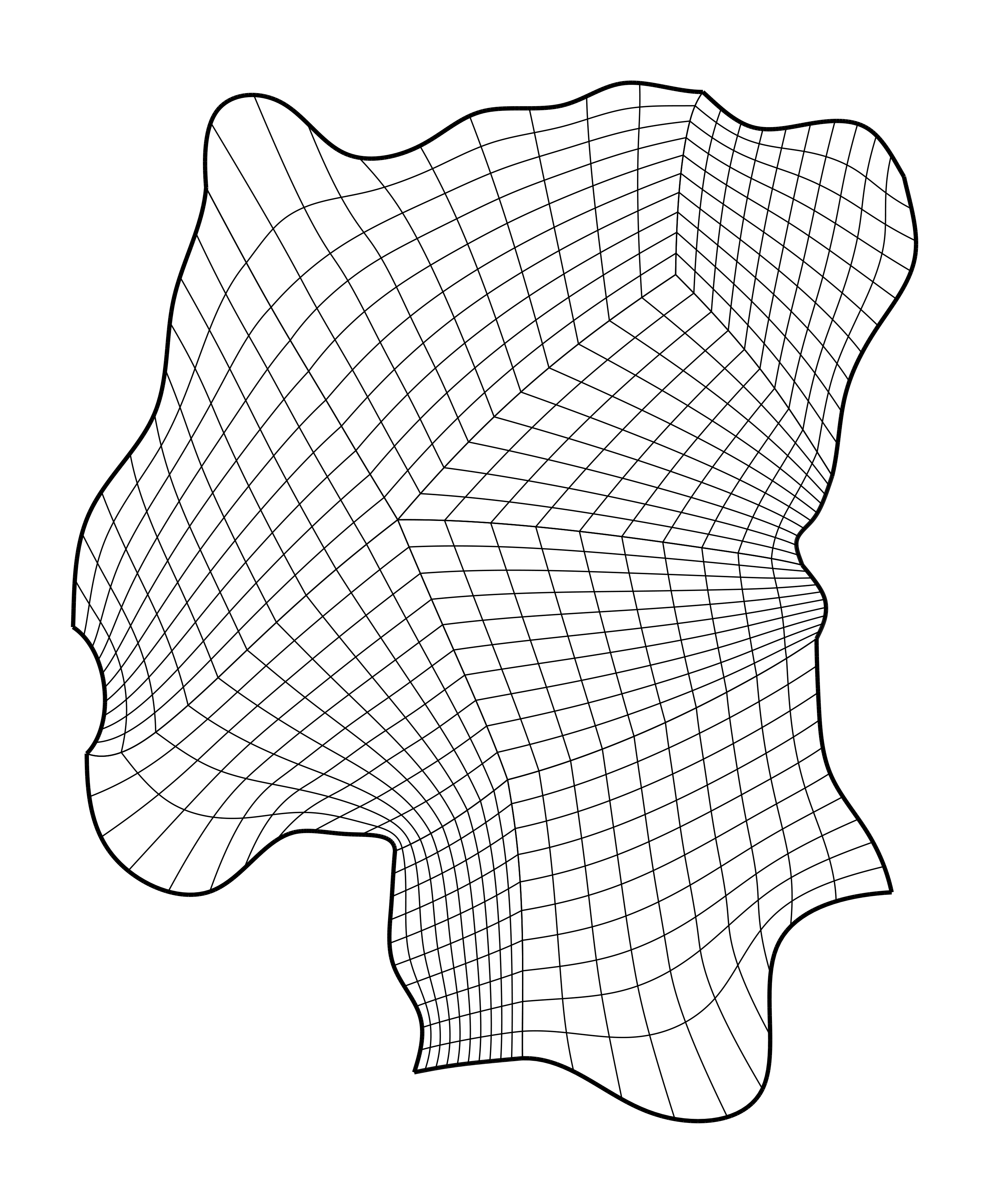}
    \caption{}
\end{subfigure} \\
\begin{subfigure}[b]{0.48\textwidth}
    \centering
    \includegraphics[width=0.95\linewidth, angle=200]{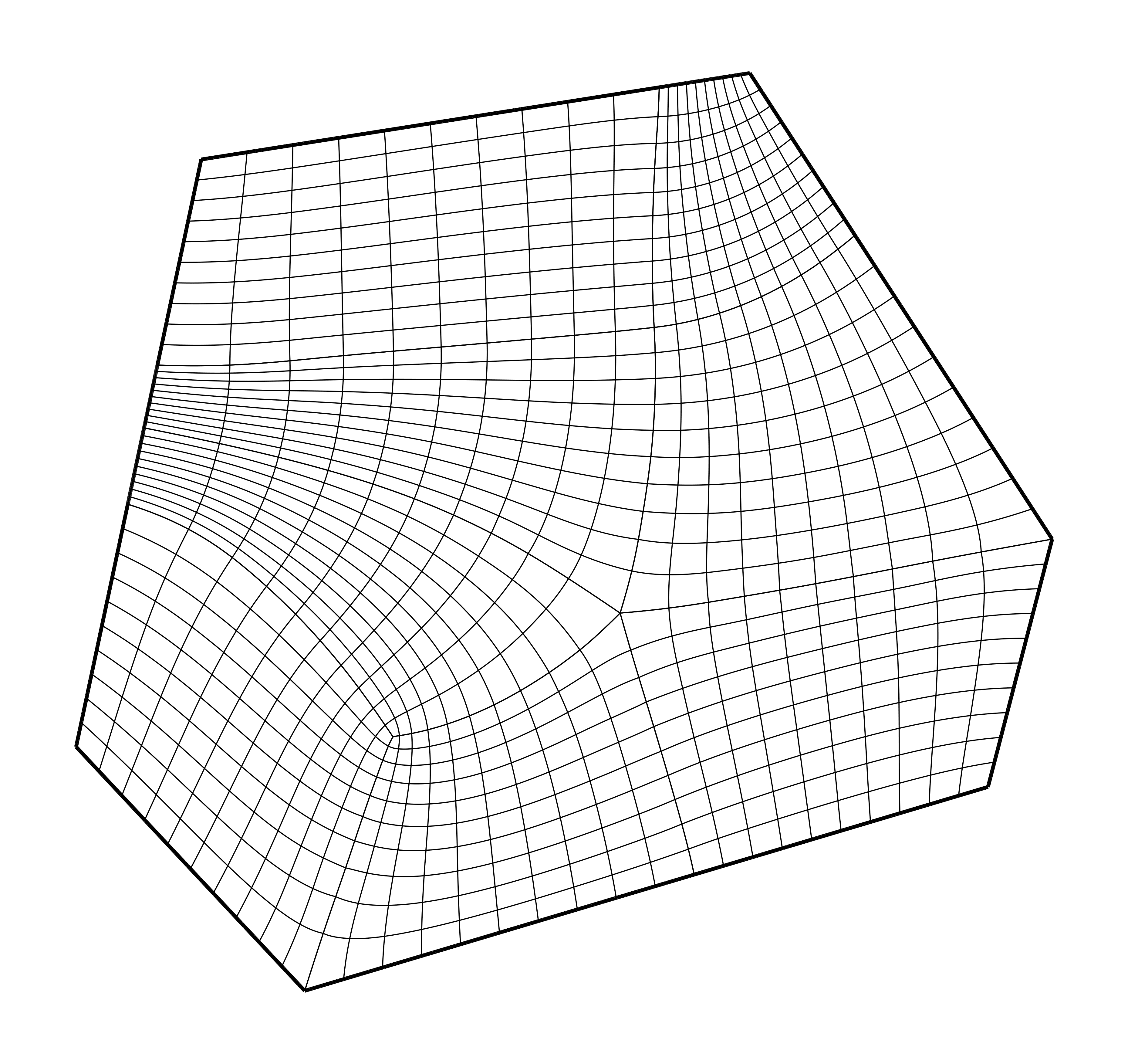}
    \caption{}
\end{subfigure}
\begin{subfigure}[b]{0.48\textwidth}
    \includegraphics[width=0.95\linewidth]{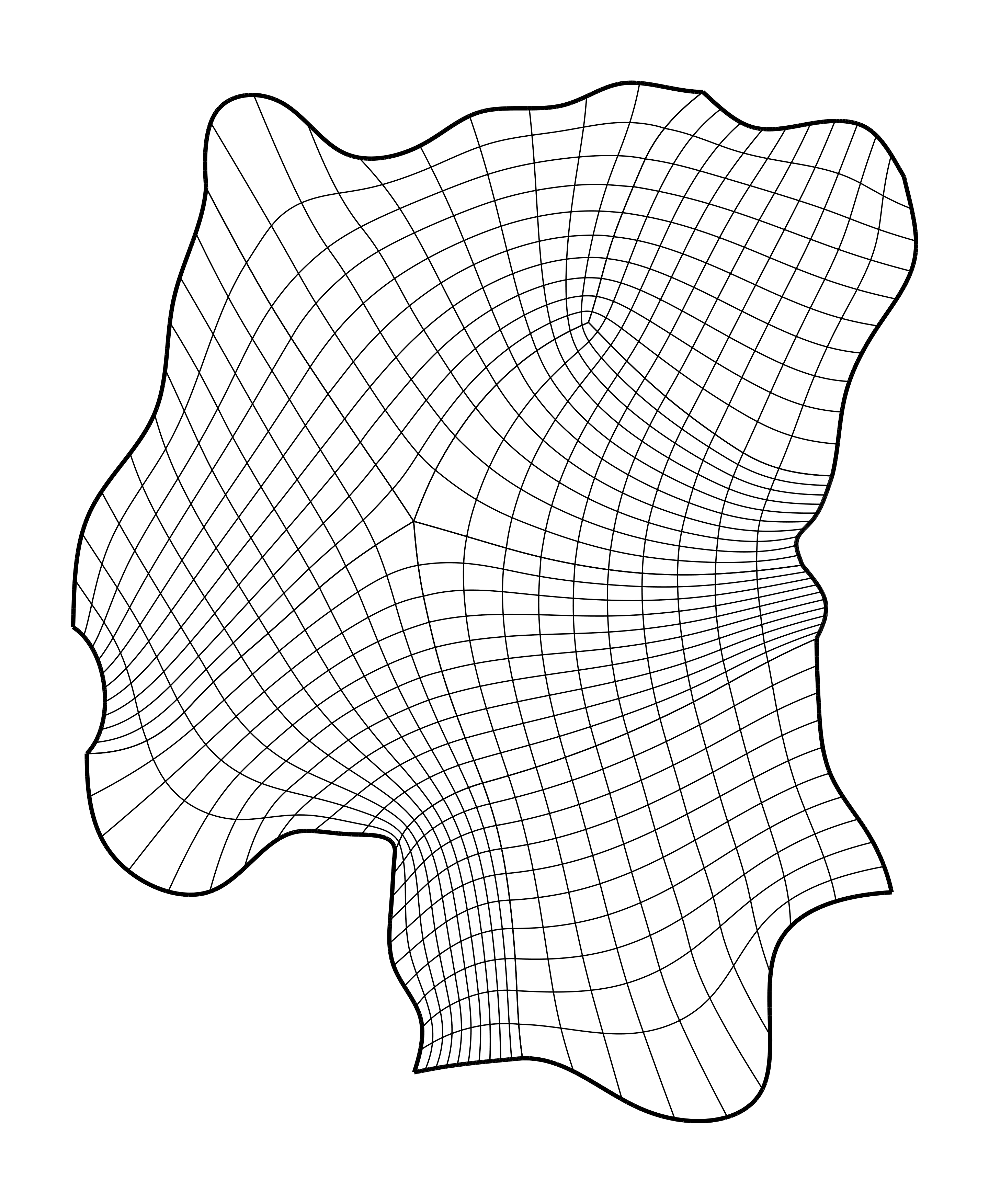}
    \caption{}
\end{subfigure}
\caption{Figure showing the reference parameterisation \textbf{(b)} of the face marked by the letter \textit{'A'} in Figure \ref{fig:cell_size_hom_africa} \textbf{(a)} as well as a parameterisation that removes the steep interface angles \textbf{(d)}. The corresponding control domain parameterisations are depicted in \textbf{(a)} and \textbf{(c)}, respectively. With respect to the reference parameterisation, the map from \textbf{(d)} shows notably reduced interface angles at the expense of a slight reduction in parametric uniformity.}
\label{fig:interface_removal_face}
\end{figure}

\noindent For point 2. (patch interface removal), we note that since an (inversely) harmonic $\bx: \hOm^\br \rightarrow \Om^S$ is a diffeomorphism in $\hOm^\br$, a covering of $\hOm^\br$ that removes the steep angles formed by the patch interfaces will, by extension, remove them in the harmonic map $\bx^\br: \hOm^\br \rightarrow \Om^S$ and in the pullback map $\bx: \hOm \rightarrow \Om^S$. Given the reference controlmap $\br: \hOm \rightarrow \hOm^\br$, we may compute a new controlmap $\bs: \hOm^\br \rightarrow \hOm^\br$ as the solution of the following PDE problem
\begin{align}
\label{eq:patch_interface_removal}
    i \in \{1, 2\}: \quad \nabla \cdot \left(D^\perp \nabla \bs_i \right) = 0, \quad \text{s.t.} \quad \bs = \br \text{ on } \partial \hOm^\br,
\end{align}
where
\begin{align}
    D^\perp(\br) = \widehat{\partial}_{\mu_1} \br \otimes \widehat{\partial}_{\mu_1} \br + \widehat{\partial}_{\mu_2} \br \otimes \widehat{\partial}_{\mu_2} \br \quad \text{on } \br(\hOm_i), \quad \text{with } \widehat{\partial}_{\mu_i} \br := \frac{\partial_{\mu_i} \br}{\| \partial_{\mu_i} \br \|},
\end{align}
the normalised counterpart of the $i$-th component of the tangent bundle $\partial_{\bmu} \br$ in the local $\bmu$ coordinate system. We note that~\eqref{eq:patch_interface_removal} need not be formulated as an inverted PDE since $\bs: \hOm^\br \rightarrow \hOm^\br$ maps into a convex domain. The PDE~\eqref{eq:patch_interface_removal} is discretised over $\mathcal{V}_h$ and solved in the usual way. Depending on the template layout $T \in \mathcal{T}$ associated with the covering of $\hOm$, the exact solution of~\eqref{eq:patch_interface_removal} may contain singularities $\det \partial_\br \bs \rightarrow 0$ or $\det \partial_\br \bs \rightarrow \infty$. While such singularities are avoided by discrete approximations, unbounded growth / shrinkage of $\det \partial_\br \bs$ will be observable in a refinement study, which may require adding a stabilisation term to $D^\perp$. For details, we refer to \cite{hinz2024use}. The discrete solution $\bs(\br): \hOm^\br \rightarrow \hOm^\br$ of~\eqref{eq:patch_interface_removal} may be reinterpreted as a map $\bs(\bxi): \hOm \rightarrow \hOm^\br$ in the original coordinate system via a pullback. The map $\bs(\bxi)$ then takes the role of the original controlmap $\br(\bxi): \hOm \rightarrow \hOm^\br$ and may be utilised in conjunction with~\eqref{eq:inverse_elliptic_weak_pullback} to compute a novel map $\bx^\bs(\bxi): \hOm \rightarrow \Om^S$ that (largely) removes the steep interface angles. \\
Figure \ref{fig:interface_removal_face} shows the parameterisation of the face in Figure \ref{fig:cell_size_hom_africa} \textbf{(a)} marked by the letter \textit{'A'} with and without interface removal, along with the corresponding parameterisations of the control domain $\hOm^\br$. \\

\begin{figure}[h!]
\centering
\begin{subfigure}[b]{0.48\textwidth}
    \includegraphics[width=0.95\linewidth]{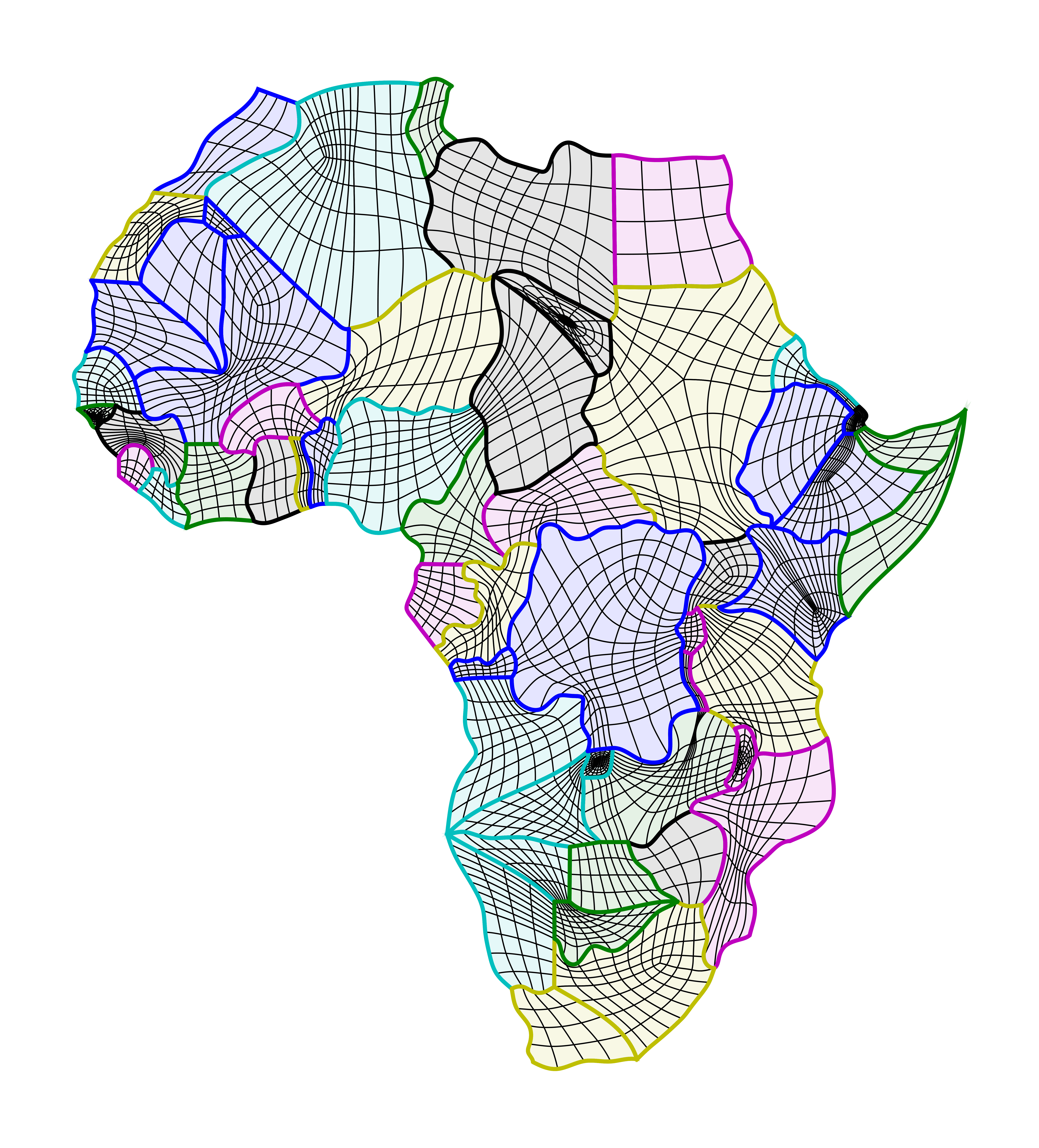}
    \caption{}
\end{subfigure}
\begin{subfigure}[b]{0.48\textwidth}
    \includegraphics[width=0.95\linewidth]{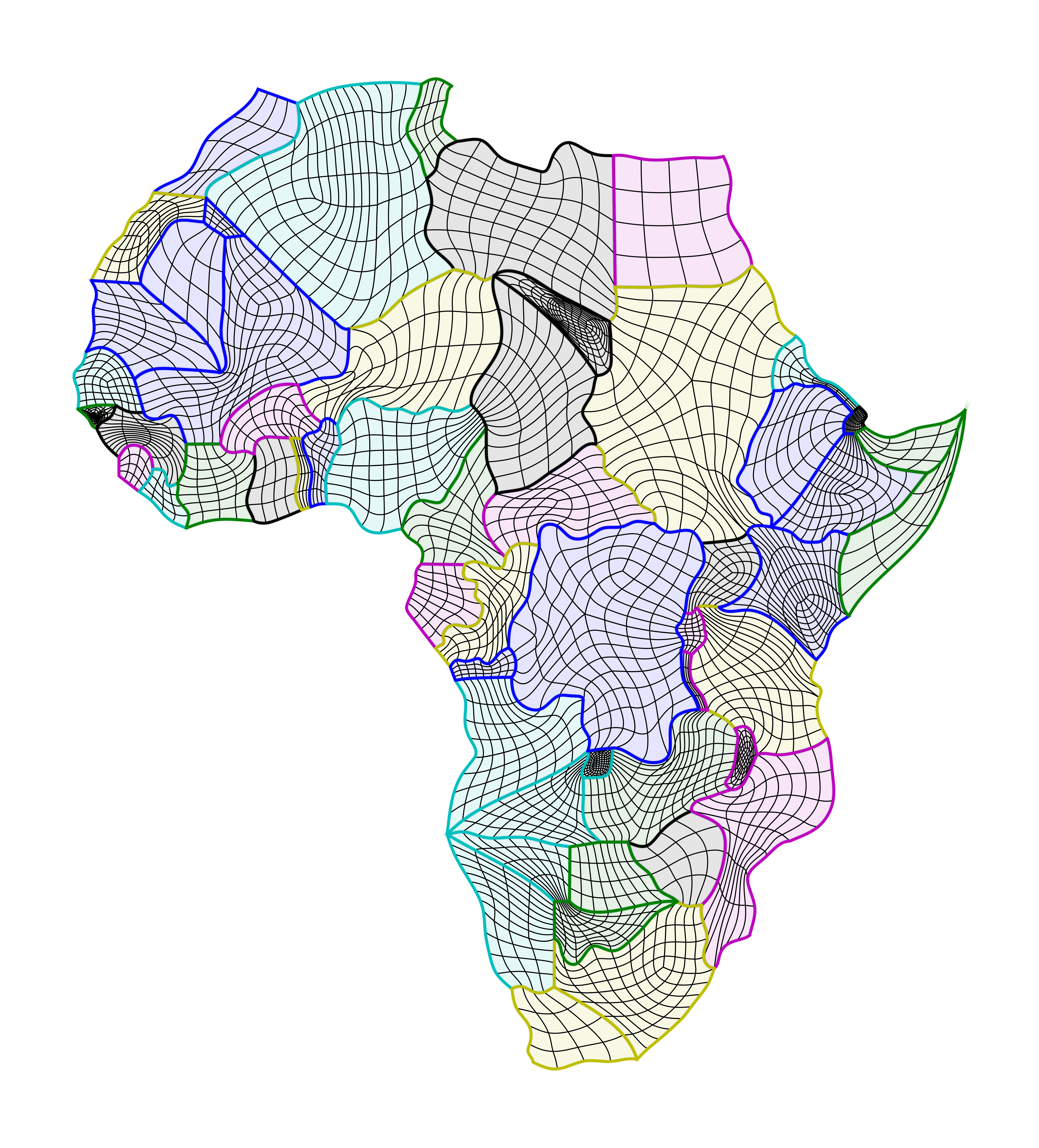}
    \caption{}
\end{subfigure}
\caption{Figure showing an example of performing interface removal in the absence of cell size homogenisation \textbf{(a)} and with homogenisation using $k=3$ \textbf{(b)}.}
\label{fig:interface_removal_africa}
\end{figure}

\noindent Figure \ref{fig:interface_removal_africa} \textbf{(a)} shows the result of performing interface removal on the plane graph from Figure \ref{fig:results_africa}.

\noindent Figure \ref{fig:interface_removal_africa} \textbf{(b)} combines interface removal with cell size homogenisation based on~\eqref{eq:inverse_elliptic_weak_D} and~\eqref{eq:homogeneity_diffusivity} using $k=3$. The result is a spline parameterised plane graph wherein each face's internal interface angles are largely absent, while Figure \ref{fig:interface_removal_africa} \textbf{(b)} shows notably improved cell size homogeneity over Figure \ref{fig:interface_removal_africa} \textbf{(a)}. \\
While the proposed techniques allow for parametric control over individual faces, achieving control across multiple faces simultaneously remains a major challenge. This limitation arises because elliptic PDEs do not allow for the simultaneous imposition of Dirichlet and Neumann boundary conditions. Addressing this challenge presents an topic for future research.
\section{Conclusion}
We have presented a framework for the spline-based parameterisation of plane graphs based on the concept of harmonic maps. The framework performs several preliminary operations on the input graph and then selects a suitable quad layout for each face from a catalogue of $\sim 6 \times 10^4$ templates for various numbers of boundary edges. For the template selection, we presented a total of three differing strategies. Hereby, computationally more demanding outperformed the computationally less demanding strategies in the essayed examples both from a robustness and from a parameterisation quality standpoint. \\
The framework was able to autonomously parameterise input plane graphs comprised of up to $42$ faces ($55$ after the removal of concave corners) and we presented several post-processing  techniques for the removal of fold-overs and for further tuning each face's parametric properties. \\
As it stands, the framework may still require manual intervention in the removal of concave corners since the proposed technique based on Hermite interpolation provides no guarantee that the resulting curve is strictly contained in the face's interior. For this, it may be desirable to find a more robust solution, for instance based on the techniques from \cite{nguyen2016isogeometric, haberleitner2017isogeometric} or by the use of a triangular approximation of the face. \\
Moreover, while the post-processing methods discussed in Section \ref{sect:post_processing} enable precise parametric control over the individual faces of the plane graph, expanding this capability to encompass parametric control across several faces simultaneously may constitute a topic for future research.

\section*{Declarations}
The authors gratefully acknowledge the support of the Swiss National Science Foundation through the project ‘‘Design-through-Analysis (of PDEs): the litmus test’’ n. 40B2-0 187094 (BRIDGE Discovery 2019). \\
The authors have no relevant financial or non-financial interests to disclose.
The authors have no competing interests to declare that are relevant to the content of this article.
All authors certify that they have no affiliations with or involvement in any organization or entity with any financial interest or non-financial interest in the subject matter or materials discussed in this manuscript.
The authors have no financial or proprietary interests in any material discussed in this article.
\appendix
\section{Floater's Algorithm}
\label{sect:appendix_floater}
\noindent Floater's algorithm (see Section \ref{sect:harmonic_maps}) is a computationally inexpensive method to acquire a decent approximation of the harmonic map $\mathbf{x}^{-1}: \Omega^S \rightarrow \hOm^{\mathbf{r}}$ over a polygonal domain $\Omega^S_h \approx \Omega^S$. The approach samples ordered point sets $X = \left(x_1, \ldots, x_n \right)$ and $\hat{X} = \left(\hat{x}_1^r, \ldots, \hat{x}_i^r \right)$ from $\partial \Omega^S$ and $\partial \hOm^{\mathbf{r}}$, respectively, leading to a canonical discrete correspondence $(\mathbf{f}_h^{\mathbf{r} \rightarrow \mathbf{x}})^{-1}: X \rightarrow \hat{X}$, with $(\mathbf{f}_h^{\mathbf{r} \rightarrow \mathbf{x}})^{-1}(x_i) = \hat{x}_i^r$. The boundary of the polygon $\Omega_h^S$ is forwarded to a meshing routine which returns the triangulation $\mathcal{T}_h$ with boundary vertices $X \subset \mathbb{R}^2$ and a total number of $N > n$ vertices. An approximately harmonic map is then reconstructed from the canonical piecewise-linear $C^0$-continuous Lagrange basis $P_h$ over $\mathcal{T}_h$, which leads to a linear problem of the form $A c = B f^r$. Here, the solution vector $c \in \mathbb{R}^{2(N - n)}$ contains the approximation's weights (w.r.t. to $P_h \times P_h$) that are not fixed from the boundary correspondence, while the vector $f^r \in \mathbb{R}^{2 n}$ is the concatination of the $\hat{x}_i^r \in \hat{X}$. For a computationally inexpensive method to assemble the SPD $A \in \mathbb{R}^{2N \times 2N}$ and the sparse $B \in \mathbb{R}^{2N \times 2n}$, we refer to \cite{floater2005surface}. The approximation $\left( \mathbf{x}_h^{\mathbf{r}} \right)^{-1} \approx \left( \mathbf{x}^{\mathbf{r}} \right)^{-1}$ maps the triangulation $\mathcal{T}_h$ onto a triangulation $\hat{\mathcal{T}}_h^r$ of $\hOm^{\mathbf{r}}_h \approx \hOm^\br$ which automatically yields a piecewise linear approximation $\mathbf{x}_h^\br: \hOm^\br_h \rightarrow \Om_h^S$ of an inversely harmonic map from the correspondence of the vertices $\hat{\mathcal{T}}^r_h \ni \hat{v}_i \rightarrow v_i \in \mathcal{T}_h$.

\section{Template Creation using Patch Adjacency Graphs}
\label{sect:appendix_PAG}
Approaches 1. and 2. from Section \ref{sect:harmonic_maps} operate on a multipatch spline space $\mathcal{V}_h$ defined over a quadrangulation of $\hOm_i$ which is itself represented by a template $T_i \in \mathbb{T}$. \\
To construct the catalogue of templates $\mathbb{T}$, this paper adopts the concept of patch adjacency graphs (PAGs) introduced in \cite{buchegger2017planar}. A PAG $P = (V_p, E_p)$ is an undirected graph that represents the patch connectivity of the quadrangulation of an $N$-sided ($N$ even) domain. The PAG places $N$ boundary nodes $v_q^b$ uniformly distributed on the unit circle, each representing a boundary vertex $v_\alpha \in \mathbb{V}(\partial E_i)$ of the $N$-sided template $T_i = (V_i, E_i, \mathcal{Q}_i)$. Then, it places an additional $N_p$ vertices $v_q$ in the interior of the resulting polygon, each representing one of the $N_p$ patches of the layout. Then, edges are drawn between the vertices, each representing the connection of (one side of) an interior patch to a boundary edge of $T_i$ or an adjacency between two patches. Here, an edge between $v_\alpha$ and $v_q^b$ means that one side of the patch represented by $v_{\alpha}$ is given by the edge between $(v_q, v_{q+1}) \subset \mathbb{V}(T_i)$. Drawing edges between adjacent boundary vertices $v_{q}^b, v_{q+1}^b$, this implies that $\operatorname{val}(v_q^b) = 3$ while $\operatorname{val}(v_q) = 4$ for interior nodes. \\

\begin{figure}[h!]
\centering
    \begin{subfigure}[b]{0.98\textwidth}
        \centering
        \includegraphics[align=c, width=.2\textwidth]{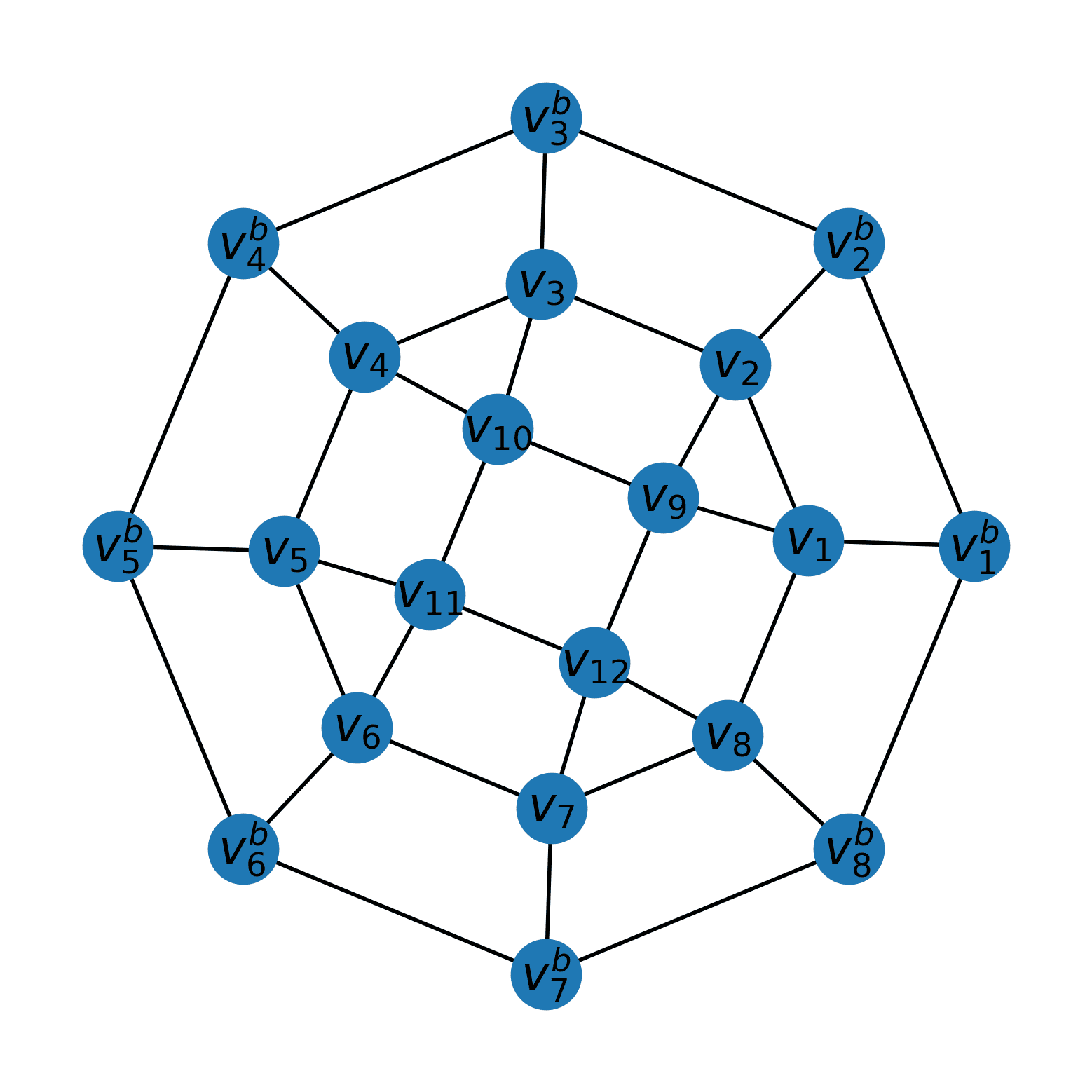}
        \includegraphics[align=c, width=.2\textwidth]{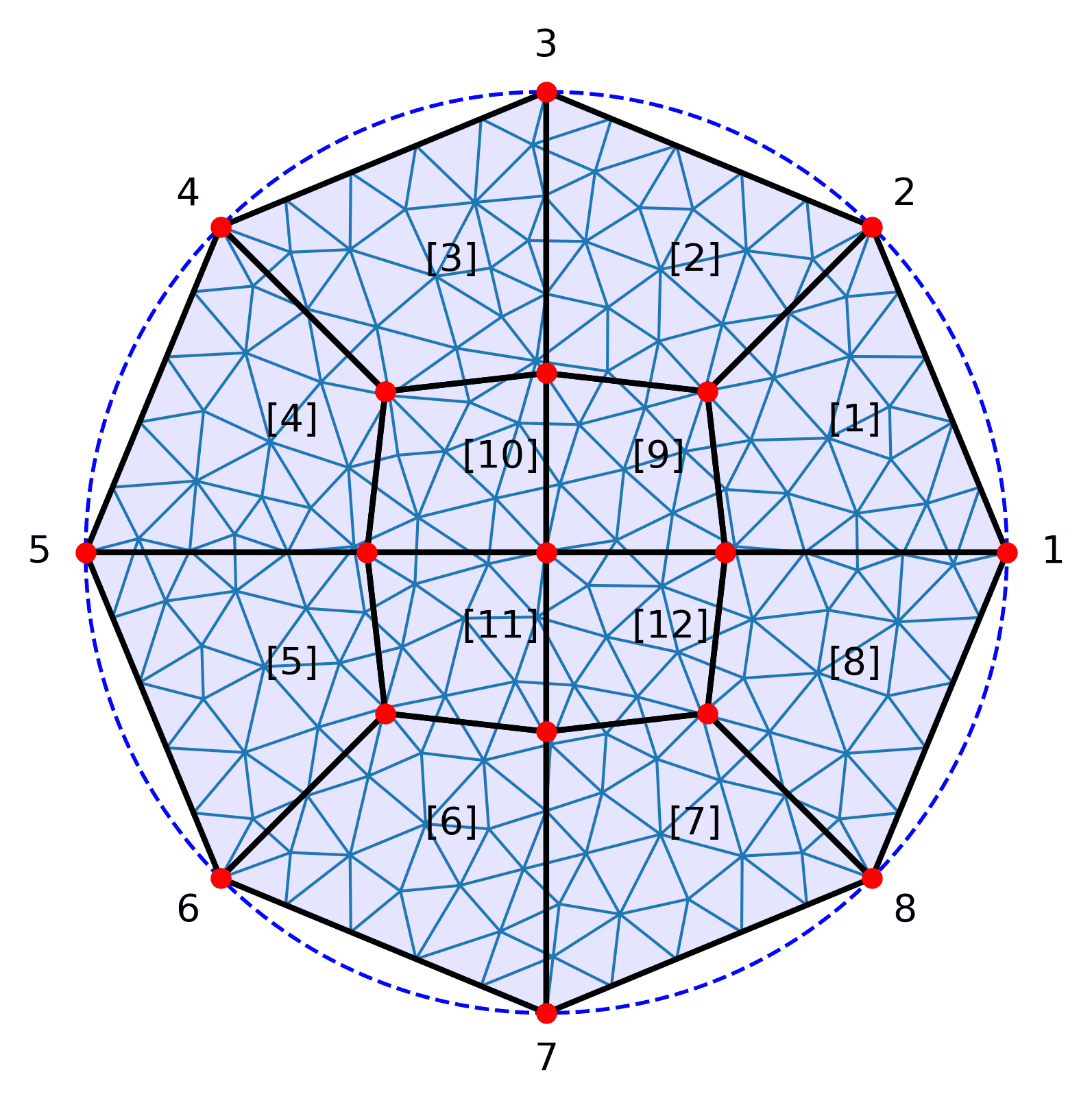}
        \includegraphics[align=c, width=.3\textwidth]{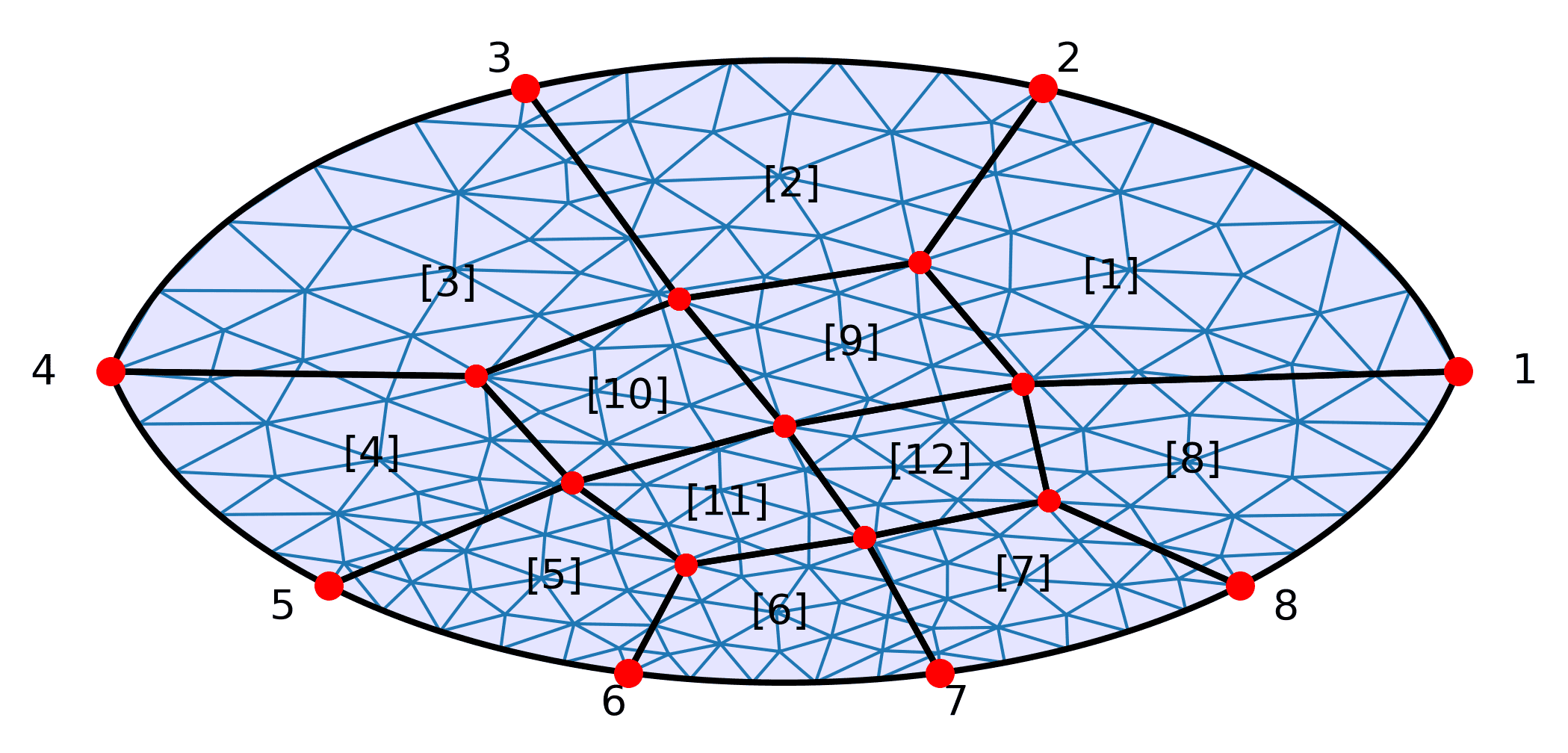} \\
        \includegraphics[align=c, width=.3\textwidth]{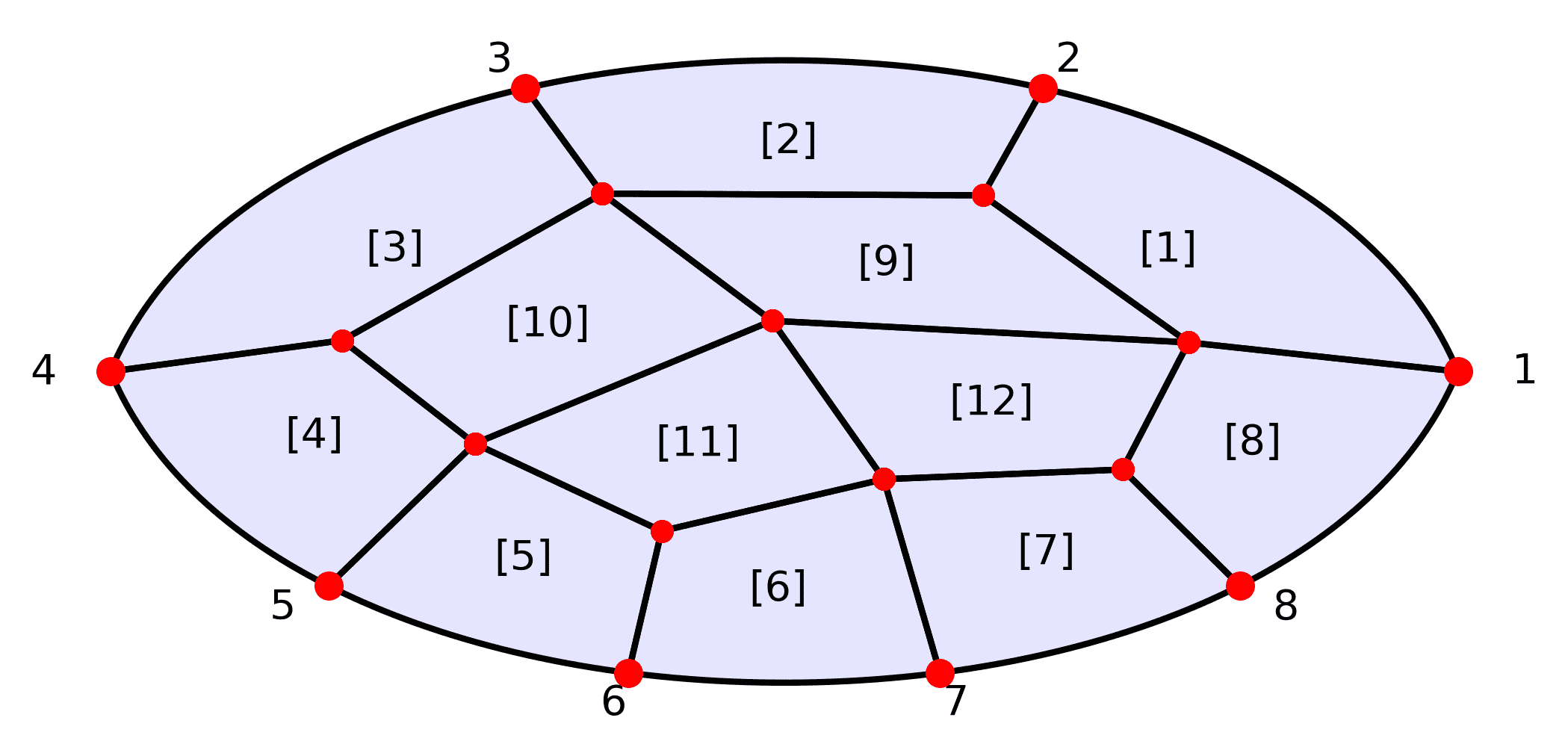}
        \includegraphics[align=c, width=.3\textwidth]{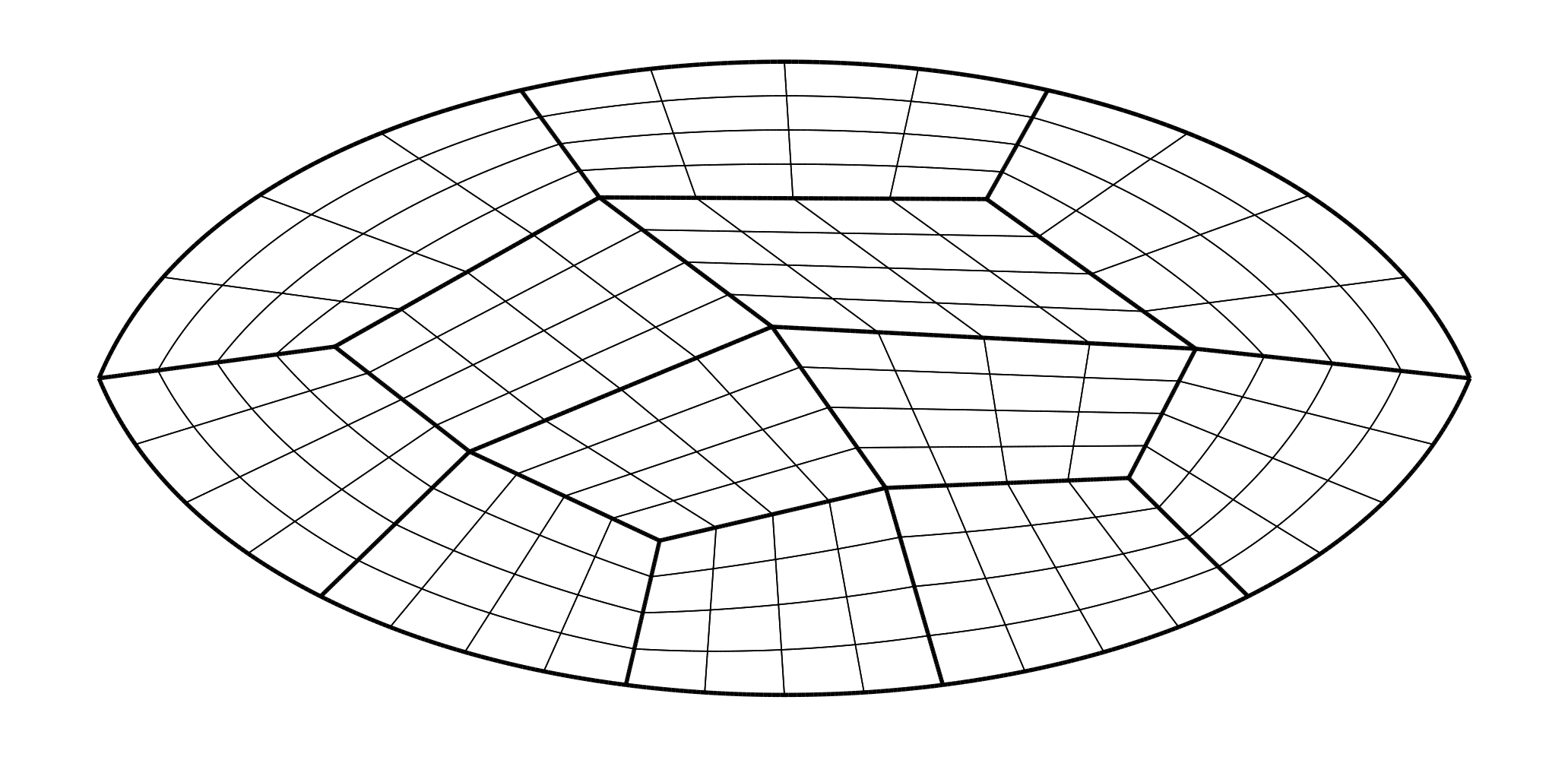}
    \end{subfigure}
\caption{Figure showing the process of controlmap creation starting from a PAG $P(8, 12)$. First, the PAG is utilised to create a covering of a regular eight-sided polygon via an approximately harmonic map from $\Om^{\square}(8)$ into $\hOm$. Then, the covering is triangulated and the triangulation is mapped onto an approximation $\hOm^{\mathbf{r}}_h$ of the lens domain $\hOm^{\mathbf{r}}$ using Floater's algorithm, with specified segmentation of the boundary. The inner vertices associated with the layout of $\hOm^{\mathbf{r}}$ follow from the correspondence of the two triangulations. Then, the inner vertex positions are optimised for quad layout area homogenisation. Finally, the map $\mathbf{r}: \hOm \rightarrow \hOm^{\mathbf{r}}$ follows from applying the bilinearly blended Coons' patch approach patchwise.}
\label{fig:adjgraph}
\end{figure}

\noindent Techniques for constructing all possible PAGs for different layouts $(N, N_p)$ are discussed in \cite{buchegger2017planar}. Given an $N$-sided PAG $P(N, N_p) = (V_p, E_p)$ with $N_p$ patches, we construct a concrete layout $T_i$ by introducing the set $\Om^{\square}(N_p) \subset \mathbb{R}^2 \times \{1, \ldots, N_p \}$ with $\Om^{\square}(N_p) := (\Om^{\square} \times \{1\}) \cup \ldots \cup (\Om^{\square} \times \{N_p\})$ and $\Om^{\square} := (0, 1)^2$. Here, each $\Om^\square \times \{i \}$ represents one of the $N_p$ patches. We construct the space
$$\mathcal{V}^{\square} := \{\phi \in \mathcal{Q}(\Om^{\square} \times \{i \}) \, \forall \, i \in \{1, \ldots, N_p\} \, \vert \, \phi \text{ assumes the same value on equivalent sides of } \partial \Om^{\square}(N_p) \},$$
where $\mathcal{Q}( \, \cdot \,)$ represents the canonical four-dimensional bilinear polynomial space over a quadrilateral two-dimensional manifold $Q \subset \mathbb{R}^3$ while the equivalence relation between the various sides of $\partial \Om^{\square}(N_p)$ follows from the edges of $P(N, N_p)$. Then, we approximate a harmonic map $\mathbf{x}^\square: \Om^{\square}(N_p) \rightarrow \hOm(N)$ by discretising the associated PDE over $\mathcal{V}^{\square}$ in the usual way. Again, the boundary correspondence $\mathbf{f}^{\square}: \gamma_D \rightarrow \partial \hOm$, with $\gamma_D \subseteq \partial \Om^{\square}(N_p)$ follows from the PAG $P(N, N_p)$. The interior vertices $v \in V_i$ of $T_i = (V_i, E_i, \mathcal{Q}_i)$ now follow from evaluating $\mathbf{x}^\square: \Om^{\square}(N_p) \rightarrow \hOm$ in the vertices of $\partial \Om^{\square}(N_p)$. Due to truncation, the resulting quad layout of $\hOm$ may be invalid in rare cases. If so, we repair the layout using an untangling routine. For an example of a PAG along with the associated template, we refer to Figure \ref{fig:adjgraph}.
Using this methodology, we constructed a template catalogue $\mathbb{T}$ comprised of $\approx 6 \times 10^4$ templates for various configurations of $(N, N_p)$, with $N \in \{4, 6, \ldots, 16\}$. We note that this is an \textit{offline} computation, which does not contribute to the computational costs of what follows. \\

\noindent Let $\mathbb{T} \ni T_i = (V_i, E_i, \mathcal{Q}_i)$ be a template. As discussed in Section \ref{sect:harmonic_maps}, various geometry types require introducing various control domains $\hOm^{\mathbf{r}}$. To construct the operator $\mathbf{r}: \hOm \rightarrow \hOm^{\mathbf{r}}$, we employ a three-step approach. First, we map $\hOm$ approximately harmonically onto a (convex) polygonal approximation $\hOm^{\mathbf{r}}_h$ of $\hOm^{\mathbf{r}}$ using Floater's algorithm (see Section \ref{sect:harmonic_maps}). This yields the piecewise linear map $\br_h: \hOm \rightarrow \hOm^\br_h$ between the triangulations of both domains. By $T^\br = (V^\br, E^\br, \mathcal{Q}^\br)$, we denote the template that results from replacing $V_i \ni v \rightarrow \br_h(v)$ in $T_i$. The template $T^\br$ represents a quadrangulation of a polygonal approximation $\hOm^\br_h$ of $\hOm^\br$. We define:
\begin{align}
\label{eq:controlmap_cross_product}
    \boldsymbol{\nu}(v^\br, q^\br) := (\mathbf{u}_{i_1} \times \mathbf{u}_{i_2})_3
\end{align}
wherein $\mathbf{u}_{i_1}$ and $\mathbf{u}_{i_2}$ represent the directions associated with the two edges $\{e_{i_1}^\br, e_{i_2}^\br \} \subset q^\br$ incident to $v^\br \in V^\br$ on $q^\br \in \mathcal{Q}^\br$ with their lengths corresponding to the distance between the edges' two end points. Furthermore, $\mathbf{a} \times \mathbf{b}$ denotes the cross product between vectors $\mathbf{a}$ and $\mathbf{b}$. A valid quadrangulation is one for which all $\boldsymbol{\nu}(v^\br, q^\br) > 0$. We optimise the quadrangulation by minimising the cost function
\begin{align}
\label{eq:ReLU_template_optimisation}
    C(\mathcal{T}^\br, \tau) := \sum \limits_{q^\br \in \mathcal{Q}^\br} \sum \limits_{v^\br \in \mathbb{V}(q^\br)} \operatorname{ReLU}(-\boldsymbol{\nu}(v^\br, q^\br) + \tau)
\end{align}
over all interior vertices $V^\br \ni v^\br \in V^\br_{\text{int}}$. Here, $\operatorname{ReLU}(\, \cdot \,)$ denotes the ReLU activation function $\operatorname{ReLU}(x) := \max \{0, x \}$, while $\tau \geq 0$ is an offset. We note that $C(\, \cdot \,, \tau)$ vanishes once all $\boldsymbol{\nu}(v^\br, q^\br)$ are above the threshold $\tau$ and the minimisation is terminated. Minimisation based on~\eqref{eq:ReLU_template_optimisation} is nonconvex and the cost function is nonsmooth. However, a gradient-based algorithm nevertheless converges after a small number of iterations as long as the initial guess represents a layout for which all $\boldsymbol{\nu}(v^\br, q^\br)$ are close to or exceed $\tau$. Starting on $\tau = 0$, its value is incremented in an outer loop until optimisation no longer yields a minimiser for which $C(\mathcal{T}^\br, \tau) = 0$. In practice, the optimisation terminates in a fraction of a second thanks to the low number of interior vertices. Reusing the previously computed map $\br_h: \hOm \rightarrow \hOm^\br_h$ to create an initial layout for optimisation based on~\eqref{eq:ReLU_template_optimisation}, in this way, hundreds of templates can be optimised in a matter of seconds. \\
Once a properly optimised template $T^\br$ has been chosen, it is converted into a controlmap $\br: \hOm \rightarrow \hOm^\br$ by replacing boundary edges by the corresponding sides of $\partial \hOm^\br$ while applying the bilinearly-blended Coons' patch approach to the resulting (curvilinear) quad cells. We note that thanks to the convexity of $\hOm^\br$, the curvilinear quads formed by replacing boundary edges remain convex, provided each $q^\br \in \mathcal{Q}^\br$ is convex. \\
In practice, the final step is only carried out for a single template which is selected from a large catalogue of optimised templates.\\
All required steps are summarised by Figure \ref{fig:adjgraph}.

\bibliographystyle{abbrv}
\bibliography{bibliography}


\section*{List of Symbols}
\begin{table}[H]
    \centering
    \begin{tabular}{l p{15cm}}
        \textbf{Section 1} & \\
        \hline
         Symbol & Property \\
         \hline 
         $H^s(\Omega), H^s(\Omega, \mathbb{R}^n)$ & (Vectorial) Sobolev space of order $s \geq 0$ \\
         $G(V, E, \mathcal{F}, \mathcal{T})$ & Graph comprised of vertices $V$, edges $E$, faces $\mathcal{F}$ and optionally templates $\mathcal{T}$\\
         $\partial E$ & Set of edges $e \in E$ located on the graph's boundary \\
         $\iota(\, \cdot \,)$ & {$\iota(e), \, e \in E$ returns the vertex indices $(i, j)$ to which $e$ is incident, $\iota(F), \, F \in \mathcal{F}$ returns a tuple of index pairs, one for each $e \in F$.} \\
         $U(e)$ & Returns $e := (i, j)$ if $e \in E$ and $-e := (j, i)$ if $-e \in E$ \\
         $w(\, \cdot \,)$ & Weight function returning point sets $w(e) = p = \{p_1, \ldots, p_{N} \}$, $w(-e)$ reverses the order of the point set $p \subset \mathbb{R}^2$ \\
         $\Om_i$ and $\Om^S_i$ & Dense polygon resulting from combining all $p = w(e)$ with $e \in F$, $F \in \mathcal{F}$ and its spline approximation\\
         $\mathbb{V}(\, \cdot \,)$ & $\mathbb{V}(F),\, F \in \mathcal{F}$ gives the $v \in V$ as they appear in $\iota(F) = \{(p, q), (q, r), \ldots \}$\\
         $\hOm(N)$ & $N$-sided regular polygon of radius one \\
         $w^S(\, \cdot \,)$ & Weight function assigning to each $\pm e \in E$ a spline $s \approx w(e)$ \\
         $\mathbf{f}_i: \partial \hOm_i \rightarrow \partial \Om^S_i$ & Spline-based boundary correspondence\\
         $\hat{t}_{-}{v}$, $\hat{t}_{+}{v}$ & Discrete unit tangent on either side of the vertex $v \in V$ \\
         $\hat{t}(v)$, $\hat{n}(v)$ & Discrete average unit tangent and normal at $v \in V$ \\
         $\angle(v)$ & Discrete interior angle at vertex $v \in \mathbb{V}(F)$ \\
         $\mu_{\angle} \ll 1$ & User-defined threshold that models a vertex $v \in \mathbb{V}(F)$ as a convex corner if $\angle(v) < \pi - \mu_\angle$\\
         $L(\, \cdot \,)$ & $L(\partial \Om)$ denotes the total circumference of $\Om$ while $L(e)$ denotes the length of an edge \\ 
         $\br: \hOm \rightarrow \hOm^\br$ & Control map from $\hOm$ into the control domain $\hOm^\br$ \\
         \hline
    \end{tabular}
\end{table}
\begin{table}[H]
    \centering
    \begin{tabular}{l p{15cm}}
        \textbf{Section 2} & \\
        \hline
         Symbol & Property \\
         \hline 
         $\varepsilon_\angle$ & User-defined threshold that flags a vertex $v \in F$ as concave if $\angle(v) \geq \pi + \varepsilon_\angle$ \\
         $\mathbb{V}^{\varepsilon_\angle}_{\text{conc}}(\, \cdot \,)$ & $\mathbb{V}^{\varepsilon_\angle}_{\text{conc}}(F)$ returns all $v \in \mathbb{V}(F)$ that have been flagged concave according to $\varepsilon_\angle$ \\
         $\mathcal{Q}( \, \cdot \,)$ & Operator providing a measure $\mathcal{Q}(C)$ for the straightness of a curve $C \in C^2([0, 1])$ \\ 
         $\mathcal{Q}^\mu(\, \cdot \,)$ & Operator favouring curves $C$ that connect two concave corners\\
         $T_i = (V_i, E_i \mathcal{Q}_i)$ & Template graph with vertices $V_i$, edges $E_i$ and quadrangular faces $\mathcal{Q}_i$\\
         $\mathbb{T}$ & Complete pre-computed catalogue of templates $T_i$ \\
         $\phi_i: \partial E_i \rightarrow F_i$ & Boundary correspondence between boundary edges $\partial E_i$ of $T_i \in \mathcal{T}$ and the face $F_i \in \mathcal{F}$ \\
         $\operatorname{val}(\, \cdot \,)$ & Operator returning the valence of a vertex $v \in V$ in $G = (V, E, \mathcal{F})$ \\
         $F_e$ & The set of faces $F \in \mathcal{F}$ with $\pm e \in F$ \\
         $F_{\mathcal{T}}$ & Subset of faces to which a template has been assigned \\
         $G^\square = (V^\square, E^\square, \mathcal{F}^\square)$ & Canonical template skeleton graph of $G$ \\
         $E_\mathcal{T}$ & Subset of edges associated with at least one template $T_i \in \mathcal{T}$ \\
         $L(e)$ & Length of the piecewise linear curve resulting from the points $p = w(e)$ \\
         $L^{\mu_{\mathcal{T}}, \mu_\partial}(\, \cdot \,)$ & Scaled length function, assigning larger values to edges $e \notin E_\mathcal{T}$ and $e \in \partial E$ \\
         $\varepsilon_L$ & Parameter $0 < \varepsilon_L \leq 1$ that marks an edge eligible for splitting if $L(e) \geq \varepsilon_L L_{\max}$ \\
         $\bx_h^\br: \hOm_h^\br \rightarrow \Om_h$ & Surrogate harmonic map between $\hOm_h^\br \approx \hOm^\br$ and $\Om_h \approx \Om^S$ computed using Floater's algorithm\\
         $\theta(\hat{v})$ & Preferred angle created by point sets incident to $\bx_h \circ \br(\hat{v})$\\
         \hline
    \end{tabular}
\end{table}

\begin{table}[H]
    \centering
    \begin{tabular}{l p{15cm}}
        \textbf{Section 3} & \\
        \hline
         Symbol & Property \\
         \hline 
         $\Xi^0$ & Base knotvector with specified number of interior knots\\
         $r_j := \|s(\xi_j) - p_j\|$ & $l^2$ mismatch between the spline fit and the $j$-th fitting point $p_j$ \\
         $\mu_{LS}$ & Threshold value that flags a knotspan for refinement if $r_j \geq \mu_{LS}$ \\
         $w^\Xi(\, \cdot \,)$ & Weight function assigning to $e \in E$ the knotvector $\Xi$ associated with $s = w^S(e)$ \\
         $w^\Xi_i(\, \cdot \,)$ & Weight function assigning a knotvector to each $\hat{e} \in \partial E_i$ of $\mathcal{T}_i = (V_i, E_i, \mathcal{Q}_i)$ \\
         $\mathcal{V}_{h, i}$ & Canonical spline space on $\hOm_i$ under the layout $T_i \in \mathcal{T}$ and the knotvectors assigned by $w_i^\Xi(\, \cdot \,)$\\
         $\mathcal{L}_\eta^\mu(\, \cdot \,, \, \cdot\, )$ & Semi-linear form used for computing an inversely harmonic map\\
         \hline
    \end{tabular}
\end{table}

\begin{table}[H]
    \centering
    \begin{tabular}{l p{15cm}}
        \textbf{Section 5} & \\
        \hline
         Symbol & Property \\
         \hline 
         $W(\, \cdot \,)$ & Winslow function\\
         $W_\varepsilon(\, \cdot \,)$ & Regularised Winslow function \\
         $\mathcal{R}_\varepsilon(\, \cdot \,)$ & Jacobian determinant regulariser \\
         $\mathcal{L}_\varepsilon^W$ & Regularised weak form discretisation semi-linear form \\
         \hline
    \end{tabular}
\end{table}

\end{document}